\documentclass[a4paper,leqno]{article}

\usepackage{amsmath}
\usepackage{amssymb}
\usepackage{amsfonts}
\usepackage{enumerate}
\usepackage{theorem}
\usepackage{esint}
\usepackage{bbm}

\theoremstyle{plain}
\newtheorem{teo}{Theorem}[section]
\newtheorem{lem}[teo]{Lemma}
\newtheorem{cor}[teo]{Corollary}
\newtheorem{prop}[teo]{Proposition}

{\theorembodyfont{\rmfamily}
\newtheorem{defin}[teo]{Definition}
\newtheorem{oss}[teo]{Remark}

}

\renewcommand{\eqref}[1]{\textnormal{(\ref{#1})}}

\numberwithin{equation}{section}

\newcommand{\cvd}{\hfill$\square$}
\newcommand{\proof}[1]{\noindent\textsc{Proof#1}}

\newcommand{\R}{\mathbb{R}}
\newcommand{\N}{\mathbb{N}}

\title{Full discretization and regularization for the Calder\'on problem}

\author{Alessandro Felisi\thanks{Dipartimento di Matematica, Università degli Studi di Genova, Italy.\newline E-mail:
\texttt{felisi@dima.unige.it}} \and\
Luca Rondi\thanks{Dipartimento di Matematica,
Universit\`a degli Studi di Pavia, Italy.\newline E-mail: \texttt{luca.rondi@unipv.it}} }

\date{%Draft of \today
}

\begin{document}

\maketitle

\setcounter{section}{0}
\setcounter{secnumdepth}{2}

\begin{abstract}
We consider the inverse conductivity problem with discontinuous conductivities.
We show in a rigorous way, by a convergence analysis, that one can construct a completely discrete minimization problem whose solution is a good approximation of a solution to the inverse problem. The minimization problem contains a regularization term which is given by a total variation penalization and is characterized by a regularization parameter. The discretization involves at the same time the boundary measurements, by the use of the complete electrode model, the unknown conductivity and the solution to the direct problem. The electrodes are characterized by a parameter related to their size, which in turn controls the number of electrodes to be used.
The discretization of the unknown and of the solution to the direct problem is characterized by another parameter related to the size of the mesh involved. In our analysis we
show how to  precisely choose the regularization, electrodes size and mesh size parameters with respect to the noise level in such a way that the solution to the discrete regularized problem is meaningful. In particular we obtain that the electrodes and mesh size parameters should decay polynomially with respect to the noise level.

\medskip

\noindent\textbf{AMS 2020 Mathematics Subject Classification} 35R30 (primary); 49J45 65N21 (secondary)

\medskip

\noindent \textbf{Keywords} inverse problems, regularization, discretization, total variation, finite elements.
\end{abstract}

\section{Introduction}

Let $\Omega\subset\mathbb{R}^N$, $N\geq 2$, be a bounded domain with Lipschitz boundary. Let $\sigma_0$ be a conductivity in $\Omega$. We assume that such a conductivity is unknown and we wish to determine and, possibly, reconstruct it by performing boundary measurements of current and voltage kind. In the ideal case of performing infinitely many measurements, this corresponds to the Calder\'on problem, \cite{Cal}, or inverse conductivity problem. This problem attracted a great attention for its numerous applications in various fields, such as nondestructive evaluation in engineering and medical imaging, just to mention a few. It is also an extremely interesting and challenging problem from a mathematical point of view and it can be seen as the prototype of several other significant inverse boundary value problems.

Namely, let $g\in L^2_{\ast}(\partial\Omega)=\{\tilde{g}\in L^2(\partial\Omega):\ \int_{\partial\Omega}\tilde{g}=0\}$ be the applied current density on $\partial\Omega$. Then the electrostatic potential $u$ is the solution to the following Neumann boundary value problem
\begin{equation}\label{Neupbmintro}
\left\{\begin{array}{ll}
-\mathrm{div}(\sigma_0\nabla u)=0 &\text{in }\Omega\\
\sigma_0\nabla u\cdot\nu=g&\text{on }\partial\Omega\\
\int_{\partial\Omega} u=0, &
\end{array}\right.
\end{equation}
$\nu$ being the outer normal.
We measure such a potential still at the boundary, that is, we measure $u|_{\partial\Omega}\in L^2_{\ast}(\partial\Omega)$. Therefore, in principle, one could measure the so-called \emph{Neumann-to-Dirichlet} map $\mathcal{N}(\sigma_0):L^2_{\ast}(\partial\Omega)\to L^2_{\ast}(\partial\Omega)$ such that for any $g\in L^2_{\ast}(\partial\Omega)$
$$\mathcal{N}(\sigma_0)[g]=u|_{\partial\Omega}\quad\text{where }u\text{ solves \eqref{Neupbmintro}}.$$
We have that $\mathcal{N}(\sigma_0)$ is a bounded and linear operator.
The inverse conductivity problem is to find $\sigma_0$ from (a partial knowledge of) $\mathcal{N}(\sigma_0)$. Since $\mathcal{N}(\sigma_0)$ comes
from measurements, the available data are, still in the best scenario, a perturbation of $\mathcal{N}(\sigma_0)$, which we call $\mathcal{N}_{\varepsilon}$, where $\varepsilon>0$ represents the noise level.

Even when uniqueness holds, reconstructing $\sigma_0$ from $\mathcal{N}(\sigma_0)$, or better $\mathcal{N}_{\varepsilon}$, is an extremely challenging task because this problem is ill-posed. For example, the simple least-squares approach of solving, for some class $\mathcal{M}$ of admissible conductivities,
$$\min\left\{\|\mathcal{N}(\sigma)-\mathcal{N}_{\varepsilon}\|^2:\ \sigma\in\mathcal{M}\right\}$$
might lead to serious issues, see for instance a review of the related instabilities in \cite{Ron16}.
Here and in what follows, the norm is that of bounded and linear operators belonging to $\mathcal{L}(L^2_{\ast}(\partial\Omega),L^2_{\ast}(\partial\Omega))$, which in short we call $L^2\text{-}L^2$ norm.

In order to recover stability, some a priori information on the unknown $\sigma_0$ or, correspondingly, some regularization of the minimization problem is needed. However, even with strong a priori assumptions on the unknown, the problem is still severely ill-posed, as the example in \cite{Man} clearly shows.

For simplicity we restrict ourselves to regularizations of Tikhonov type. 
A general introduction to their use in inverse problem is in \cite{Eng-et-al}. In the smooth case, the theory of convergence of Tikhonov regularized solutions for nonlinear operators, with convergence estimate, was developed in \cite{Eng-Kun-Neu}, see also \cite{Eng-et-al}. For the inverse conductivity problem, the Tikhonov regularization in the smooth case was treated in \cite{LMP}, see also \cite{Ji-Ma}

We are interested in the more challenging case of discontinuous conductivities. In the nonsmooth case, a reasonable regularization is that of a total variation penalization, or other $BV$-related variants. For its use in solving  ill-posed linear problems, see \cite{AV,CK} and \cite{Vas1,Vas2} and the references therein.
 
 For the inverse conductivity problem, where linearity is not available any more, we still consider the a priori assumption that the unknown conductivity $\sigma_0$ belongs to $BV(\Omega)$, that is, it has bounded variation. We call
$$|\sigma|_{BV(\Omega)}=TV(\sigma)=|D\sigma|(\Omega)$$
the \emph{total variation} of a conductivity $\sigma$. Therefore, let us assume that $|\sigma_0|_{BV}$ is finite. The corresponding regularization is given by a
total variation penalization, that is, by solving the regularized problem
\begin{equation}\label{regminintro}
\min\left\{\|\mathcal{N}(\sigma)-\mathcal{N}_{\varepsilon}\|^2+a|\sigma|_{BV(\Omega)}:\ \sigma\in\mathcal{M}\right\}
\end{equation}
where the regularization parameter $a=a(\varepsilon)>0$ should be suitably chosen with respect to the noise $\varepsilon$.

This kind of regularization has been proven to be effective with many different numerical methods, see \cite{Dob-San94,Chan-Tai,Chu-Chan-Tai}. Only later on, in \cite{Ron08}, the validity of this approach was rigorously proved by showing the following convergence result. Let $\sigma_{\varepsilon}$ be a solution to the regularized problem \eqref{regminintro}. Then, provided $a(\varepsilon)$ is suitably chosen, and up to a subsequence, $\sigma_{\varepsilon}$ converges in $L^1(\Omega)$ to a conductivity $\hat{\sigma}$ which belongs to the following set $\hat{S}$ of optimal solutions of the inverse problem. Let $S=\{\sigma\in\mathcal{M}:\ \mathcal{N}(\sigma)=\mathcal{N}(\sigma_0\}$ be the set of solutions of the inverse problems. Then
$$\hat{S}=\left\{\sigma\in S:\ |\sigma|_{BV(\Omega)}=\min_{\tilde{\sigma}\in S}|\tilde{\sigma}|_{BV(\Omega)}\right\}.$$
That is, $\hat{\sigma}\in \hat{S}$ if
\begin{multline}\label{optimalsol}
\mathcal{N}(\hat{\sigma})=\mathcal{N}(\sigma_0)\quad\text{and}\\
|\hat{\sigma}|_{BV(\Omega)}=\min\{|\sigma|_{BV(\Omega)}:\ \sigma\in\mathcal{M}\text{ such that }\mathcal{N}(\sigma)=\mathcal{N}(\sigma_0)\}.
\end{multline}
In other words, $\hat{\sigma}$ is a solution to the inverse problem with minimal total variation.
When uniqueness holds, that is, when $S=\{\sigma_0\}$, then we have convergence, without passing to subsequences, of $\sigma_{\varepsilon}$ to  $\sigma_0$, still in the $L^1(\Omega)$ norm.

Other $BV$ related penalizations have been used for the Calder\'on problem, such as the Mumford-Shah functional introduced in \cite{Mu-Sh}. The numerical part was developed in \cite{Ron-San},
the convergence analysis is again in \cite{Ron08}. Further developments can be found in \cite{JMP} and in \cite{Ron16}, where the results of \cite{Che} are presented.

We give a brief account of the main uniqueness results for this inverse problem. The first uniqueness results, for scalar conductivities in dimension $N\geq 3$, were proved in \cite{Koh-Vog84:1,Koh-Vog85}. They showed that the conductivity at the boundary is determined by boundary measurements and treated the analytic case. This regularity was lowered to $C^2$ in \cite{Syl-Uhl87}. The piecewise smooth case was treated in \cite{Isak88}.
More recently, the regularity has been reduced to $C^1$ or Lipschitz but close to a constant, in \cite{Hab-Tat}, and to Lipschitz, in \cite{Car-Rog}. The up to date result is given in \cite{Hab}, where conductivities with unbounded gradient are allowed
and uniqueness is shown for $W^{1,N}$ conductivities, at least for $N=3,4$.

The two-dimensional case, that is, when $N=2$, is slightly different. The first uniqueness result for smooth conductivities was proved in \cite{Nac}. By complex analytic techniques, the two dimensional case is now completely solved, since uniqueness holds for $L^{\infty}$ scalar conductivities,
\cite{Ast-Pai}.

The anisotropic case, that is, when the scalar conductivity is replaced by a symmetric conductivity tensor, has a natural obstruction to uniqueness. Namely, a change of variables that keeps fixed the boundary does not change the boundary data of the transformed equation. In dimension $2$, this is the only obstruction, as shown in \cite{Syl} in the smooth case and finally in \cite{Ast-Pai-Las} for the $L^{\infty}$ case.

Unfortunately, in $\R^N$ with $N\geq 3$, the required a priori assumptions to have uniqueness are still much stronger than assuming the unknown to be a function with bounded variation.

The convergence result in \cite{Ron08} does not take into account two major difficulties.

The first one is that the Neumann-to-Dirichlet map involves infinitely many measurements. First numerical approaches, of a variational type, with a finite number of measurements are in \cite{Yor,K-V,McK}. However,
 the measurements that can be really obtained in the experiments are those describes in
 \cite{Som e Che e Isa}, which are called \emph{experimental measurements} or \emph{Complete Electrode Model} (CEM). They are encoded in an $M\times M$ matrix
 $R(\sigma_0)$, called \emph{resistance matrix}, where $M$ is the number of electrodes used. A first numerical investigation of CEM may be found in \cite{L-R}. Then,
in a series of papers, \cite{Hyv1,Hyv2,Hyv-et-al}, it has been shown that, provided the electrodes are suitable chosen, the resistance matrix may be seen as a good approximation of the Neumann-to-Dirichlet map. Therefore, it is reasonable to replace the Neumann-to-Dirichlet map with a perturbed resistance matrix $R_{\varepsilon}$, where again $\varepsilon$ denotes the noisel level.
 
 The second major difficulty is the discretization, which can be seen from two different viewpoints.
 
 The first one, the most important, is the discretization of the unknown. In fact, discretizing the unknown $\sigma_0$ produces a discretization error that adds to the noise error on the measurements and it is very difficult to tackle, see for example \cite{Riv-Bar-Ob}. The main issue is in finding a good balance between the desire of a good resolution (which requires a finer discretization) and that of stability (which requires a coarser discretization). An answer to this issue for the inverse conductivity problem was first given in \cite{Ron16}, with a precise indication of how to choose both the regularization parameter and the mesh size of the discretization with respect to the noise level. However, in \cite{Ron16}, the full Neumann-to-Dirichlet map was used, rather than the Complete Electrore Model, and the next discretization issue was not taken into account.
 
 Such a second discretization issue concerns the direct problem. In most numerical methods of reconstruction, one needs to solve, maybe several times, the direct problem, which is a boundary value problem for an elliptic equation like \eqref{Neupbmintro}. In practice, a discretization of such a boundary value problem is used and it would be interesting to understand how this affects the convergence of the overall method. A preliminary analysis of this problem is in \cite{Ge-Ji-Lu}. They use the Complete Electrode Model and study the convergence of a discretized model, in the polyhedral case and
 in the curved case, which is the main contribution of the paper. What is missing in their analysis is the fact that both the regularization parameter and the electrodes are kept fixed and only the discretization is allowed to change.
 
The aim of this paper is to provide a full discretization of the inverse problem, with a suitable regularization, for which a convergence results such as in \cite{Ron08} still holds. The regularization used is the total variation penalization and the discretization involves the measurements, with the use of CEM, the unknown conductivity and the direct problem. In other words, we combine all previous approximations together and simultaneously to obtain a fully discretized and regularized minimum problem whose solution is a good approximation of an optimal solution to the inverse problem. We mention that
a preliminary analysis towards this goal can be found in \cite{Fel}, whose results are here considerably sharpened.

In order to be more precise, let us illustrate our approach. In this introduction we consider for simplicity just scalar conductivities, but all our results carry over to the anisotropic case of symmetric conductivity tensors.

Let us assume that $\sigma_0$ is the unknown conductivity in $\Omega$.
We apply suitable $M$ electrodes on $\partial\Omega$. These electrodes are characterized by a parameter $\delta>0$ which describes the size of the electrodes. If the measurements are noise-free, then we could measure the $M\times M$ resistance matrix $R^{\delta}_0=R^{\delta}(\sigma_0)$, as described in \cite{Som e Che e Isa}. In practice, only an approximation of the resistance matrix can be obtained, namely our available data is 
encoded in an $M\times M$ matrix $R^{\delta}_{\varepsilon}$ such that
$$\|R^{\delta}_{\varepsilon}-R^{\delta}_0\|_{M\times M}\leq M\varepsilon.$$
Here $\varepsilon>0$ is the noise level and the norm $\|\cdot\|_{M\times M}$ is the Euclidean norm in $\R^{M\times M}$, with the $M\times M$ matrix identified with a vector in this space. In this case, we assume that the measurement at each electrode is marred by an error of order $\varepsilon$. We
assume that such an error is independent on the size of the electrode and, more importantly, on the value of the voltage as well, that is, it is not a relative error. We believe this to be the worst case, however, see Remark~\ref{osserv1} for the simpler case of a relative error $\varepsilon$.
We note that the number of electrodes $M$ grows, as $\delta\to 0^+$, like $\delta^{-(N-1)}$, see \eqref{Lest}.

With the same electrodes, we consider a simplified version of the resistance matrix, which, for any conductivity $\sigma$,
we call $\hat{R}^{\delta}(\sigma)$. In this simplified version we neglect the contact impedance on the electrodes, thus we just employ the usual direct problem \eqref{Neupbmintro}. This allows us to simplify the numerical treatment of the direct problem which is a completely standard one.

We then discretize the problem as following. For a parameter $h>0$, which describes the size of the mesh, we discretize our domain $\Omega$ with a suitable triangulation and we associate to it its corresponding finite element space, which we call $X^h$. $X^h$ is a finite dimensional subspace of $H^1(\Omega)$. When restricting to $X^h$, for any conductivity $\sigma$ we can find a discretized version of $\mathcal{N}(\sigma)$, of $R^{\delta}(\sigma)$ and of $\hat{R}^{\delta}(\sigma)$, which we call $\mathcal{N}_h(\sigma)$,
$R^{\delta}_h(\sigma)$ and $\hat{R}^{\delta}_h(\sigma)$, respectively.
In order to make it a really discrete problem, we should guarantee that each electrode is the union of a finite number of elements of the corresponding triangulation of $\partial\Omega$, see Remark~\ref{discreteelectrorem}.

We then consider the following completely discretized and regularized minimization problem
\begin{equation}\label{minregdiscr}
\min\left\{\| \hat{R}^{\delta}_h(\sigma) - R^{\delta}_{\varepsilon}  \|^2+a|\sigma|_{BV(\Omega)}:\ \sigma\in X^h\cap\mathcal{M}\right\}.
\end{equation}
where $a=a(\varepsilon)>0$ is the regularization parameter.
We note that for any $M\times M$ matrix $R$, $\|R\|$ denotes the norm as a linear operator from $\R^M$ into itself.
Such a minimization problem is completely discrete and admits a solution which, dropping the dependence on $\varepsilon$, we call $\sigma^{\delta}_{h,a}$. Even if the problem is completely discrete, solving numerically \eqref{minregdiscr} could be still a daunting task and we will not discuss these numerical difficulties in this paper.
What we are interested in is the following approximation result: we would like to find parameters $a=a(\varepsilon)$, $\delta=\delta(\varepsilon)$ and $h=h(\varepsilon)$, depending on the noise level $\varepsilon$, such that, as $\varepsilon\to 0^+$ we have
$$\sigma^{\delta(\varepsilon)}_{h(\varepsilon),a(\varepsilon)} \to \hat{\sigma}$$ 
where $\hat{\sigma}$ is an optimal solution to the Calder\'on inverse problem, that is, 
$\hat{\sigma}\in \hat{S}$ as in \eqref{optimalsol}.

In \cite{Hyv1,Hyv2,Hyv-et-al}, it has been shown how the resistance matrices approximates, as $\delta\to 0^+$, the corresponding Neumann-to-Dirichlet maps, with respect to the $L^2\text{-}L^2$ norm. Namely, they introduced
a suitable projection $Q$ and a suitable extension operator $E$, depending on the electrodes thus on $\delta$, that allow to transform any (simplified or not) resistance matrix $R$ into an operator
$E\circ R\circ Q:L^2_{\ast}(\partial\Omega)\to L^2_{\ast}(\partial\Omega)$. More precisely, instead of $R$ we should use here the corresponding operator $\mathcal{R}$ defined in \eqref{mathcalR0}.
 Next they show that, under suitable geometric assumptions on the electrodes,
for any conductivity $\sigma$,
\begin{equation}\label{Hyvon}
\|E\circ R^{\delta}(\sigma)\circ Q-\mathcal{N}(\sigma) \|_{L^2\text{-}L^2}\to 0\quad\text{as }\delta\to 0^+.
\end{equation}
If $\sigma=\sigma_0$ and 
also $\varepsilon\to 0^+$, we conclude that the same property applies to $R^{\delta}_{\varepsilon}$, that is,
$$\|E\circ R^{\delta}_{\varepsilon}\circ Q-\mathcal{N}(\sigma_0) \|_{L^2\text{-}L^2}\to 0\quad\text{as }\delta,\,\varepsilon\to 0^+.$$

Our approach is inspired by $\Gamma$-convergence techniques. The most difficult part is to find a recovery sequence. Namely, we need to find $\sigma_h\in X^h\cap\mathcal{M}$ such that $\sigma_h\to \sigma_0$ as $h\to 0^+$ and such that
$$\|E\circ \hat{R}^{\delta}_h(\sigma_h)\circ Q-\mathcal{N}(\sigma_0)\|_{L^2\text{-}L^2}\to 0$$
in a suitable way.
We split this difficult problem into the following four terms
\begin{multline*}
\|E\circ \hat{R}^{\delta}_h(\sigma_h)\circ Q-\mathcal{N}(\sigma_0)\|
\leq  \|E\circ \hat{R}^{\delta}_h(\sigma_h)\circ Q-E\circ R^{\delta}_h(\sigma_h)\circ Q\|
\\
+\|E\circ R^{\delta}_h(\sigma_h)\circ Q-\mathcal{N}_{h}(\sigma_h)\|+
\|\mathcal{N}_{h}(\sigma_h)-\mathcal{N}(\sigma_h)\|
+
\|\mathcal{N}(\sigma_h)-\mathcal{N}(\sigma_0)\|.
\end{multline*}
The fourth term is the continuity of the 
 Neumann-to-Dirichlet map with respect to the coefficient. The third term is
the approximation of solutions of elliptic equations by their discretized counterparts using finite elements methods. The second term is the analogous of \eqref{Hyvon}, just at the discrete level. The first term corresponds to showing that the contact impedance does not play a significative contribution.

As we already pointed out,
one by one, the convergence of many of these terms have already been studied. The difficulty here is that we need very precise and quantitative convergence estimates in terms of $\delta$ and $h$. For example, the first two terms need to converge as $\delta\to 0^+$ in a way that is independent from $h$.

In our main result, Theorem~\ref{mainteo}, we show that if we choose $a=a(\varepsilon)$, $h=h(\varepsilon)$ and $\delta=\delta(\varepsilon)$, with respect to $\varepsilon$, as in \eqref{hdef-deltadef} and \eqref{coeffdef}, then the solutions to the corresponding regularized and discretized minimization problems converge, up to subsequences, to an optimal solution to our inverse problem. Therefore our result rigorously shows that, if we choose in the right way the regularization, the electrodes and the discretization, we end up with a fully discrete problem whose solution is, provided the noise error is small enough, a good approximation of an optimal solution to the inverse problem.

The plan of the paper is the following. In Section~\ref{sec2}, we introduce the notation and present some preliminary result. In particular, in Subsection~\ref{Gammasec} we recall the definition and main properties of $\Gamma$-convergence.
After introducing the notion of conductivity tensors, Subsection~\ref{tensorssubs},
we study 
several properties of Lipschitz domains and of their corresponding Sobolev spaces, Subsection~\ref{Sobolevsec}.
The proof of most of these results is sketched in the Appendix.
In Subsection~\ref{finiteelsec} we consider the discretization of the domain and recall the definition of the corresponding finite element space.
Then, Subsection~\ref{continuumsec}, we describe the continuum direct problem and we state its continuity with respect to the coefficient of the equation, Theorem~\ref{pcontinuity}, which allows us to deal with the fourth term above.
In Section~\ref{secCEM}, we review the experimental measurements and we present two different results. First, Subsection~\ref{L2L2ressubs}, we show that the experimental measurements are controlled by the Neumann-to-Dirichlet map, by slightly improving an analogous result of \cite{Ron15}. Then, Subsection~\ref{approxNtoDsubs}, we show, following
\cite{Hyv1,Hyv2,Hyv-et-al}, that the experimental measurements are a good approximation of the Neumann-to-Dirichlet map, which is the key to treat the second term above. The main novelty here is a careful handling of all the constants involved and the fact that the estimates are performed also at the discrete level by replacing the $H^1(\Omega)$ space with any subspace $X$ containing constants. Finally, Subsection~\ref{simplsec}, we introduce the simplified resistance matrix and estimate its difference with the resistance matrix, that is, we study the role of the contact impedance at the electrodes. This allows us to treat the first term above. In Section~\ref{discrsec} we deal with the discretization error, that is, with the third term above.
Following \cite{Ron16}, we approximate, in a suitable way, the unknown conductivity by piecewise linear ones, Proposition~\ref{Adiscrprop}. The study of the discretization error for the solution to the direct problem requires an even more delicate analysis and it is completely new in several aspects. This is carried over in Proposition~\ref{vdiscrprop} and it is applied to the Neumann-to-Dirichlet map in Corollary~\ref{Neumanncor}. In Section~\ref{mainsec},
we finally state and prove our main approximation result, Theorem~\ref{mainteo}. In Section~\ref{conclsec}, some final remarks and perspectives are presented. In the Appendix, we prove most of the results of Subsection~\ref{Sobolevsec} as well as Theorem~\ref{Meyers}.

\bigskip

\noindent
\textbf{Acknowledgement}\\
Luca Rondi acknowledges support by GNAMPA, INdAM.

\bigskip

\noindent
\textbf{Conflict of interest statement}\\
The authors have no conflict of interest to declare.

\section{Preliminaries}\label{sec2}

The integer $N\geq 2$ denotes the space dimension and we recall that we usually drop the dependence of any constant on $N$. 
For any Borel set $E\subset\mathbb{R}^N$, we denote with
$|E|$  its Lebesgue measure, whereas $\mathcal{H}^{N-1}(E)$ denotes its $(N-1)$-dimensional Hausdorff measure.
For any $x\in\mathbb{R}^N$ and any $s>0$, $B_s(x)$ denotes the open ball with center $x$ and radius $s$. Usually $B_s$ stands for $B_s(0)$.
For any $E\subset\mathbb{R}^N$, we call $B_s(E)=\bigcup_{x\in E}B_s(x)$.

For any two Banach spaces $B$, $B_1$, $\mathcal{L}(B,B_1)$
denotes the Banach space of bounded linear operators from $B$ to $B_1$ with the usual operator norm. The dual of $B$, that is,  
$\mathcal{L}(B,\R)$, is denoted by $B'$. For $x'\in B'$  and $x\in B$, we denote 
with $\langle x',x \rangle_{B',B}$ the usual duality, that is, $\langle x',x \rangle_{B',B}=x'[x]$. If there is no risk of confusion, we omit the subscript $B',B$.
We call $\mathcal{L}(B):=\mathcal{L}(B,B)$ and we denote its identity operator by $\mathbbm{1}$.
If $B\subset B_1$, we still call $\mathbbm{1}$ the natural immersion. We say that the immersion is continuous if $\mathbbm{1}$ is a bounded linear operator from $B$ into $B_1$, namely for some constant $C$ we have
$$\|x\|_{B_1}\leq C\|x\|_B\quad\text{for any }x\in B,$$
and that it is compact if $\mathbbm{1}$ is a compact linear operator from $B$ into $B_1$.
We note that $B\subset B_1$ with continuous (respectively compact) immersion implies that $B_1'\subset B'$ with continuous (respectively compact) immersion. Moreover, for the same constant $C$ we have
$$\|x'\|_{B'}\leq C\|x'\|_{B_1'}\quad\text{for any }x'\in B_1'.$$
Finally, if $B\subset B_1$ with continuous immersion, for some constant $C$, and $\tilde{B}\subset \tilde{B}_1$ with continuous immersion, for some constant $\tilde{C}$, then $\mathcal{L}(B_1,\tilde{B})   \subset \mathcal{L}(B,\tilde{B}_1)$ with continuous immersion, namely
$$\|L\|_{\mathcal{L}(B,\tilde{B}_1)}\leq C\tilde{C} \|L\|_{\mathcal{L}(B_1,\tilde{B})}\quad\text{for any }L\in \mathcal{L}(B_1,\tilde{B}).$$

For any $p$, $1\leq p\leq +\infty$, we denote with $p'$ its conjugate exponent, that is $1/p+1/p'=1$.

For any $M\in\N$,
we call $\mathbb{M}^{M\times M}(\mathbb{R})$ the space of real valued $M\times M$ matrices and
$\mathbb{M}_{sym}^{M\times M}(\mathbb{R})$ the subspace of real valued symmetric $M\times M$ matrices.
For any $A\in \mathbb{M}^{M\times M}(\mathbb{R})$, $\|A\|$ denotes its norm as a linear operator from $\R^M$ into itself, whereas
$\|A\|_{M\times M}$ denotes the Euclidean norm in $\R^{M\times M}$ where $A$ is seen as a vector in this space. We note that
\begin{equation}\label{matrixnorms}
\|A\|\leq \|A\|_{M\times M}   \leq  \sqrt{M}\|A\|\quad\text{for any }A\in \mathbb{M}^{M\times M}(\mathbb{R}).
\end{equation}
 
 For any measure $m$ and any measurable set $E$ such that $0<m(E)<+\infty$,
 $\displaystyle{\fint_{E}f}$ denotes the mean value of a function $f$ over the set $E$, that is,
$$\fint_{E}f=\frac{1}{m(E)}\int_Ef(x)dm(x).$$
Unless it is not clear from the context, we do not specify the measure with respect to which the mean is taken.
For any $1\leq p\leq +\infty$,
we call
$$L^p_{\ast}(E)=\left\{f\in L^p(E):\ \fint_{E}f=0\right\}.$$

We say that a subset of $\R^N$ is a \emph{domain} if it is open and connected.

\subsection{$\Gamma$-convergence}\label{Gammasec}

Our convergence result is essentially based on $\Gamma$-convergence. Its applications to the regularization of inverse problems, and in particular of Calder\'on problem, goes back to \cite{Ron08}.
We recall the definition and basic properties of
$\Gamma$-convergence, see \cite{DM} for a more detailed introduction.

Let $(X,d)$ be a metric space. Then a sequence
$F_n:X\to [-\infty,+\infty]$, $n\in\mathbb{N}$, $\Gamma$-converges as $n\to+\infty$
to a function $F:X\to [-\infty,+\infty]$ if for every $x\in X$ we have
\begin{align}\label{liminf}
&\text{for every sequence $\{x_n\}_{n\in\mathbb{N}}$ converging to $x$ we have}\\
&\hspace{2cm}F(x)\leq \liminf_n F_n(x_n);\nonumber\\
\label{limsup}
& \text{there exists a sequence $\{x_n\}_{n\in\mathbb{N}}$ converging to $x$ such that}\\
&\hspace{2cm}F(x)=\lim_n F_n(x_n).\nonumber
\end{align}
The function $F$ will be called the $\Gamma$-limit of the sequence $\{F_n\}_{n\in\mathbb{N}}$ as
$n\to+\infty$ with respect to the metric $d$ and we denote it by
$F=\Gamma\textrm{-}\!\lim_n F_n$.
We recall that condition \eqref{liminf} above is usually called the $\Gamma$-liminf inequality, whereas condition \eqref{limsup} is usually referred to as the existence of a recovery sequence.

We say that the functionals $F_n$, $n\in\mathbb{N}$, are \emph{equicoercive}
if there exists a compact set
$K\subset X$ such that $\inf_K F_n=\inf_X F_n$ for any $n\in\mathbb{N}$.

The Fundamental Theorem of
$\Gamma$-convergence is the following.

\begin{teo}\label{fundthm}
Let $(X,d)$ be a metric space and let $F_n:X\to [-\infty,+\infty]$,
$n\in\mathbb{N}$, 
be a sequence of functions defined on $X$. If the functionals $F_n$, $n\in\mathbb{N}$, are equicoercive
and
$F=\Gamma\textrm{-}\!\lim_n F_n$, then $F$ admits a minimum over $X$ and we
have
$$\min_X F=\lim_n\inf_X F_n.$$
Furthermore, if $\{x_n\}_{n\in\mathbb{N}}$ is a sequence of points in $X$ which
converges to a point $x\in X$ and
satisfies $\lim_n F_n(x_n)=\lim_n\inf_X F_n$, then $x$ is a minimum
point
for $F$.
\end{teo}

The definition of $\Gamma$-convergence may be extended in a natural way
to families depending on a continuous parameter.
The family of functions $F_{\varepsilon}$, defined for every
$\varepsilon>0$, $\Gamma$-converges to a function $F$ as
$\varepsilon\to 0^+$ if for every sequence $\{\varepsilon_n\}_{n\in\mathbb{N}}$ of positive numbers 
converging to $0$ as $n\to+\infty$, we have $F=\Gamma\textrm{-}\!\lim_n F_{\varepsilon_n}$.

\subsection{Conductivity tensors}\label{tensorssubs}

Let $\Omega\subset\R^N$ be a domain.
We say that $A$ is a \emph{symmetric conductivity tensor} in $\Omega$ if $A\in L^{\infty}(\Omega,\mathbb{M}_{sym}^{N\times N}(\mathbb{R}))$ and, for some constants
$0<\lambda_0\leq\lambda_1$,
$$\lambda_0\|\xi\|^2\leq A(x)\xi\cdot\xi\leq \lambda_1\|\xi\|^2\quad\text{ for any }\xi\in\mathbb{R}^N\text{and for a.e. }x\in\Omega.$$
Such a property is in short written as
\begin{equation}\label{ell}
\lambda_0 I_N\leq A(x)\leq \lambda_1 I_N\quad\text{for a.e. }x\in\Omega,
\end{equation}
where $I_N$ denotes the $N\times N$ identity matrix.

We say that a conductivity tensor $\sigma$ is a \emph{scalar conductivity}
if $A=\sigma I_N$ where $\sigma\in L^{\infty}(\Omega)$ and satisfies, for some constants
$0<\lambda_0\leq\lambda_1$,
\begin{equation}\label{ell-scalar}
\lambda_0\leq\sigma(x)\leq \lambda_1\quad\text{for a.e. }x\in\Omega.
\end{equation}
In this case, we often identify $A$ with the scalar function $\sigma$.

For constants $0<\lambda_0\leq \lambda_1$, we call
$$\mathcal{M}(\lambda_0,\lambda_1)=\{A:\ A\text{ is a conductivity tensor satisfying }\eqref{ell}\}$$
and
\begin{multline*}
\mathcal{M}_{scal}(\lambda_0,\lambda_1)=\{A\in\mathcal{M}(\lambda_0,\lambda_1):\ A \text{ is a scalar conductivity}\}\\=\{\sigma\in L^{\infty}(\Omega):\ \sigma\text{ satisfies }\eqref{ell-scalar}\}.
\end{multline*}

Since $\mathcal{M}(\lambda_0,\lambda_1)\subset L^{\infty}(\Omega,\mathbb{M}_{sym}^{N\times N}(\mathbb{R}))$, we measure
the distance between any two conductivity tensors $A_1$ and $A_2$ in $\Omega$ with an $L^p$ metric, for any $p$, $1\leq p\leq +\infty$,
as follows
$$\|A_1-A_2\|_{L^p(\Omega,\mathbb{M})}=\|(\|A_1-A_2\|)\|_{L^p(\Omega)},$$
where $\|A_1-A_2\|(x)=\|A_1(x)-A_2(x)\|$ for any $x\in \Omega$.
Clearly, for scalar conductivities
$A_i=\sigma_iI_N$, $i=1,2$, this reduces to the usual $L^p(\Omega)$ norm of $\sigma_1-\sigma_2$. 
Let us note that for any of these $L^p$ metrics, any of the classes $\mathcal{M}(\lambda_0,\lambda_1)$ and
$\mathcal{M}_{scal}(\lambda_0,\lambda_1)$ defined above
is a complete metric space.
 
 Let us here note that, analogously, for any $F\in L^p(\Omega,\mathbb{R}^N)$, $1\leq p\leq +\infty$, we set
 $$\|F\|_{L^p(\Omega,\R^N)}=\|(\|F\|)\|_{L^p(\Omega)},$$
 where $\|F\|(x)=\|F(x)\|$ for any $x\in \Omega$.
 
\subsection{Properties of Lipschitz domains and of Sobolev spaces}\label{Sobolevsec}
We say that a domain $\Omega\subset\mathbb{R}^N$
has a \emph{Lipschitz boundary} if for any $x\in\partial\Omega$ there exist a neighbourhood $U_x$  of $x$ and a
Lipschitz function
$\varphi:\mathbb{R}^{N-1}\to \mathbb{R}$ such that, up to a rigid change of coordinates, we have
$$\Omega\cap U_x=\{y=(y_1,\ldots,y_{N-1},y_N)\in U_x:\ y_N<\varphi(y_1,\ldots,y_{N-1})\}.$$
We usually denote with $\nu$ the exterior unit normal to $\partial\Omega$.

We say that a domain $\Omega\subset\mathbb{R}^N$
 belongs to the class $\mathcal{A}(r,L,R)$ if 
$\Omega\subset B_R$ and its boundary is Lipschitz with constants $r$ and $L$ in the following sense. For any $x\in\partial\Omega$
we can choose in the previous definition $U_x=B_r(x)$ and $\varphi$ with Lipschitz constant bounded by $L$.

\begin{oss}\label{oss1} For any bounded domain $\Omega$ with Lipschitz boundary, there exist constants $r$, $L$ and $R$ such that $\Omega\in \mathcal{A}(r,L,R)$.
\end{oss}

Let $\Omega\subset \R^N$ be a domain belonging to the class $\mathcal{A}(r,L,R)$.
First of all, we note that there exist positive constants $\tilde{c}_1$ and $\tilde{c}_2$, depending on $r$, $L$ and $R$ only, such that
\begin{equation}\label{measure}
\tilde{c}_1\leq |\Omega|,\,\mathcal{H}^{N-1}(\partial \Omega)\leq \tilde{c}_2.
\end{equation}

We say that $u\in L^1(\Omega)$ is a \emph{function of bounded variation} if its distributional derivative $Du$ is a vector valued Radon measure with finite total variation. We call $BV(\Omega)=\{u\in L^1(\Omega):\ u\text{ is of bounded variation}\}$, which is a Banach space if equipped with the norm
$$\|u\|_{BV(\Omega)}=\|u\|_{L^1(\Omega)}+|u|_{BV(\Omega)}$$
where the $BV$-seminorm corresponds to the \emph{total variation} of $u$, that is,
$$|u|_{BV(\Omega)}=TV(u)=|Du|(\Omega).$$
It is important to note that $|\cdot|_{BV(\Omega)}$ is lower semicontinuous with respect to convergence in $L^1(\Omega)$ and that, provided $\Omega$ is a bounded Lipschitz domain, the immersion of $BV(\Omega)$ into $L^1(\Omega)$ is compact,
see \cite[Theorem~3.23]{Amb-Fus-Pal} for instance.
For more information on $BV$ functions, we refer for instance to \cite{Giu,Amb-Fus-Pal}.

In order to deal with conductivity tensors with $BV$ regularity, we say that $A\in L^1(\Omega,\mathbb{M}^{N\times N}(\R))$ belongs to 
$BV(\Omega,\mathbb{M}^{N\times N}(\R))$ if each element of the matrix, $A_{ij}$ for any $i,\,j=1,\ldots,N$, belongs to $BV(\Omega)$. 
$BV(\Omega,\mathbb{M}^{N\times N}(\R))$ is a Banach space if endowed with the norm
$$\|A\|_{BV(\Omega,\mathbb{M})}=\|A\|_{L^1(\Omega,\mathbb{M})}+|A|_{BV(\Omega,\mathbb{M})}$$
where the $BV$-seminorm is defined as follow
$$|A|_{BV(\Omega,\mathbb{M})}=\|TV(A)\|.$$
where $TV(A)$ is the matrix whose elements are $(TV(A))_{ij}=TV(A_{ij})$, for any $i,\,j=1,\ldots,N$. 
We note that when $A=\sigma I_N$, then $A\in BV(\Omega,\mathbb{M}^{N\times N}(\R))$ if and only if $\sigma\in BV(\Omega)$ and
$$\|A\|_{BV(\Omega,\mathbb{M})}=\|\sigma\|_{BV(\Omega)}.$$
As before, one can show that
$|\cdot|_{BV(\Omega,\mathbb{M})}$ is lower semicontinuous with respect to convergence in $L^1(\Omega,\mathbb{M}^{N\times N}(\R))$ and that, provided $\Omega$ is a bounded Lipschitz domain, the immersion of $BV(\Omega,\mathbb{M}^{N\times N}(\R))$ into $L^1(\Omega,\mathbb{M}^{N\times N}(\R))$ is compact.

For the remaining part of this section we fix $p$, $1\leq p< +\infty$. We consider the usual Sobolev space $W^{1,p}(\Omega)$ which we endow with the norm
$$\|u\|_{W^{1,p}(\Omega)}:=\left(\|u\|_{L^p(\Omega)}^p+\|\nabla u\|_{L^p(\Omega,\R^N)}^p\right)^{1/p}\quad\text{for any }u\in W^{1,p}(\Omega).$$
Recall that $W^{1,1}(\Omega)\subset BV(\Omega)$ and that for any $u\in W^{1,1}(\Omega)$ we have that
$\|u\|_{W^{1,1}(\Omega)}=\|u\|_{BV(\Omega)}$. For more information on Sobolev spaces, we refer for instance to \cite{Ada,Leo}.

We let $W^{1,p}_0(\Omega)$ be the closure in $W^{1,p}(\Omega)$ of $C^{\infty}_0(\Omega)$. For these functions, the \emph{Poincar\'e inequality} holds, namely
there exists a constant $c_P$, depending on $p$ and $R$ only, such that
\begin{equation}\label{Poin0}
\|u\|_{L^p(\Omega)}\leq c_P\| \nabla u\|_{L^p(\Omega,\mathbb{R}^N)}\quad\text{for any }u\in W^{1,p}_0(\Omega).
\end{equation}
We define on $W^{1,p}_0(\Omega)$ the norm
$$\|u\|_{W^{1,p}_0(\Omega)}:=\|\nabla u\|_{L^p(\Omega,\R^N)}\quad\text{for any }u\in W^{1,p}_0(\Omega).$$
Such a norm, on $W^{1,p}_0(\Omega)$, is topologically equivalent to the usual $W^{1,p}(\Omega)$ norm, since
$$\|u\|_{W^{1,p}_0(\Omega)}\leq (1+c_P^p)^{1/p}\|\nabla u\|_{L^p(\Omega,\R^N)}\quad\text{for any }u\in W^{1,p}_0(\Omega).$$

For $p>1$, we consider the Besov space
$B^{1-1/p,p}(\partial\Omega)$ which is endowed with the following norm.
For any $\varphi\in B^{1-1/p,p}(\partial \Omega)$ we have
$$\|\varphi\|_{B^{1-1/p,p}(\partial \Omega)}=\left(\|\varphi\|_{L^p(\partial\Omega)}^p+| \varphi |_{B^{1-1/p,p}(\partial\Omega)}^p\right)^{1/p}$$
where the $B^{1-1/p,p}$ seminorm is defined as
\begin{equation}\label{Besovseminorm}
| \varphi |_{B^{1-1/p,p}(\partial\Omega)}=\left(\int_{\partial\Omega}\left(\int_{\partial\Omega}\frac{|\varphi(x)-\varphi(y)|^p}{|x-y|^{N+p-2}}d\mathcal{H}^{N-1}(x)\right)d\mathcal{H}^{N-1}(y)\right)^{1/p}.
\end{equation}

The following continuous and compact immersions hold, see for instance \cite[Theorem~7.57]{Ada}. We consider the following cases
\begin{equation}\label{pqcases}
\left\{
\begin{array}{ll}
\text{a)} &\text{$1<p<N$ and $1\leq q\leq \dfrac{(N-1)p}{N-p}$;}\\
\text{b)} &\text{
$p=N$ and $1\leq q<+\infty$;}\\
\text{c)} &\text{$p>N$ and $1\leq q\leq+\infty$.}
\end{array}\right.
\end{equation}
Then there exists a constant $c_{p,q}$, depending on $p$, $q$, $r$, $L$ and $R$ only, such that
\begin{equation}\label{Besovimmersion}
\|\varphi\|_{L^q(\partial\Omega)}\leq c_{p,q}\|\varphi\|_{B^{1-1/p,p}(\partial\Omega)}\quad\text{for any }\varphi\in B^{1-1/p,p}(\partial\Omega).
\end{equation}
We also note that, in case a), immersion is compact for any $1\leq q<\dfrac{(N-1)p}{N-p}$, whereas, in cases b) and c),
immersion is compact for any $1\leq q<+\infty$. In particular, for any $1<p<+\infty$, $B^{1-1/p,p}(\partial\Omega)$ is compactly immersed in $L^p(\partial\Omega)$.
The only novelty in this result is the fact that the constant $c_{p,q}$ in \eqref{Besovimmersion} depends on $\Omega$ only through the geometric constants $r$, $L$ and $R$. This fact can be deduced by the estimates in \cite[Theorem~7.57]{Ada}, via the geometric construction described
 in the Appendix and the invariance of the norms involved through bi-Lipschitz transformations.

We call $TW^{1,p}(\Omega)$ the space of traces of $W^{1,p}(\Omega)$ functions on $\partial\Omega$.
It is well known that 
$TW^{1,p}(\Omega)\subset L^p(\partial\Omega)$.
We call $W^{1,p}_{\ast}(\Omega)$ the set of functions in $W^{1,p}(\Omega)$ whose trace has zero mean on $\partial\Omega$, that is, whose trace belongs to 
$$TW^{1,p}_{\ast}(\Omega) :=TW^{1,p}(\Omega)\cap L^p_{\ast}(\partial\Omega).$$
With the help of the Poincar\'e inequality \eqref{Poin1}, which we discuss below, for any $1\leq p<+\infty$ we define on $W^{1,p}_{\ast}(\Omega)$ the norm
$$\|u\|_{W^{1,p}_{\ast}(\Omega)}:=\|\nabla u\|_{L^p(\Omega,\R^N)}\quad\text{for any }u\in W^{1,p}_{\ast}(\Omega).$$
Such a norm is, on $W^{1,p}_{\ast}(\Omega)$, topologically equivalent to the usual $W^{1,p}(\Omega)$ norm, since
$$\|u\|_{W^{1,p}_
{\ast}(\Omega)}\leq (1+C_P^p)^{1/p}\|\nabla u\|_{L^p(\Omega,\R^N)}\quad\text{for any }u\in W^{1,p}_{\ast}(\Omega).$$
We finally note that, for any $u\in W^{1,p}(\Omega)$, there exist $u_{\ast}\in W_{\ast}^{1,p}(\Omega)$ and a constant $c$ such that
$u=u_{\ast}+c$. Such a decomposition is actually unique and $c=\displaystyle{\fint_{\partial\Omega}u}$.

Any trace space will be endowed with the quotient norm, namely
$$\|\varphi\|_{TW^{1,p}(\Omega)}=\inf\left\{\|u\|_{W^{1,p}(\Omega)}:\ u\in W^{1,p}(\Omega)\text{ and } u|_{\partial\Omega}=\varphi\right\}$$
and
$$\|\varphi\|_{TW_{\ast}^{1,p}(\Omega)}=\inf\left\{\|u\|_{W^{1,p}_{\ast}(\Omega)}:\ u\in W^{1,p}_{\ast}(\Omega)\text{ and } u|_{\partial\Omega}=\varphi\right\}.$$

\begin{oss}
For $1<p<+\infty$, in both cases such an infimum is actually a minimum.
\end{oss}

We state the following \emph{trace inequality} and a characterization of the norm of traces, a sketch of the proof may be found in the Appendix.

\begin{teo}\label{traceineqteo}
Let $\Omega\in \mathcal{A}(r,L,R)$.
For any $1\leq p<+\infty$, there exists a constant $C_T$, depending on $p$, $r$, $L$ and $R$ only, such that the following holds.

If $p>1$, then
\begin{equation}\label{traceineq}
\|u\|_{B^{1-1/p,p}(\partial\Omega)}\leq C_T \|u\|_{W^{1,p}(\Omega)}\quad\text{for every }u\in W^{1,p}(\Omega).
\end{equation}
On the other hand, for any $\varphi\in B^{1-1/p,p}(\partial\Omega)$, there exists a function $u\in W^{1,p}(\Omega)$ such that
$u|_{\partial\Omega}=\varphi$ and
\begin{equation}\label{reversetraceineq}
\|u\|_{W^{1,p}(\Omega)} \leq C_T \|\varphi\|_{B^{1-1/p,p}(\partial\Omega)}.
\end{equation}

When $p=1$, we have
\begin{equation}\label{traceineq2}
\|u\|_{L^1(\partial\Omega)}\leq C_T \|u\|_{BV(\Omega)}\quad\text{for every }u\in BV(\Omega).
\end{equation}
Clearly, this last formula still holds if we replace $BV(\Omega)$ with $W^{1,1}(\Omega)$.
On the other hand, for any $\varphi\in L^1(\partial\Omega)$, there exists a function $u\in W^{1,1}(\Omega)$ such that
$u|_{\partial\Omega}=\varphi$ and
\begin{equation}\label{reversetraceineq2}
\|u\|_{W^{1,1}(\Omega)} \leq C_T \|\varphi\|_{L^1(\partial\Omega)}.
\end{equation}
\end{teo}

\begin{oss} We conclude that, for $p=1$, $TW^{1,1}(\Omega)$ and $L^1(\partial\Omega)$ and, for $1<p<+\infty$,
$TW^{1,p}(\Omega)$ and $B^{1-1/p,p}(\partial\Omega)$ coincide and they have topologically equivalent norms, since all immersions are continuous with constant $C_T$.
\end{oss}

It is also well-known that any bounded Lipschitz domain is an extension domain. We just show that the norm of the extension operator depend on $p$, $r$, $L$ and $R$ only. Also the proof of this result is sketched in the Appendix.

\begin{teo}\label{extension}
Let $\Omega\in \mathcal{A}(r,L,R)$. Then there exists a linear map $S:BV(\Omega)\to BV(\R^N)$ such that,
for any $u\in BV(\Omega)$, $Su\equiv 0$ outside $B_{R+1}$ and $Su=u$ in $\Omega$. Furthermore, we have that
$|D(Su)|(\partial\Omega)=0$ and
there exists a constant $C_E$, depending on $r$, $L$ and $R$ only, such that, for any $u\in BV(\Omega)$,
$$\|Su\|_{BV(\R^N)}\leq C_E\| u\|_{BV(\Omega)}.$$

Furthermore, for any $1\leq p<+\infty$, we have that if $u\in W^{1,p}(\Omega)$, then $Su\in W^{1,p}(\Omega)$ as well and 
there exists a constant $C_E$, depending on $r$, $L$, $R$ and $p$ only, such that
$$\|Su\|_{W^{1,p}(\R^N)}\leq C_E\| u\|_{W^{1,p}(\Omega)}\quad\text{for any }u\in W^{1,p}(\Omega).$$
\end{teo}

The following versions of the \emph{Poincar\'e inequality} are useful.

\begin{teo}\label{Pointeo}
For any $1\leq p<+\infty$, there exist constants $C_P$, $\tilde{C}_P$ and $\hat{C}_P$, depending on $p$, $r$, $L$ and $R$ only, such that, for any $u\in W^{1,p}(\Omega)$, we have
\begin{equation}\label{Poin1}
\left\|u-\fint_{\partial \Omega}u\right\|_{L^p(\Omega)}\leq C_P\| \nabla u\|_{L^p(\Omega,\mathbb{R}^N)},
\end{equation}
and
\begin{equation}\label{Poin2}
\left\|u-\fint_{\Omega}u\right\|_{L^p(\Omega)}\leq \tilde{C}_P\| \nabla u\|_{L^p(\Omega,\mathbb{R}^N)},
\end{equation}
and, finally,
\begin{equation}\label{Poin3}
\left\|u-\fint_{\partial \Omega}u\right\|_{L^p(\partial \Omega)}\leq \hat{C}_P\| \nabla u\|_{L^p(\Omega,\mathbb{R}^N)}.
\end{equation}
\end{teo}

The proof of this result is an easy modification of the argument used to prove \cite[Proposition~3.2]{A-Mor-Ros} where
\eqref{Poin2} and \eqref{Poin3} are proved in the case $p=2$. We sketch the proof in the Appendix.

We call $(W^{1,p}(\Omega))'$ the dual to $W^{1,p}(\Omega)$ and we call $(W^{1,p}(\Omega))'_{\ast}$ the subspace of $F\in (W^{1,p}(\Omega))'$ such that
$$\langle F,1\rangle
=0.$$
For any $u\in W^{1,p}(\Omega)$, we have
$\langle F,u\rangle= \langle F,u_{\ast}\rangle$,
hence
$$\langle F,u\rangle\leq \|F\|_{(W^{1,p}_{\ast}(\Omega))'}\|\nabla u\|_{L^p(\Omega)},$$
and $(W^{1,p}(\Omega))'_{\ast}$ is endowed with the $(W^{1,p}_{\ast}(\Omega))'$ norm.

We call $(TW^{1,p}(\Omega))'$ the dual to $TW^{1,p}(\Omega)$.
We call $(TW^{1,p}(\Omega))'_{\ast}$  the subspace of $g\in (TW^{1,p}(\Omega))'$ such that
$$\langle g,1\rangle
=0.$$
For the same reason as before, we endow
$(TW^{1,p}(\Omega))'_{\ast}$ with the $(TW^{1,p}_{\ast}(\Omega))'$ norm.

We note that $(TW^{1,p}(\Omega))'$ is immersed in $(W^{1,p}(\Omega))'$, with immersion constant $1$, 
if we identify any $g\in (TW^{1,p}(\Omega))'$ with $F_g\in  (W^{1,p}(\Omega))'$ given by
\begin{equation}\label{F_gdefin}
\langle F_g,u\rangle=\langle g,u|_{\partial\Omega}\rangle\quad\text{for any }u\in W^{1,p}(\Omega).
\end{equation}
Analogously, $(TW^{1,p}(\Omega))'_{\ast}$ is immersed in $(W^{1,p}(\Omega))'_{\ast}$, with immersion constant $1$.

For $p=1$, $(TW^{1,1}(\Omega))'$ can be identified with $L^{\infty}(\partial\Omega)$ and 
$(TW^{1,1}(\Omega))'_{\ast}$ can be identified with $L^{\infty}_{\ast}(\partial\Omega)$,
with norms which are equivalent for some constants depending on $r$, $L$ and $R$ only.

For $1<p<+\infty$, we note that $L^{q'}(\partial\Omega)$ is continuously, and compactly respectively, immersed in 
$(TW^{1,p}(\Omega))'$ whenever $TW^{1,p}(\Omega)$ is continuously, and compactly respectively, immersed in $L^q(\partial\Omega)$. Continuous immersion holds for all cases in \eqref{pqcases}, in particular, this is true for $q=p$, where compact immersion holds. The immersion constants are exactly given by $c_{p,q}C_T$. In this case, any $g\in L^{q'}(\partial\Omega)$ is identified with the linear operator
\begin{equation}
\langle g,\varphi\rangle=
\int_{\partial\Omega}g\varphi\quad\text{for any }\varphi\in TW^{1,p}(\Omega).
\end{equation}
We note that
$L^{q'}_{\ast}(\partial\Omega)\subset (TW^{1,p}(\Omega))'_{\ast}$, with continuous or compact immersion, for the same values of $p$ and $q$ as before.

When $p=2$, we use following alternative notation:
$H^1(\Omega)=W^{1,2}(\Omega)$, $H^1_0(\Omega)=W^{1,2}_0(\Omega)$, $H^{1/2}(\partial\Omega)=TW^{1,2}(\Omega)$,
$H^{-1/2}(\partial\Omega)=(TW^{1,2}(\Omega))'$.
We also use
$H^1_{\ast}(\Omega)$, $H^{1/2}_{\ast}(\partial\Omega)$, $(H^1(\Omega))'_{\ast}$ and $H^{-1/2}_{\ast}(\partial\Omega)=(H^{1/2}(\partial\Omega))'_{\ast}$.

\subsection{Discretizations of the domain and corresponding finite element spaces}\label{finiteelsec}

\begin{defin}
We say that a bounded domain $\Omega\subset\R^N$ is \emph{polyhedral} if
$\Omega$ coincides with the interior of $\overline{\Omega}$ and
its boundary is the finite union of cells, any \emph{cell} being the closure of an open connected subset of an $(N-1)$-dimensional hyperplane.
\end{defin}

We use standard conforming piecewise linear finite elements, for which we refer for instance to \cite[Chapter~2]{Cia}. Let, for the time being, $\Omega$ be just a bounded domain of $\R^N$ such that $\Omega$ coincides with the interior of $\overline{\Omega}$.

\begin{defin}
A finite set $\mathcal{T}$ of subsets of $\overline{\Omega}$ is a \emph{triangulation} of $\overline{\Omega}$
if the following holds
\begin{itemize}
\item $\overline{\Omega}=\displaystyle{\bigcup_{K\in \mathcal{T}}K}$;
\item each $K\in \mathcal{T}$ is a closed
$N$-simplex which is nondegenerate, that is, it has nonempty interior;
\item the intersection of two different elements of $\mathcal{T}$ is either empty or consists of a common face.
 \end{itemize}
\end{defin}

To any triangulation $\mathcal{T}$ of $\overline{\Omega}$, we can associate the following finite element space $X^{\mathcal{T}}$ given by
$$X^{\mathcal{T}}=\{v\in C(\overline{\Omega}):\ v|_K\in P_1(K)\text{ for any }K\in\mathcal{T}\}$$
where $P_1(K)$ is the space of polynomials of order at most $1$ restricted to $K$.
By \cite[Theorem~2.2.3]{Cia} we have that $X^{\mathcal{T}}\subset C(\overline{\Omega})\cap H^1(\Omega)$. It is also clear that $X^{\mathcal{T}}_{0}=\{v\in X^{\mathcal{T}}:\ v|_{\partial\Omega}=0\}$
is contained in $H^1_0(\Omega)$. 

For any $K\in\mathcal{T}$ we call
$$h_K=\mathrm{diam}(K)\quad\text{and}\quad\rho_K=\sup
\{\mathrm{diam}(B):\ B\text{ is a ball contained in }K\}.$$ 
We say that  $\mathcal{T}$ is a
\emph{regular triangulation} of $\overline{\Omega}$ with positive constants $s$ and $h$ if
\begin{equation}\label{regularfinite}
h_K\leq h\text{ and }h_K\leq s\rho_K\quad\text{for any }K\in\mathcal{T}.
\end{equation}

Since for any ball $B$ we have $|B|=|B_1|(\mathrm{diam}(B)/2)^N$, we have, for any $K\in\mathcal{T}$,
$$h_K^N\leq s^N\rho_K^N
\quad\text{and}\quad
|B_1|\left(\frac{\rho_K}2\right)^N\leq|K|\leq |B_1|h_K^N.$$

It is also important to note that $\mathcal{T}$ induces an $(N-1)$ dimensional triangulation $\mathcal{T}|_{\partial\Omega}$ of the boundary $\partial\Omega$. Furthermore, for any $s>0$ there exists $\tilde{s}>0$, depending on $s$ and $N$ only, such that if
$\mathcal{T}$ is a
regular triangulation of $\overline{\Omega}$ with positive constants $s$ and $h$, then $\mathcal{T}|_{\partial\Omega}$ is a
regular triangulation of $\partial\Omega$ with positive constants $\tilde{s}$ and $h$.

We call $\Pi_{\mathcal{T}}$ the associated interpolation operator
defined on $C(\overline{\Omega})$.
The following estimate is an immediate consequence of \cite[Theorem~3.1.6]{Cia}.

\begin{teo}\label{Ciarletestimteo}
Let $\mathcal{T}$ be a
regular triangulation of $\overline{\Omega}$ with positive constants $s$ and $h$.

Let us consider a real number $q$ such that $q>N/2$.
Then there exists a constant $\tilde{C}(q)$, depending on $q$ only, such that
for any $u\in W^{2,q}(\Omega)$ we have
\begin{multline}\label{Ciarletest}
\|u-\Pi_{\mathcal{T}} u\|_{L^{q}(\Omega)}\leq \tilde{C}(q)h^2\|D^2u\|_{L^q(\Omega,\mathbb{M})}\quad\text{and}\\
\|\nabla(u-\Pi_{\mathcal{T}} u)\|_{L^{q}(\Omega,\R^N)}\leq \tilde{C}(q)sh\|D^2u\|_{L^q(\Omega,\mathbb{M})}.
\end{multline}
 \end{teo}

\begin{defin}\label{sdiscretizable}
We say that a bounded domain $\Omega$, such that
$\Omega$ coincides with the interior of $\overline{\Omega}$, is \emph{discretizable} with constant $s>0$ if for any $h\in (0,1]$ there exists
a regular triangulation $\mathcal{T}_h$ of $\overline{\Omega}$ with constants $s$ and $h$. We call $X^h=X^{\mathcal{T}_h}$
and $\Pi_h=\Pi_{\mathcal{T}_h}$, for any $h\in (0,1]$.
\end{defin}

\begin{oss}\label{oss2} If $\Omega$ is a bounded domain such that
$\Omega$ coincides with the interior of $\overline{\Omega}$ and
$\overline{\Omega}$ admits a triangulation, then $\Omega$ must be polyhedral.
On the other hand, for any bounded domain $\Omega$ which is polyhedral, there exists a constant $s>0$, depending on $\Omega$, such that
$\Omega$ is discretizable with constant $s$. Clearly the dependence on $s$ from $\Omega$ is rather involved and finding the best possible constant might not be an easy task. In any case, the sketch of the proof of this claim is the following.

It is well-known that any convex polytope admits a triangulation. This can be proved by induction. For $N=2$, just take an interior point and connect it with the vertices. For a convex polytope of dimension $N$, assume that we have a triangulation of its boundary. Then we take an interior point and connect it with the vertices of the triangulation of the boundary.  We note that the faces of a convex polytope are, in turn, convex polytopes of their corresponding dimensions. Hence we argue in the following way. We start by a triangulation of all $2$ dimensional faces, with the technique introduced before. Then we move on by constructing a triangulation of $3$ dimensional faces, until we obtain a triangulation of the $N$-dimensional polytope.

By taking all hyperplanes containing one of the cells forming the boundary of $\Omega$, we can decompose $\overline{\Omega}$ in the union of a finite number of convex polytopes. Intersection of two of these polytopes is through a common face only.
A little care is needed to obtain a triangulation of $\overline{\Omega}$ through suitable triangulations of these convex polytopes. As before,
we begin by triangulating all $2$-dimensional faces of this collection of polytopes. We then proceed to triangulate all $3$-dimensional faces by choosing one of their interior points and using the previous triangulation of their boundaries. If we proceed iteratively in this way, when we reach the dimension $N$, we indeed obtain
a regular triangulation $\mathcal{T}_0$, with some positive constants $s_0$ and $h_0$, of $\overline{\Omega}$. We assume that actually $h_0$ is reached for some $K_0\in \mathcal{T}_0$.

Let us consider such an initial triangulation $\mathcal{T}_0$. We can subdivide any of its
$N$-simplices into $2^N$ simplices with same $N$-dimensional measure  by the procedure given  in \cite{Edel1,Edel2}. Again, some care is needed to obtain a new global triangulation of $\overline{\Omega}$,
but this can be done by following, for instance, \cite[Theorem~3.6]{Gra}. We call $\mathcal{T}_1$ this new triangulation. We iterate this construction to obtain a sequence of triangulations $\{\mathcal{T}_n\}_{n\geq 0}$. By using the results in \cite{Edel1,Edel2}, it can be shown that there exists a constant
$s$, $0<s\leq s_0$, such that for any $n\geq 0$ we have that
$$h_K\leq s\rho_K\quad\text{for any }K\in\mathcal{T}_n.$$

Let us call, for any $n\geq 0$,
$$h_n=\max\{h_K:\ K\in\mathcal{T}_n\}\quad\text{and}\quad\tilde{h}_n=\min\{h_K:\ K\in\mathcal{T}_n\}.$$
Then
$$|B_1|\left(\frac{\tilde{h}_0}{2s}\right)^N\leq |K|\leq |B_1|h_0^N\quad\text{for any }K\in \mathcal{T}_0,$$
therefore, for any $n\geq 1$
$$|B_1|\left(\frac{\tilde{h}_0}{2^{n+1}s}\right)^N\leq |K|\leq |B_1|\left(\frac{h_0}{2^{n}}\right)^N\quad\text{for any }K\in \mathcal{T}_n.$$
We conclude that, for any $n\geq 1$,
\begin{equation}\label{htri}
\frac{1}{2s}\frac{\tilde{h}_0}{2^{n}}\leq \tilde{h}_n\leq h_n\leq 2s\frac{h_0}{2^{n}}.
\end{equation}
The proof of the claim can therefore be easily concluded. We finally note that the same argument can be applied to $(\mathcal{T}_n)|_{\partial\Omega}$, where $s$ is replaced by $\tilde{s}$, and completely analogous estimates hold with $N$ replaced by $N-1$.
\end{oss}

We just recall the following classical result, which is usually known as C\'ea's Lemma.
\begin{teo}\label{cea} Let $H_{\ast}$ be a Hilbert space and $X_{\ast}$ be a closed subspace of $H_{\ast}$. Let $b:H_{\ast}\times H_{\ast}$ be a bilinear symmetric form which is coercive and continuous, that is, for some constants $0<c_1\leq c_2$ and any $u$, $w\in H_{\ast}$
$$b(u,w)\leq c_2\|u\|_{H_{\ast}}\|w\|_{H_{\ast}}\quad\text{and}\quad b(u,u)\geq c_1\|u\|_{H_{\ast}}^2.$$

Let $F\in (H_{\ast})'$ and let $u\in H_{\ast}$ and $u_X\in X_{\ast}$ solve, respectively,
$$b(u,w)=F[w]\quad\text{for any }w\in H_{\ast}$$ and
$$b(u_X,w)=F[w]\quad\text{for any }w\in X_{\ast}.$$

Then
\begin{equation}\label{ceaeq}
\|u-u_X\|_{H_{\ast}}\leq \sqrt{\frac{c_2}{c_1}}\inf\{\|u-v\|_{H_{\ast}}:\ v\in X_{\ast}\}.
\end{equation}
\end{teo}

\subsection{The continuum models}\label{continuumsec}
Let $\Omega\subset\mathbb{R}^N$ belong to $\mathcal{A}(r,L,R)$ and let $A$ be a symmetric conductivity tensor in $\Omega$ belonging to $\mathcal{M}(\lambda_0,\lambda_1)$. In order to consider discretized version of boundary value problems, we also consider 
$X$ which is a closed subspace of $H^1(\Omega)$ containing the constant functions. Usually $X$ will be taken as the finite element space $X^{\mathcal{T}}$ corresponding to a triangulation $\mathcal{T}$ of $\overline{\Omega}$.
For such subspace $X$,
we call $X_{\ast}=X\cap H^1_{\ast}(\Omega)$ and $X_0=X\cap H^1_0(\Omega)$. We have that $X$ is endowed with the $H^1(\Omega)$ norm, whereas $X_{\ast}$ and $X_0$ are endowed with the $H^1_{\ast}(\Omega)$ and $H^1_0(\Omega)$ norm, respectively.
Analogously, we define $TX$ and $TX_{\ast}$ as the corresponding trace spaces, with the usual quotient norm.
The dual to $X$ is $X'$, whereas $X'_{\ast}$ is the subspace of $X'$ consisting of functional which are zero on constants. The space $X'_{\ast}$
is endowed with the $(X_{\ast})'$ norm. The dual to $TX$ is $(TX)'$, whereas $(TX)'_{\ast}$ is the subspace of $(TX)'$ consisting of functional which are zero on constants. The space $(TX)'_{\ast}$
is endowed with the $(TX_{\ast})'$ norm.

\begin{oss}\label{tracenormsM}
For any $u\in TX$, we have that $u\in H^{1/2}(\partial\Omega)$ and
$$\|u\|_{H^{1/2}(\partial\Omega)}\leq \|u\|_{TX}\quad\text{for any }u\in TX.$$
Correspondingly, for any $g\in H^{-1/2}(\partial\Omega)$, we have that $g\in (TX)'$ and
$$\|g\|_{(TX)'}\leq \|g\|_{H^{-1/2}(\partial\Omega)}\quad\text{for any }g\in H^{-1/2}(\partial\Omega).$$
The same properties and inequalities apply to the corresponding zero mean spaces as well.

Whether reserve inequalities holds and, if this is the case, with which constants, that could strongly depend on $X$.
\end{oss}

We consider the bilinear form
$$b_{A}(u,v):=\int_{\Omega}A\nabla u\cdot \nabla v\quad\text{for any }u,\,v\in H^1(\Omega).$$
It is immediate to show that such a form is symmetric and continuous on $H^1(\Omega)\times H^1(\Omega)$. In fact,
for any $u,\, v\in H^1(\Omega)$, we have
$$\int_{\Omega}A\nabla u\cdot \nabla v\leq \lambda_1\|\nabla u\|_{L^2(\Omega,\R^N)}\|\nabla v\|_{L^2(\Omega,\R^N)}$$
and
$$\lambda_0\|\nabla u\|^2_{L^2(\Omega,\R^N)}\leq \int_{\Omega}A\nabla u\cdot \nabla u.$$
Therefore, for any $X$, such a form is coercive both on $X_{\ast}$ and $X_0$.

\subsubsection{Dirichlet problems and Dirichlet-to-Neumann maps}
For any $\varphi\in X$ and any $F\in (H^1(\Omega))'$, let $u=u(F,\varphi)\in X$ be
such that $u-\varphi\in X_0$ and
$$\int_{\Omega}A\nabla u\cdot\nabla w=F[w]\quad\text{for any }w\in X_0.$$
We have that such a $u$ exists and is unique, since $u-\varphi$ satisfies
$$\int_{\Omega}A\nabla (u-\varphi)\cdot \nabla w=F(w)-\int_{\Omega}A\nabla \varphi\cdot \nabla w\quad\text{for any }w\in X_0.$$
We conclude that
$$\lambda_0\|u-\varphi\|_{X_0}^2\leq \|F\|_{(H^1(\Omega))'}\|u-\varphi\|_{X}+\lambda_1\|\nabla \varphi\|_{L^2(\Omega,\R^N)}\|u-\varphi\|_{X_0}$$
hence
$$\|u-\varphi\|_{X_0}\leq \frac{1}{\lambda_0}\left((1+c_P^2)^{1/2}\|F\|_{(H^1(\Omega))'}+\lambda_1\|\nabla \varphi\|_{L^2(\Omega,\R^N)}\right)$$
and, finally,
$$\|\nabla u\|_{L^2(\Omega,\R^N)}\leq \frac{1}{\lambda_0}\left((1+c_P^2)^{1/2}\|F\|_{(H^1(\Omega))'}+(\lambda_1+\lambda_0)\|\nabla \varphi\|_{L^2(\Omega,\R^N)}\right)$$
and
$$\|u\|_{X}\leq \|\varphi\|_X+\frac{(1+c_P^2)^{1/2}}{\lambda_0}\left((1+c_P^2)^{1/2}\|F\|_{(H^1(\Omega))'}+\lambda_1\|\nabla \varphi\|_{L^2(\Omega,\R^N)}\right).$$
It would be enough to consider $F\in (X_0)'$ and one can replace in the last three inequalities
$(1+c_P^2)^{1/2}\|F\|_{(H^1(\Omega))'}$ with $\|F\|_{(X_0)'}$.
 
When $X=H^1(\Omega)$, then $u$ is the weak solution to
\begin{equation}\label{Dirpbm}
\left\{\begin{array}{ll}
-\mathrm{div(}A\nabla u)=F &\text{in }\Omega\\
u=\varphi  & \text{on }\partial\Omega.
\end{array}
\right.
\end{equation}

When $F=0$, let $u(\varphi)=u(0,\varphi)$ and 
let us denote $\Lambda_X(A):TX\to (TX)'$ be the linear operator such that, for any $\varphi$, $\psi\in TX$,
$$\langle \Lambda_X(A)[\varphi],\psi\rangle_{TX',TX}=\int_{\Omega}A\nabla u(\varphi)\cdot \nabla w,$$
with $w\in X$ such that $w|_{\partial\Omega}=\psi$. In fact,
for any $\varphi_0\in X_0$ and any $w_0\in X_0$, we have
$$\int_{\Omega}A\nabla u(\varphi)\cdot \nabla w=\int_{\Omega}A\nabla u(\varphi+\varphi_0)\cdot \nabla (w+w_0).$$
We note that $\Lambda_X(A)$ is actually a map from $TX$ into $(TX)'_{\ast}$.
The norm of $\Lambda_X(A)$ can be estimated as follows
\begin{multline*}
\|\Lambda_X(A)\|_{\mathcal{L}(TX,(TX)')}\leq \|\Lambda_X(A)\|_{\mathcal{L}(TX,(TX)'_{\ast})}\\\leq
\|\Lambda_X(A)\|_{\mathcal{L}(TX_{\ast},(TX)'_{\ast})}\leq \frac{\lambda_1}{\lambda_0}(\lambda_1+\lambda_0).
\end{multline*}

When $X=H^1(\Omega)$, we call $\Lambda(A)=\Lambda_{H^1(\Omega)}(A):H^{1/2}(\partial\Omega)\to H_{\ast}^{-1/2}(\partial\Omega)$ which is the so-called \emph{Dirichlet-to-Neumann} map associated to the conductivity $A$. In fact, for any Dirichlet datum $\varphi\in H^{1/2}(\partial\Omega)$,
$\Lambda(A)[\varphi]$ formally corresponds to $A\nabla u\cdot\nu$ on $\partial\Omega$, which is the Neumann datum of $u=u(\varphi)$.
We call $$\|\Lambda(A)\|_{nat}:=\|\Lambda(A)\|_{\mathcal{L}(H^{1/2}_{\ast}(\partial\Omega),H^{-1/2}_{\ast}(\partial\Omega))}$$ the \emph{natural norm} on Dirichlet-to-Neumann maps.

\subsubsection{Neumann problems and Neumann-to-Dirichlet maps}
For any $F\in (H^1(\Omega))'_{\ast}$,
let $v=v(F)\in X_{\ast}$ be such that
\begin{equation}\label{Neuweak}
\int_{\Omega}A\nabla v\cdot\nabla w=F[w]\quad\text{for any }w\in X.
\end{equation}
We have that 
such a $v$ exists and it is unique.
We note that
$$\|v\|_{X_{\ast}}\leq \frac{1}{\lambda_0}(1+C_P^2)^{1/2}\|F\|_{(H^1(\Omega))'}.$$
It would be enough to consider $F\in X'_{\ast}$ and one can replace in the previous formula
$(1+C_P^2)^{1/2}\|F\|_{(H^1(\Omega))'}$ with $\|F\|_{X'_{\ast}}$.

When $X=H^1(\Omega)$, then $v=v(F)$ is the weak solution to
\begin{equation}\label{Neupbm}
\left\{\begin{array}{ll}
-\mathrm{div(}A\nabla v)=F &\text{in }\Omega\\
A\nabla v\cdot \nu=0  & \text{on }\partial\Omega.
\end{array}
\right.
\end{equation}

Let us now consider $g\in (TX)'_{\ast}$
 and 
let us define $F_g\in X'_{\ast}$ analogously as in \eqref{F_gdefin}. One can choose,  for instance, $g\in H^{-1/2}_{\ast}(\partial\Omega)$ and in this case $F_g\in (H^1(\Omega))'_{\ast}$.
Let $v(g)=v(F_g)$ and let us consider
 the map $\mathcal{N}_X(A):(TX)'_{\ast}\to TX_{\ast}$ such that
 $$\mathcal{N}_X(g)=v(g)|_{\partial\Omega}\quad\text{for any }g\in (TX)'_{\ast}.$$
 We immediately note that $\mathcal{N}_X(A)$ is a linear bounded operator and that
 $$\mathcal{N}_X(A)_{\mathcal{L}((TX)'_{\ast},TX_{\ast})}\leq \frac{1}{\lambda_0}.$$
We also have that $\mathcal{N}_X(A)$ is the inverse of $\Lambda_X(A)|_{TX_{\ast}}$.

Since $TX_{\ast}$ is continuously immersed in $H^{1/2}_{\ast}(\partial\Omega)$ with constant $1$, thus in $L^2_{\ast}(\partial\Omega)$ with constant $\hat{C}_P$ as in \eqref{Poin3}, we have that
\begin{equation}\label{L2L2N0}
\|\mathcal{N}_X\|_{\mathcal{L}(L^2_{\ast}(\partial\Omega),H^{1/2}_{\ast}(\partial\Omega))}\leq \frac{\hat{C}_P}{\lambda_0}
\quad\text{and}\quad \|\mathcal{N}_X\|_{\mathcal{L}(H^{-1/2}_{\ast}(\partial\Omega),L^2_{\ast}(\partial\Omega),)}\leq \frac{\hat{C}_P}{\lambda_0}
\end{equation}
and
\begin{equation}\label{L2L2N}
\|\mathcal{N}_X\|_{\mathcal{L}(L^2_{\ast}(\partial\Omega))}\leq \frac{\hat{C}_P^2}{\lambda_0}. 
\end{equation}

When $X=H^1(\Omega)$, we call $\mathcal{N}(A)=\mathcal{N}_{H^1(\Omega)}(A):H^{-1/2}_{\ast}(\partial\Omega)\to H^{1/2}_{\ast}(\partial\Omega)$ which is the so-called \emph{Neumann-to-Dirichlet} map associated to the conductivity $A$. In fact, for any Neumann datum $g\in H^{-1/2}_{\ast}(\partial\Omega)$,
$\mathcal{N}(A)[\varphi]$ is the corresponding Dirichlet datum of the solution $v=v(g)$.
We call $$\|\mathcal{N}(A)\|_{nat}:=\|\mathcal{N}(A)\|_{\mathcal{L}(H^{-1/2}_{\ast}(\partial\Omega),H^{1/2}_{\ast}(\partial\Omega))}$$ the \emph{natural norm} on Neumann-to-Dirichlet maps. However, as we shall see, for the Neumann-to-Dirichlet map is much more convenient to adopt what we call the $L^2\text{-}L^2$-norm, that is
\begin{equation}
\|\mathcal{N}(A)\|_{L^2\text{-}L^2}:=\|\mathcal{N}(A)\|_{\mathcal{L}(L^{2}_{\ast}(\partial\Omega))}.
\end{equation}
We have that
$$\|\mathcal{N}(A)\|_{L^2\text{-}L^2}\leq  \hat{C}_P^2 \|\mathcal{N}(A)\|_{nat}.$$
On the other hand, since $L^{2}_{\ast}(\partial\Omega)$ is dense in $H^{-1/2}_{\ast}(\partial\Omega)$, we have that, for any two conductivity tensors in $\Omega$,
$\|\mathcal{N}(A_1)-\mathcal{N}(A_2)\|_{L^2\text{-}L^2}=0$ if and only if $\|\mathcal{N}(A_1)-\mathcal{N}(A_2)\|_{nat}=0$ which is equivalent to
$\|\Lambda(A_1)-\Lambda(A_2)\|_{nat}=0$. Hence, from the point of view of uniqueness of the inverse problem, the use of the Dirichlet-to-Neumann map or the Neumann-to-Dirichlet map, even when restricted to $L^{2}_{\ast}(\partial\Omega)$, is completely equivalent.
From the point of view of stability, if we use natural norms, the use of the Dirichlet-to-Neumann map or the Neumann-to-Dirichlet map is completely equivalent, as shown in the next remark.

\begin{oss} We have that
$$\Lambda(A_1)-\Lambda(A_2)=\Lambda(A_1)(\mathcal{N}(A_2)-\mathcal{N}(A_1))\Lambda(A_2)$$
and
$$\mathcal{N}(A_1)-\mathcal{N}(A_2)=\mathcal{N}(A_1)(\Lambda(A_2)-\Lambda(A_1))\mathcal{N}(A_2).$$
Therefore, for a constant $C\geq 1$ depending on $r$, $L$, $R$, $\lambda_0$ and $\lambda_1$ only, we have
$$\frac1{C}\|\mathcal{N}(A_1)-\mathcal{N}(A_2)\|_{nat}\leq \|\Lambda(A_1)-\Lambda(A_2)\|_{nat}   \leq   C\|\mathcal{N}(A_1)-\mathcal{N}(A_2)\|_{nat}.$$
\end{oss}

We conclude this part with the following important result, which is essentially due to Meyers, \cite{NMey}. The extension of Meyers result to Neumann problems is in \cite{Gal-Mon}. Here we recall these classsical results, the only difference is that we, as usual, specify the dependence of involved constants on $\Omega$ through the geometric constants characterizing it. We postpone the proof to the Appendix.

\begin{teo}\label{Meyers}
Let $\Omega\in \mathcal{A}(r,L,R)$ and let $A\in \mathcal{M}(\lambda_0,\lambda_1)$. Then
there exists a constant $Q_1>2$, depending on $r$, $L$, $R$, $\lambda_0$ and $\lambda_1$ only, such that
for any $p$, $2\leq p\leq Q_1$, the following holds.

There exists a constant $\tilde{D}$, depending on
$r$, $L$, $R$, $\lambda_0$, $\lambda_1$ and $p$ only,
such that for any 
$\varphi\in W^{1,p}(\Omega)$ and $F\in (W^{1,p}_0(\Omega))'$ we have
for $u$ solution to \eqref{Dirpbm}
\begin{equation}\label{MeyDir}
\|u\|_{W^{1,p}(\Omega)}\leq \tilde{D}(\|\varphi\|_{W^{1,p}(\Omega)}+\|F\|_{(W^{1,p}_0(\Omega))'})
\end{equation}
and for any $F\in (W^{1,p}(\Omega))'_{\ast}$ we have for  $v$ solution to \eqref{Neupbm}
\begin{equation}\label{MeyNeu}
\|v\|_{W^{1,p}(\Omega)}\leq \tilde{D}\|F\|_{(W^{1,p}(\Omega))'_{\ast}}.
\end{equation}
\end{teo}

As already pointed out in \cite{Ron08,Ron15}, an important consequence is the following.
Under the assumptions of Theorem~\ref{Meyers},
for any $p$, $2\leq p\leq Q_1$, let $q$,
$2< q\leq +\infty$, be such that
\begin{equation}\label{qdef}
\frac{1}{q}+\frac{1}{p}+\frac{1}{2}=1.
\end{equation}
Then there exists a constant
$\tilde{D}_1$, depending on $r$, $L$, $R$, $\lambda_0$, $\lambda_1$ and $p$ only,
such that the following holds.

For any $A_1$, $A_2\in \mathcal{M}(\lambda_0,\lambda_1)$, any
$\varphi\in W^{1,p}(\Omega)$ and $F\in (W^{1,p}_0(\Omega))'$, we have
\begin{equation}\label{Direst}
\|u_1-u_2\|_{H^1(\Omega)}\leq \tilde{D}_1\left(\|\varphi\|_{W^{1,p}(\Omega)}+\|F\|_{(W_0^{1,p}\Omega))'}\right)\|A_1-A_2\|_{L^q(\Omega)},
\end{equation}
where $u_i$, $i=1,2$, is the solution to \eqref{Dirpbm} with $A$ replaced by $A_i$,
and, for any $F\in (W^{1,p}(\Omega))'_{\ast}$ we have
\begin{equation}\label{Neuest}
\|v_1-v_2\|_{H^1(\Omega)}\leq \tilde{D}_1\|F\|_{(W^{1,p}(\Omega))'_{\ast}}\|A_1-A_2\|_{L^q(\Omega)},
\end{equation}
where $v_i$, $i=1,2$, is the solution to \eqref{Neupbm} with $A$ replaced by $A_i$.

We easily note that in both previous inequalities, \eqref{Direst} and \eqref{Neuest}, if $2<p\leq Q_1$ then $q<+\infty$ and we can replace
$\|A_1-A_2\|_{L^q(\Omega)}$ by $(2\lambda_1)^{1-\beta}\|A_1-A_2\|_{L^1(\Omega)}^{\beta}$
where $\beta=1/q$, hence
\begin{equation}\label{betadef}
\beta=1-\frac1p-\frac12=\frac12-\frac1p=\frac{p-2}{2p}.
\end{equation}

Finally,  picking $p=Q_1$, 
for any $g\in L^2_{\ast}(\partial \Omega)$, we have that $F_g\in (W^{1,Q_1}\Omega))'_{\ast}$ and
$$\|F_g\|_{(W^{1,Q_1}\Omega))'_{\ast}}\leq C(Q_1)\|g\|_{L^2_{\ast}(\partial\Omega)}\quad\text{for any }g\in L^2_{\ast}(\partial\Omega).$$
Here $C(Q_1)$ depends on $Q_1$, $r$, $L$ and $R$ only.
Therefore, we conclude that
\begin{multline}\label{Neuest2}
\|v_1-v_2\|_{L^2_{\ast}(\partial\Omega)} \leq C_T(1+C_P^2)^{1/2}\|v_1-v_2\|_{H^{1/2}_{\ast}(\partial\Omega)}\\\leq 
C_T(1+C_P^2)^{1/2}\|v_1-v_2\|_{H^1(\Omega)}\\
\leq C_T(1+C_P^2)^{1/2}\tilde{D}_1C(Q_1)(2\lambda_1)^{1-\beta_1}\|g\|_{L^2_{\ast}(\partial\Omega)}\|A_1-A_2\|^{\beta_1}_{L^1(\Omega)}.
\end{multline}
where
\begin{equation}\label{beta1def}
\beta_1=\frac12-\frac1{Q_1}=\frac{Q_1-2}{2Q_1}.
\end{equation}

Therefore we can easily deduce the following Lipschitz and H\"older continuity results.

\begin{teo}\label{pcontinuity}
Let $\Omega\in\mathcal{A}(r,L,R)$.
Let us consider two conductivity tensors $A_1$ and $A_2\in \mathcal{M}(\lambda_0,\lambda_1)$.

We have the following continuity properties. First, $\Lambda$ and $\mathcal{N}$ are Lipschitz continuous with respect to the $L^{\infty}(\Omega)$ norm on $\mathcal{M}(\lambda_0,\lambda_1)$ and the natural operator norms, that is
$$\|\Lambda(A_1)-\Lambda(A_2)\|_{\mathcal{L}(H^{1/2}_{\ast}(\partial \Omega),H^{-1/2}_{\ast}(\partial \Omega))}\leq C
\|A_1-A_2\|_{L^{\infty}(\Omega)}$$
and
$$\|\mathcal{N}(A_1)-\mathcal{N}(A_2)\|_{\mathcal{L}(H^{-1/2}_{\ast}(\partial \Omega),H^{1/2}_{\ast}(\partial \Omega))}\leq C
\|A_1-A_2\|_{L^{\infty}(\Omega)}$$
where $C$ depends on $r$, $L$, $R$, $\lambda_0$ and $\lambda_1$ only.

Then, $\Lambda$ and $\mathcal{N}$ are H\"older continuous with respect to the $L^1(\Omega)$ norm on $\mathcal{M}(\lambda_0,\lambda_1)$ and the following norms. Fixed $p$, $2<p\leq Q_1$,
$$\|\Lambda(A_1)-\Lambda(A_2)\|_{\mathcal{L}(TW^{1,p}_{\ast}(\partial \Omega),H^{-1/2}_{\ast}(\partial \Omega))}\leq C
\|A_1-A_2\|^{\beta}_{L^1(\Omega)}$$
and
$$\|\mathcal{N}(A_1)-\mathcal{N}(A_2)\|_{\mathcal{L}((TW^{1,p}(\partial \Omega))'_{\ast},H^{1/2}_{\ast}(\partial \Omega))}\leq C
\|A_1-A_2\|^{\beta}_{L^1(\Omega)}$$
where $C$ depends on $r$, $L$, $R$, $\lambda_0$, $\lambda_1$ and $p$ only, whereas $\beta$ is given in \eqref{betadef}, thus it depends on 
$p$ only.

In particular, choosing $p=Q_1$, we have for the $L^2\text{-}L^2$ norm
\begin{multline*}
\|\mathcal{N}(A_1)-\mathcal{N}(A_2)\|_{\mathcal{L}((L^2_{\ast}(\partial\Omega),L^2_{\ast}(\partial \Omega))}\\\leq C_T(1+C_P^2)^{1/2}
\|\mathcal{N}(A_1)-\mathcal{N}(A_2)\|_{\mathcal{L}((L^2_{\ast}(\partial\Omega),H^{1/2}_{\ast}(\partial \Omega))}\\
\leq 
C_T(1+C_P^2)^{1/2}
\tilde{D}_1C(Q_1)(2\lambda_1)^{1-\beta_1}
\|A_1-A_2\|^{\beta_1}_{L^1(\Omega,\mathbb{M}^{N\times N}(\R))}.
\end{multline*}
Here all constants depend on $r$, $L$, $R$, $\lambda_0$ and $\lambda_1$ only, whereas $\beta_1$ is given in \eqref{beta1def}, thus it also depends on 
$r$, $L$, $R$, $\lambda_0$ and $\lambda_1$ only.
\end{teo}

\section{Experimental measurements (Complete Electrode Model)}\label{secCEM}

In order to discretize the measurements, we use the so-called \emph{experimental measurements}, which have been introduced in \cite{Som e Che e Isa}. These measurements, also known as the \emph{Complete Electrode Model} or CEM, correspond to the actual data one can obtain from the experiments. These discrete data are encoded into a matrix, called the \emph{resistance matrix}, which depends on the conductivity tensor.

There are strict and important relationships between Neumann-to-Dirichlet maps and resistance matrices. On the one hand, the error on the Neumann-to-Dirichlet maps, with respect to the $L^2\text{-}L^2$ norm, controls the error on the resistance matrices. This is essentially proved in \cite{Ron15}, here we just make the estimate more precise, see Subsection~\ref{L2L2ressubs}. On the other hand, provided we position the electrodes in a suitable way and we increase its number, the corresponding resistance matrices approximate the continuum measurements of the Neumann-to-Dirichlet map. Such a result has been proved in \cite{Hyv1,Hyv2,Hyv-et-al}, and is just slightly generalized in Subsection~\ref{approxNtoDsubs}.

We describe the CEM, for more details we refer to \cite{Som e Che e Isa}. Let $\Omega\in \mathcal{A}(r,L,R)$ and let $A$ be a symmetric conductivity tensor in $\Omega$ belonging to $\mathcal{M}(\lambda_0,\lambda_1)$.

We apply on the boundary of the conductor $M$ \emph{electrodes}, $e_m$, $m$ being here and in what follows an index such that $m=1,\ldots,M$.  The electrodes are identified with their \emph{contact regions}, that is, with subsets of $\partial\Omega$. These subsets $e_m$
are open and nonempty and pairwise disjoint. Further hypotheses on the electrodes will be introduced when needed. 
A current is sent to the body through the electrodes and the corresponding voltages
are measured on the same electrodes. 
The current applied to the electrode $e_m$ is denoted by
$I_m$ and the voltage measured on the electrode is denoted by $V_m$.
The \emph{current pattern} is given by 
 the column vector $I\in \R^M$ whose
components are $I_m$. $I$ has to satisfy the compatibility condition $\sum_{m=1}^MI_m=0$.
The corresponding \emph{voltage pattern}, the column vector $V\in\R^M$ whose components are $V_m$,
is determined up to an additive constant and we always choose to normalize it in such a way that
$\sum_{m=1}^MV_m=0$. The voltage pattern depends on the current pattern in a linear way, that is,
$V=RI$ where $R=R(A)\in \mathbb{M}_{sym}^{M\times M}(\R)$ is the \emph{resistance matrix}. Without loss of generality we assume that 
$R[1]=0$, where $[1]$ denotes the column vector whose components are all equal to $1$.

The contact between each electrode $e_m$ and the boundary is described by a \emph{surface impedance} $z_m$. For any $m$, we assume $z_m$ to be constant and 
such that
\begin{equation}\label{zl}
Z_1\leq z_m\leq Z_2,
\end{equation}
 for $0<Z_1<Z_2$ given constants which we assume to be independent on the size of the electrodes.

We call $e$ the union of all electrodes, that is,
\begin{equation}\label{edef}
e:=\bigcup_{m=1}^Me_m.
\end{equation}
Important geometric quantities related to the electrodes are the following
\begin{equation}\label{deltadef}
\delta:=\max_{m=1,\ldots,M}\mathrm{diam}(e_m)
\end{equation}
and
\begin{equation}\label{mudef}
\mu:=
\frac{\displaystyle \max_{m=1,\ldots,M}\{\mathcal{H}^{N-1}(e_m)\}}{ \displaystyle\min_{m=1,\ldots,M}\{\mathcal{H}^{N-1}(e_m)\}}.
\end{equation}

If the current pattern $I$ is applied on the electrodes, then the voltage $u$ inside $\Omega$ is the solution to the following
boundary value problem
\begin{equation}\label{expmeaspbm}
\left\{\begin{array}{ll}
-\mathrm{div}(A\nabla u)=0 &\text{in }\Omega,\\
u+z_m A\nabla u\cdot\nu=U_m &\text{on }e_m,\ m=1,\ldots,M,\\
A\nabla u\cdot\nu=0 &\text{on }\partial\Omega\backslash\bigcup_{m=1}^Me_m,\\
\int_{e_m}A\nabla u\cdot\nu=I_m&\text{for any }m=1,\ldots,M,
\end{array}\right.
\end{equation}
where $U_m$ are constants to be determined. We call $U\in \R^M$ the column vector whose
components are $U_m$.

For any $m$ we have
$V_m=\int_{e_m}u$, thus, by \eqref{expmeaspbm},
$$V_m=\mathcal{H}^{N-1}(e_m)U_m-z_mI_m.$$

By \cite[Theorem~3.3]{Som e Che e Isa}, there exists a unique pair $(u,U)\in H^1(\Omega)\times \R^M$,
satisfying
$\sum_{m=1}^M\left(\mathcal{H}^{N-1}(e_m)U_m-z_mI_m\right)=0$, such that \eqref{expmeaspbm} is satisfied. Thus the current pattern
$I$ uniquely determines the voltage pattern $V$ with the given normalization.
Furthermore, it has been also proved in \cite{Som e Che e Isa} that the relation
between $I$ and $V$ is linear, thus the resistance matrix $R(A)$ is well defined, and that $R(A)$ is symmetric. In what follows, the fact that $R(A)$ is symmetric will not play a significant role.

We recall the argument, because in the sequel we need some precise estimates.

We need to consider the discretized version of the problem for the CEM, thus we consider all estimates also for any  $X$, $X$ being a closed subspace of $H^1(\Omega)$ containing the constants. Then we call $H=X\times \R^M$ and $H_{\ast}=X_{\ast}\times \R^M$.
For any
$(u,U)$, $(w,W)\in H$, we let 
$$B_{A}((u,U),(w,W))=\int_{\Omega}A\nabla u\cdot\nabla w+
\sum_{m=1}^M\frac{1}{z_m}\int_{e_m}(U_m-u)(W_m-w).$$
We prove that $B_{A}$ is a continuous bilinear form on $H\times H$ and that it is coercive on $H_{\ast}\times H_{\ast}$. In fact,
\begin{multline*}
B_{A}((u,U),(w,W))\leq \lambda_1\|\nabla u\|_{L^2(\Omega,\R^N)}\|\nabla v\|_{L^2(\Omega,\R^N)}
\\+\frac1{Z_1}\sum_{m=1}^M\left(\int_{e_m}(U_m-u)^2\right)^{1/2}\left(\int_{e_m}(W_m-w)^2\right)^{1/2}.
\end{multline*}
By Cauchy-Schwarz inequality,
\begin{multline*}
\sum_{m=1}^M\left(\int_{e_m}(U_m-u)^2\right)^{1/2}\left(\int_{e_m}(W_m-w)^2\right)^{1/2}\\\leq 
\left(\sum_{m=1}^M\int_{e_m}(U_m-u)^2\right)^{1/2}\left(\sum_{m=1}^M\int_{e_m}(W_m-w)^2\right)^{1/2}.
\end{multline*}
We note that
\begin{multline*}\sum_{m=1}^M\int_{e_m}(U_m-u)^2\leq 2\sum_{m=1}^M\left(\int_{e_m}u^2+\mathcal{H}^{N-1}(e_m)|U_m|^2\right)\\\leq
2\left( \|u\|_{L^2(\partial\Omega)}^2+\max_{m=1,\ldots,M}\{\mathcal{H}^{N-1}(e_m)\}\|U\|^2_{\R^M}\right)\\\leq
2\max\{C^2_T,\mathcal{H}^{N-1}(\partial\Omega)\}\left(\|u\|_{H^1(\Omega)}^2+\|U\|^2_{\R^M}\right)\\=
2\max\{C^2_T,\mathcal{H}^{N-1}(\partial\Omega)\}\|(u,U)\|^2_{H}.
\end{multline*}
We conclude that for any $(u,U)$ and $(w,W)$ in $H$ we have
$$B_{A}((u,U),(w,W))\leq C_1\|(u,U)\|_{H}\|(w,W)\|_{H}$$
where
$$C_1=\lambda_1+\frac{2\max\{C^2_T,\mathcal{H}^{N-1}(\partial\Omega)\}}{Z_1}.$$

On the other hand, let $u\in H_{\ast}$. Then
$$\int_{\Omega}A\nabla u\cdot\nabla u\geq \lambda_0\|\nabla u\|_{L^2(\Omega,\R^N)}^2 
\geq \frac{\lambda_0}{1+C_P^2}\| u\|_{H^1(\Omega)}^2.$$
We take half of the right-hand side and note that
$$\int_{\Omega}A\nabla u\cdot\nabla u\geq 
 \frac{\lambda_0}{2(1+C_P^2)}\| u\|_{H^1(\Omega)}^2
+\frac{\lambda_0}{2C_T^2(1+C_P^2)}\|u\|_{L^2(\partial\Omega)}^2
.$$
The second term of the bilinear form can be estimated as follows
$$\sum_{m=1}^M\frac{1}{z_m}\int_{e_m}(U_m-u)^2\geq \frac1{Z_2}\sum_{m=1}^M\left(\int_{e_m}u^2+\mathcal{H}^{N-1}(e_m)|U_m|^2
-2U_m\int_{e_m}u\right).$$
But for any $\varepsilon>0$ we have
$$2U_m\int_{e_m}u\leq \varepsilon\mathcal{H}^{N-1}(e_m)|U_m|^2+\frac1{\varepsilon}\int_{e_m}u^2,$$
therefore, if we pick $\varepsilon$, $0<\varepsilon<1$ such that
$$\frac1{\varepsilon}=1+\dfrac{\lambda_0Z_2}{2C_T^2(1+C_P^2)},$$
then
$$
B_{A}((u,U),(u,U))\geq
 \frac{\lambda_0}{2(1+C_P^2)}\| u\|_{H^1(\Omega)}^2+\frac{1-\varepsilon}{Z_2}\min_{m=1,\ldots,M}\{\mathcal{H}^{N-1}(e_m)\}
 \|U\|^2_{\R^M}.
 $$
Hence, 
for any $(u,U)\in H_{\ast}$ we have
$$B_{A}((u,U),(w,W))\geq C_2\|(u,U)\|_{H}^2$$
where
\begin{equation}\label{coerccost}
C_2=\min\left\{\frac{\lambda_0}{2(1+C_P^2)},\frac{1-\varepsilon}{Z_2}\min_{m=1,\ldots,M}\{\mathcal{H}^{N-1}(e_m)\}   \right\}.
\end{equation}

Then, for any $g\in H^{-1/2}_{\ast}(\partial\Omega)$, or $g\in (TX)'_{\ast}$, and any $I\in\mathbb{R}^M$
such that $\sum_{m=1}^MI_m=0$, by using the Lax-Milgram Theorem in $H_{\ast}$, 
there exists a unique pair $(u_{\ast},U_{\ast})\in H_{\ast}$ satisfying, for any $(w,W)\in H_{\ast}$,
$$B_{A}((u_{\ast},U_{\ast}),(w,W))=\langle g, w|_{\partial\Omega}\rangle +\sum_{m=1}^MI_mW_m.$$
By the properties of $g$ and $I$, the previous formula actually holds for any $(w,W)\in H$ and, if $(u,U)\in H$ 
satisfies
\begin{equation}\label{eqform}
B_{A}((u,U),(w,W))=\langle g, w|_{\partial\Omega}\rangle+\sum_{m=1}^MI_mW_m\quad\text{for any }(w,W)\in H,
\end{equation}
then we have, for some constant $c$, $u=u_{\ast}+c$ and $U_m=(U_{\ast})_m+c$ for any $m$.
 
 Let $(u,U)\in H_{\ast}$ solve \eqref{eqform}. Then
\begin{multline*}
\|(u,U)\|_{H_{\ast}}^2\leq
\frac1{C_2}\mathcal{B}_A((u,U),(u,U))\\\leq  \frac1{C_2}\left(\|g\|_{H^{-1/2}_{\ast}(\partial\Omega)}\|\nabla u\|_{L^2(\Omega,\R^N)}+\|I\|_{\R^M}\|U\|_{\R^M}\right)
 \end{multline*}
hence
\begin{equation}\label{coerc}
\|\nabla u\|_{L^2(\Omega,\R^N)}\leq\|(u,U)\|_{H_{\ast}}\leq \frac{1}{C_2}\left( \|g\|^2_{H^{-1/2}_{\ast}(\partial\Omega)}+\|I\|^2_{\R^M}  \right)^{1/2}.
\end{equation}
Since $\lambda_0\|\nabla u\|^2_{L^2(\Omega,\R^N)}\leq \mathcal{B}_A((u,U),(u,U))$, when $I=0$, we have that 
\begin{equation}\label{coercbis}
\|\nabla u\|_{L^2(\Omega,\R^N)}\leq \frac1{\lambda_0}\|g\|_{H^{-1/2}_{\ast}(\partial\Omega)}.
\end{equation}
In the last two formulas, if $g\in (TX)'_{\ast}$, we can replace the $H^{-1/2}_{\ast}(\partial\Omega)$ norm with the $(TX)'_{\ast}$ norm.

We note that if $(u,U)\in H$ solves \eqref{eqform}, then, by choosing $W=0$, we have
\begin{equation}\label{diveq}
\int_{\Omega}A\nabla u\cdot\nabla w=0\quad \text{for any } w\in X_0,
\end{equation}
which, when $X=H^1(\Omega)$, reads as
\begin{equation}\label{diveq2}
-\mathrm{div}(A\nabla u)=0\quad\text{in }\Omega.
\end{equation}

We recall that $\Lambda_X(A)$ and 
$\mathcal{N}_X(A)$ are the Dirichlet-to-Neumann map and the Neumann-to-Dirichlet map associated to $A$ and to the subspace $X$.
Then $(u,U)\in H_{\ast}$
solves \eqref{eqform} if and only if $u$ satisfies \eqref{diveq}, we have
\begin{equation}\label{Uexp}
\mathcal{H}^{N-1}(e_m)U_m=z_mI_m+
\int_{e_m}u,
\quad\text{for any }m=1,\ldots,M,
\end{equation}
and the following equation holds in $(TX)'_{\ast}$ for $\phi=\Lambda_X(A)[u|_{\partial\Omega}]$
\begin{multline}\label{boundaryeqform}
\phi+\sum_{m=1}^M\frac{1}{z_m}\left(\mathcal{N}_X(A)[\phi]-
\frac{1}{\mathcal{H}^{N-1}(e_m)}\int_{e_m}\mathcal{N}_X(A)[\phi]\right)\chi_{e_m}
=\\ g+\sum_{m=1}^M\left(\frac{I_m}{\mathcal{H}^{N-1}(e_m)}\chi_{e_m}\right).
\end{multline}

As pointed out in \cite{Som e Che e Isa}, when $X=H^1(\Omega)$
we have that $(u,U)$ solves our direct problem \eqref{expmeaspbm}
for a given current pattern $I$, that is $I\in \mathbb{R}^M$ such that
$\sum_{m=1}^MI_m=0$,
if and only if
\eqref{eqform} is satisfied with $g=0$.
Therefore, we can define the following matrix $R_X(A)$ such that
$R_X(A)I=V$ where, for any $m=1,\ldots,M$,  
\begin{equation}\label{Vldef}
V_m=\int_{e_m}u+c\mathcal{H}^{N-1}(e_m),
\end{equation}
where $(u,U)$ solves \eqref{eqform} and 
$c$ is a constant which can be computed by imposing the condition
that $\sum_{m=1}^MV_m=0$, that is,
\begin{equation}\label{constantc0}
c=-\frac{\sum_{m=1}^M\int_{e_m}u}{\sum_{m=1}^M\mathcal{H}^{N-1}(e_m)}.
\end{equation}
The matrix $R_X(A)$ is then completed by setting $R_X[1]=0$.

When $X=H^1(\Omega)$, we call $R(A)=R_{H^1(\Omega)}(A)$ the \emph{resistance matrix} associated to $A$. Such a matrix corresponds to the experimental measurements that can be performed in practice.

\subsection{Estimating the resistance matrix by the Neumann-to-Dirichlet map}\label{L2L2ressubs}

In this subsection we mainly follow \cite{Ron15}.
We consider two alternative ways of computing $V_m$ in \eqref{Vldef}, which are suited to estimate the resistance matrix and its approximation of the Neumann-to-Dirichlet map.

In order to compare the discrete measurements of the resistance matrix with the continuum measurements of the Neumann-to-Dirichlet maps,
we need to introduce the following spaces and projection operators. We call
$$PC:=\left\{\tilde{I}=\sum_{m=1}^M\frac{I_m}{\mathcal{H}^{N-1}(e_m)}\chi_{e_m}:\ I\in \R^M\right\}$$
the set of functions which are constant on each electrode and zero elsewhere. We endow $PC$ with the $L^2(\partial\Omega)$ norm.
It is clear that any $\tilde{I}\in PC$ can be identified with a vector in $\R^M$ through the following linear map $\Phi:\R^M\to PC$ such that
$$\Phi(I)=\tilde{I}=\sum_{m=1}^M\frac{I_m}{\mathcal{H}^{N-1}(e_m)}\chi_{e_m}\quad\text{for any }I\in \R^M.$$
We call
$$\R_{\ast}^M=\left\{I\in \R^M:\ \sum_{m=1}^MI_m=0\right\}$$
which is the set where our discrete current densities are taken from.
We note that $\Phi$ is a bijection from $\R^M$ into $PC$ as well as from $\R_{\ast}^M$ into $PC_{\ast}:=PC\cap L^2_{\ast}(\partial\Omega)$. In particular, for any $I\in\R^M$,
\begin{multline}\label{Itilde}
\|\Phi(I)\|_{L^2(\partial\Omega)}=\|\tilde{I}\|_{L^2(\partial\Omega)}=
\left(\sum_{m=1}^M\frac{I_m^2}{\mathcal{H}^{N-1}(e_m)}\right)^{1/2}\\
\leq {\left(\displaystyle{\min_{m=1,\ldots,M}\{\mathcal{H}^{N-1}(e_m)\}}\right)^{-1/2}}\|I\|.
\end{multline}
and, for any $\tilde{I}\in PC$, since $I_m=\int_{e_m}\tilde{I}$,
\begin{equation}\label{Itildeinverse}
\|I\|=\|\Phi^{-1}(\tilde{I})\|
\leq {\left(\displaystyle{\max_{m=1,\ldots,M}\{\mathcal{H}^{N-1}(e_m)\}}\right)^{1/2}}\|\tilde{I}\|_{L^2(\partial\Omega)}.
\end{equation}

To any $M\times M$ matrix $R$, which can be seen as a linear and bounded operator $R:\R^M\to\R^M$, we can associate the linear and bounded operator $\mathcal{R}\in\mathcal{L}(PC)$ such that
\begin{equation}\label{mathcalR0}
\mathcal{R}(\tilde{I})=    (\Phi   \circ R  \circ\Phi^{-1})(\tilde{I})\quad \text{for any }\tilde{I}\in PC.
\end{equation}
Clearly,
$$
\|R\|\leq \mu^{1/2}\|\mathcal{R}\|_{\mathcal{L}(PC)}\quad\text{and}\quad  \|\mathcal{R}\|_{\mathcal{L}(PC)} \leq \mu^{1/2}\|R\|.
$$
We note that $R\in \mathcal{L}(\R^M_{\ast})$, that is, $R$ is a matrix such that $R[I]\in \R^M_{\ast}$ for any $I\in \R^M_{\ast}$ which is extended as usual to $\R^M$ by setting R[1]=0, if and only if $\mathcal{R}\in \mathcal{L}(PC_{\ast})$. In this case,
\begin{equation}
\|R\|\leq \mu^{1/2}\|\mathcal{R}\|_{\mathcal{L}(PC_{\ast})}\quad\text{and}\quad  \|\mathcal{R}\|_{\mathcal{L}(PC_{\ast})} \leq \mu^{1/2}\|R\|.
\end{equation}

 We define $P_{\ast}:L^2(\partial\Omega)\to L_*^2(\partial\Omega)$,
 $P_e:L^2(\partial\Omega)\to PC$ and $P:L^2(\partial\Omega)\to PC_{\ast}$ 
 as the orthogonal projections on $L_{\ast}^2(\partial\Omega)$, $PC$ and $PC_{\ast}$, respectively. We define $P_{e\ast}:PC\to PC_{\ast}$ as the orthogonal projection on $PC_{\ast}$. All these operators have clearly norm equal to $1$ and they are given by
   \begin{equation}
       P_e[f]=\sum_{m=1}^M \left(\fint_{e_m}f\right)\chi_{e_m}\quad\text{and}\quad P_{\ast}[f]=f-\fint_{\partial\Omega} f
 \quad\text{for any }f\in L^2(\partial\Omega)
   \end{equation}
and
\begin{equation}
      P_{e\ast}[f]=f-\left(\fint_{e} f\right)\chi_e\quad\text{for any }f\in PC,
   \end{equation}
and, finally,  $P=P_{e_*}\circ P_e$.\

Let us note that our voltage measurement $V\in\R^M_{\ast}$ is related to the solution $u$ by
\begin{equation}\label{voltageoperator}
V=\Phi^{-1}(P[u]).
\end{equation}
It follows that
\begin{equation}\label{Vtilde}
\|V\|\leq {\left(\displaystyle{\max_{m=1,\ldots,M}\{\mathcal{H}^{N-1}(e_m)\}}\right)^{1/2}}\| u  \|_{L^2(\partial \Omega)}.
\end{equation}
Therefore, we can alternatively view our measurements either as the resistance matrix, that is, the linear and bounded operator $R_X(A):\R_{\ast}^M\to \R_{\ast}^M$, such that $V=R_X(A)(I)$ for any $I\in \R_{\ast}^M$ and which is extended as usual to $\R^M$, or as the operator $\mathcal{R}_X(A):PC_{\ast}\to PC_{\ast}$ such that
\begin{equation}\label{mathcalR}
\mathcal{R}_X(A)(\tilde{I})=    (\Phi   \circ R_X(A)  \circ\Phi^{-1})(\tilde{I})\quad \text{for any }\tilde{I}\in PC_{\ast}.
\end{equation}
As before,
$$\|R_X(A)\|\leq \mu^{1/2}\|\mathcal{R}_X(A)\|_{\mathcal{L}(PC_{\ast})}\quad\text{and}\quad  \|\mathcal{R}_X(A)\|_{\mathcal{L}(PC_{\ast})} \leq \mu^{1/2}\|R_X(A)\|.$$

Let $\mathcal{K}_X(A):H^{-1/2}_{\ast}(\partial\Omega)\to L^2_{\ast}(\partial\Omega)$
be the
operator defined as follows. For any $\phi\in H^{-1/2}_{\ast}(\partial\Omega)$
$$\mathcal{K}_X(A)[\phi]=\sum_{m=1}^M\frac{1}{z_m}\left(\mathcal{N}_X(A)[\phi]-
\frac{1}{\mathcal{H}^{N-1}(e_m)}\int_{e_m}\mathcal{N}_X(A)[\phi]\right)\chi_{e_m}.$$
We can also consider $\mathcal{K}_X(A)$ to be defined on $(TX)'_{\ast}$. In both cases,
$\mathcal{K}_X(A)$ is a compact linear operator. Moreover,
$$\|\mathcal{K}_X(A)[\phi]\|^2_{L^2_{\ast}(\partial\Omega)}\leq \frac1{Z_1^2}\| \mathcal{N}_X(A)[\phi]  \|^2_{L^2(e)}
\leq \frac{\hat{C}_P^2}{Z_1^2}\frac{1}{\lambda_0^2}  \|\phi\|^2_{(TX)'_{\ast}},
$$
thus
\begin{equation}\label{normk}
\|\mathcal{K}_X(A) \|_{\mathcal{L}(H^{-1/2}_{\ast}(\partial\Omega),L^2_{\ast}(\partial\Omega))}\leq 
\|\mathcal{K}_X(A) \|_{\mathcal{L}((TX)'_{\ast},L^2_{\ast}(\partial\Omega))}\leq 
\frac{\hat{C}_P}{Z_1}\frac{1}{\lambda_0}.
\end{equation} 
Moreover, $\mathcal{K}_X(A)$ is a compact linear operator also from $H^{-1/2}_{\ast}(\partial\Omega)$
into itself, from $L^2_{\ast}(\partial\Omega)$ into itself and from $(TX)'_{\ast}$ into itself.
Since for any $g\in H^{-1/2}_{\ast}(\partial\Omega)$ and $I=0$
the equation \eqref{eqform} admits a solution, we can infer that
$\mathbbm{1}+\mathcal{K}_X(A):(TX)'_{\ast}\to (TX)'_{\ast}$
is bijective.
We deduce that 
$\mathbbm{1}+\mathcal{K}_X(A)$ is bijective as well from $H^{-1/2}_{\ast}(\partial\Omega)$
into itself and from $L^2_{\ast}(\partial\Omega)$ into itself.
We denote with $\tilde{\mathcal{K}}_X(A)$ the inverse to $\mathbbm{1}+\mathcal{K}_X(A)$.
For any $g\in (TX)'_{\ast}$, or $g\in H_{\ast}^{-1/2}(\partial\Omega)$, let $\phi=\tilde{\mathcal{K}}_X(A)[g]$. Then, by \eqref{boundaryeqform} and \eqref{coercbis}, we infer that
$$\|\phi\|_{(TX)'_{\ast}}\leq \|\phi\|_{H_{\ast}^{-1/2}(\partial\Omega)}   \leq  \frac{\lambda_1}{\lambda_0}\| g\|_{(TX)'_{\ast}}\leq \frac{\lambda_1}{\lambda_0}\| g\|_{H_{\ast}^{-1/2}(\partial\Omega)}.$$

We conclude that, if $g\in L^2_{\ast}(\partial\Omega)$, then
$$\|\tilde{\mathcal{K}}_X(A)[g]\|_{H^{-1/2}_{\ast}(\partial\Omega)}\leq C_P\frac{\lambda_1}{\lambda_0}\|g\|_{L^2_{\ast}(\partial\Omega)},
$$
hence, since $\tilde{\mathcal{K}}_X(A)[g]=g-\mathcal{K}_X(A)\big[\tilde{\mathcal{K}}_X(A)[g]\big]$,
\begin{multline}\label{L2L2}
\|\tilde{\mathcal{K}}_X(A)[g]\|_{L^2_{\ast}(\partial\Omega)}\leq\\
\|\mathcal{K}_X(A) \big[\tilde{\mathcal{K}}_X(A)[g]\big] \|_{L^2_{\ast}(\partial\Omega)}+
\|g\|_{L^2_{\ast}(\partial\Omega)}
\leq \left(
\frac{\hat{C}_PC_P}{Z_1}\frac{\lambda_1}{\lambda_0^2}+1\right)\|g\|_{L^2_{\ast}(\partial\Omega)}.
\end{multline}
We call
\begin{equation}\label{tildeC1}
\tilde{C}_1=\left(
\frac{\hat{C}_PC_P}{Z_1}\frac{\lambda_1}{\lambda_0^2}+1\right)
\end{equation}
and we notice that it depends $r$, $R$, $L$, $\lambda_0$, $\lambda_1$ and $Z_1$ only.

For any given current pattern $I$, taking $g=0$,
it is immediate to show that $\phi=\tilde{\mathcal{K}}_X(A)[\tilde{I}]$, hence
$$u=\mathcal{N}_X(A)\big[\tilde{\mathcal{K}}_X(A)[\tilde{I}]\big].$$
Therefore we can alternatively define
\begin{equation}\label{R-Kconn}
V_m=\int_{e_m}\mathcal{N}_X(A)\big[\tilde{\mathcal{K}}_X(A)[\tilde{I}]\big]+c\mathcal{H}^{N-1}(e_m),
\end{equation}
with
\begin{equation}\label{constantc}
c=-\frac{\sum_{m=1}^M\int_{e_m}\mathcal{N}_X(A)\big[\tilde{\mathcal{K}}_X(A)[\tilde{I}]\big]}{\sum_{m=1}^M\mathcal{H}^{N-1}(e_m)}.
\end{equation}

Still considering $g=0$, and taking $W\equiv 0$ in \eqref{eqform}, we
also infer that, in $(TX)'$,
$$\phi=\sum_{m=1}^M\frac{1}{z_m}(U_m-u)\chi_{e_m},$$
hence
\begin{equation}\label{R-Kconnbis}
V_m=\int_{e_m}\mathcal{N}_X(A)\left[\sum_{m=1}^M\frac{1}{z_m}(U_m-u)\chi_{e_m}\right]+c\mathcal{H}^{N-1}(e_m),
\end{equation}
with
\begin{equation}\label{constantcbis}
c=-\frac{\sum_{m=1}^M\int_{e_m}\mathcal{N}_X(A)\left[\sum_{m=1}^M\frac{1}{z_m}(U_m-u)\chi_{e_m}\right]}{\sum_{m=1}^M\mathcal{H}^{N-1}(e_m)}.
\end{equation}
We finally observe that, by \eqref{Uexp},
\begin{equation}\label{Uexp2}
I_m=\frac1{z_m}\int_{e_m}(U_m-u).
\end{equation}

The first characterization allows us to infer the next result, the second is useful in the next subsection.

\begin{prop}\label{expmeasprop}
Let $A_1$, $A_2\in \mathcal{M}(\lambda_0,\lambda_1)$.
Then there exists a constant $\hat{c}$
such that
\begin{equation}\label{expmeasest}
\|R_X(A_1)-R_X(A_2)\|\leq \hat{c}\|\mathcal{N}_X(A_1)-\mathcal{N}_X(A_2)\|_{\mathcal{L}(L^2_{\ast}(\partial\Omega),L^2_{\ast}(\partial\Omega))}.
\end{equation}
The constant $\hat{c}$ has the following form
\begin{equation}
\hat{c}=\left(\frac{\hat{C}_P^2\tilde{C}_1^2}{Z_1\lambda_0}+\tilde{C}_1
\right)\mu^{1/2}
\end{equation}
where $\mu$ is defined in \eqref{mudef}, and $\tilde{C}_1$ is as in \eqref{tildeC1}, thus it depends on $r$, $L$, $R$, $\lambda_0$, $\lambda_1$ and $Z_1$ only.
\end{prop}

\proof{.} We recall that we have set $R_X(A_1)[1]=R_X(A_2)[1]=0$.
We evaluate
$\|(R_X(A_1)-R_X(A_2))I\|$ for any $I\in\mathbb{R}^M$, with
$\sum_{m=1}^MI_m=0$.
Arguing as in the proof of \eqref{Vtilde} and using
\eqref{R-Kconn},
we have that
\begin{multline*}
\|(R_X(A_1)-R_X(A_2))I\|\\\leq C\|
\mathcal{N}_X(A_1)\big[\tilde{\mathcal{K}}_X(A_1)[\tilde{I}]\big]-\mathcal{N}_X(A_2)\big[\tilde{\mathcal{K}}_X(A_2)[\tilde{I}]\big]
\|_{L^2_{\ast}(\partial\Omega)},
\end{multline*}
where $C={\left(\displaystyle{\max_{m=1,\ldots,M}\{\mathcal{H}^{N-1}(e_m)\}}\right)^{1/2}}$.
But
\begin{multline}\label{Rest}
\|\mathcal{N}_X(A_1)\big[\tilde{\mathcal{K}}_X(A_1)[\tilde{I}]\big]-\mathcal{N}_X(A_2)\big[\tilde{\mathcal{K}}_X(A_2)[\tilde{I}]\big]
\|_{L^2_{\ast}(\partial\Omega)}\\
 \leq\|\mathcal{N}_X(A_1)[\big(\tilde{\mathcal{K}}_X(A_1)-\tilde{\mathcal{K}}_X(A_2))[
\tilde{I}]\big]\|_{L^2_{\ast}(\partial\Omega)}+\\
\|(\mathcal{N}_X(A_1)-\mathcal{N}_X(A_2))\big[\tilde{\mathcal{K}}_X(A_2)
[\tilde{I}]\big]\|_{L^2_{\ast}(\partial\Omega)} 
\\
 \leq \Big(\|\mathcal{N}_X(A_1)\|_{\mathcal{L}(L^2_{\ast}(\partial\Omega))}
\|\tilde{\mathcal{K}}_X(A_1)-\tilde{\mathcal{K}}_X(A_2)\|_{\mathcal{L}(L^2_{\ast}(\partial\Omega))}
+\\
 \| \tilde{\mathcal{K}}_X(A_2) \|_{\mathcal{L}(L^2_{\ast}(\partial\Omega))}  \|\mathcal{N}_X(A_1)-\mathcal{N}_X(A_2)\|_{\mathcal{L}(L^2_{\ast}(\partial\Omega))}\Big)
 \|\tilde{I}\|_{L^2_{\ast}(\partial\Omega)}\\
 \leq
 \Big(\frac{\hat{C}_P^2}{\lambda_0}
\|\tilde{\mathcal{K}}_X(A_1)-\tilde{\mathcal{K}}_X(A_2)\|_{\mathcal{L}(L^2_{\ast}(\partial\Omega))}
+\\
\tilde{C}_1 \|\mathcal{N}_X(A_1)-\mathcal{N}_X(A_2)\|_{\mathcal{L}(L^2_{\ast}(\partial\Omega))}\Big)
 {\left(\displaystyle{\min_{m=1,\ldots,M}\{\mathcal{H}^{N-1}(e_m)\}}\right)^{-1/2}}\|I\|,
\end{multline}
where we used \eqref{L2L2N}, \eqref{L2L2} and \eqref{Itilde}.

It remains to evaluate the term
$\|\tilde{\mathcal{K}}_X(A_1)-\tilde{\mathcal{K}}_X(A_2)\|_{\mathcal{L}(L^2_{\ast}(\partial\Omega))}$.
Using the identity
$$\tilde{\mathcal{K}}_X(A_1)-\tilde{\mathcal{K}}_X(A_2)=\tilde{\mathcal{K}}_X(A_1)[\mathcal{K}_X(A_2)-\mathcal{K}_X(A_1)]\tilde{K}_X(A_2)
$$
we obtain that, using again \eqref{L2L2},
\begin{equation}
\label{K1-K2}
\|\tilde{\mathcal{K}}_X(A_1)-\tilde{\mathcal{K}}_X(A_2)\|_{\mathcal{L}(L^2_{\ast}(\partial\Omega))}\leq \tilde{C}_1^2
\|\mathcal{K}_X(A_1)-\mathcal{K}_X(A_2)\|_{\mathcal{L}(L^2_{\ast}(\partial\Omega))}.
\end{equation}
Since
$$\|\mathcal{K}_X(A_1)-\mathcal{K}_X(A_2)\|_{\mathcal{L}(L^2_{\ast}(\partial\Omega))}\leq \frac{1}{Z_1}
\|\mathcal{N}_X(A_1)-\mathcal{N}_X(A_2)\|_{\mathcal{L}(L^2_{\ast}(\partial\Omega))},
$$
the proof is concluded.\cvd

\bigskip

\subsection{Approximating the Neumann-to-Dirichlet map by the resistance matrices}\label{approxNtoDsubs}
 In this subsection we mainly follow \cite{Hyv1,Hyv2,Hyv-et-al} and we need to consider the following further assumptions on the electrodes. Let us assume that to each electrode $e_m$, $m=1,\ldots,M$, is associated the \emph{extended electrode}
$\tilde{e}_m$, an open subset of $\partial\Omega$, such that $e_m\subset\tilde{e}_m$, the sets $\tilde{e}_m$ are pairwise disjoint and their union covers $\partial\Omega$ up to a negligible set, that is,
$$\mathcal{H}^{N-1}\left(\partial\Omega\backslash\bigcup_{m=1}^M\tilde{e}_m\right)=0.$$
Further assumptions on these sets $\tilde{e}_m$ will be made when needed.
Analogously, we call
$$\widetilde{PC}:=\left\{\tilde{I}=\sum_{m=1}^M\frac{I_m}{\mathcal{H}^{N-1}(\tilde{e}_m)}\chi_{\tilde{e}_m}:\ I\in \R^M\right\}$$
the set of functions which are constant on each extended electrode and zero elsewhere. We still endow $\widetilde{PC}$ with the $L^2(\partial\Omega)$ norm and call $\widetilde{PC}_{\ast}=\widetilde{PC}\cap L^2_{\ast}(\partial\Omega)$.

The following geometric quantities are of interest
\begin{equation}\label{thetadef}
\theta:=\max_{m=1,\ldots,M}\left\{\frac{\mathcal{H}^{N-1}(\tilde{e}_m)}{\mathcal{H}^{N-1}(e_m)}\right\}
\end{equation}
and
\begin{equation}\label{etadef}
\eta:=\max_{m=1,\ldots,M}\left\{\frac{\mathrm{diam}(\tilde{e}_m)}{(\mathcal{H}^{N-1}(\tilde{e}_m))^{1/(N-1)}}\right\},
\end{equation}
that is,
$$\mathrm{diam}(\tilde{e}_m)\leq\eta(\mathcal{H}^{N-1}(\tilde{e}_m))^{1/(N-1)}\quad\text{for any }m=1,\ldots,M.$$
It also implies that
$$\mathrm{diam}(e_m)\leq \mathrm{diam}(\tilde{e}_m)\leq\eta\theta^{1/(N-1)}(\mathcal{H}^{N-1}(e_m))^{1/(N-1)}\ \text{for any }m=1,\ldots,M.$$
We note that the reverse inequality always holds, that is, there exists a constant $\tilde{H}$, depending on $r$, $L$ and $R$ only, such that for any $E$ Borel subset of $\partial\Omega$ we have
\begin{equation}\label{measurevsdiam}
\mathcal{H}^{N-1}(E)\leq \tilde{H}(\mathrm{diam}(E))^{N-1}.
\end{equation}

These quantities also provide an upper bounds on the number of electrodes $M$. Namely, let $m_0$ be such that
$\delta=\mathrm{diam}(e_{m_0})$. Then
$$\delta\leq\eta\theta^{1/(N-1)}(\mathcal{H}^{N-1}(e_{m_0}))^{1/(N-1)}\leq \eta\left(\theta\mu \min_{m=1,\ldots,M}\{\mathcal{H}^{N-1}(e_m)\}
\right)^{1/(N-1)},$$
so we have
\begin{multline*}
M\delta\leq \eta(\theta\mu)^{1/(N-1)}\sum_{m=1}^M(\mathcal{H}^{N-1}(e_m))^{1/(N-1)}\\ \leq
\eta(\theta\mu)^{1/(N-1)}MM^{-1/(N-1)}\left(\sum_{m=1}^M\mathcal{H}^{N-1}(e_m)\right)^{1/(N-1)}\\\leq
\eta(\theta\mu)^{1/(N-1)}MM^{-1/(N-1)}\left(\mathcal{H}^{N-1}(\partial\Omega)\right)^{1/(N-1)}.
\end{multline*}
We conclude that
\begin{equation}\label{Lest}
M\leq \eta^{N-1}\theta\mu\mathcal{H}^{N-1}(\partial\Omega)\delta^{-(N-1)}\leq
\eta^{N-1}\theta\mu\tilde{c}_2\delta^{-(N-1)}
\end{equation}
where $\tilde{c}_2$ as in \eqref{measure} depends on $r$, $L$ and $R$ only.

In order to approximate the continuum Neumann-to-Dirichlet map we need the following further operators.

We define the non-orthogonal projection $Q:L^2(\partial\Omega)\to L^2(\partial\Omega)$ such that
$$Q(f):= \sum_{m=1}^M
      \frac{1}{\mathcal{H}^{N-1}(e_m)}\left( \int_{\tilde{e}_m}f \right)\chi_{e_m}\quad\text{for any }f\in L^2(\partial\Omega).$$
It is easy to check that $Q(f)\in PC$, for any $f\in L^2(\Omega)$, and that 
\begin{equation}\label{Qnorm}
\|Q\|_{\mathcal{L}(L^2(\partial\Omega))}\leq \theta^{1/2}.
\end{equation}
Moreover,
$Q:L^2_{\ast}(\partial\Omega)\to PC_{\ast}$ and $Q$ is the identity on $PC_{\ast}$.
Through $Q$, any $L^2_{\ast}(\partial\Omega)$ current density can be mapped into a current density belonging to $PC_{\ast}$.

For any voltage $u$ on the boundary, what is measured is $\Phi^{-1}(P[u])$, which is a vector in $\R^M_{\ast}$, or $P[u]\in PC_{\ast}$.
We define the extension $E:PC\to L^2(\partial\Omega)$ such that
$$E(\tilde{V})=P_{\ast}\left[\sum_{m=1}^M\frac{V_m}{\mathcal{H}^{N-1}(e_m)}\chi_{\tilde{e}_m}\right]  \quad\text{for any }\tilde{V}\in PC$$
where $V=\Phi^{-1}(\tilde{V})$.
Clearly $E(\tilde{V})\in L^2_{\ast}(\partial\Omega)$ and $E(\tilde{V})\in \widetilde{PC}_{\ast}$, for any $\tilde{V}\in PC$.
Moreover, we can decompose $E$ as $E=P_{\ast}\circ \tilde{E}\circ \Phi^{-1}$ where
$$\tilde{E}(V):=\sum_{m=1}^M\frac{V_m}{\mathcal{H}^{N-1}(e_m)}\chi_{\tilde{e}_m}  \quad\text{for any }V\in \R^M.$$
We have that
$$\|\tilde{E}(V)\|_{L^2(\partial\Omega)}\leq {\left(\displaystyle{\min_{m=1,\ldots,M}\{\mathcal{H}^{N-1}(e_m)\}}\right)^{-1/2}} \theta^{1/2}  \|V\|\quad\text{for any }V\in \R^M,$$
but it is not difficult to show that
 \begin{equation}\label{Enorm}
 \|E\|_{\mathcal{L}(PC,L^2_{\ast}(\partial\Omega))}=\|\tilde{E}\circ \Phi^{-1}\|_{\mathcal{L}(PC,L^2(\partial\Omega))}  \leq \theta^{1/2}.
 \end{equation}
The operator $E|_{PC_{\ast}}:PC_{\ast}\to \widetilde{PC}_{\ast}$ has an inverse, which is given by $E^{-1}:\widetilde{PC}_{\ast}\to PC_{\ast}$
 such that
 \begin{equation}
      E^{-1}[f] := P_{e*}\left[ \sum_{m=1}^M V_m\chi_{e_m} \right]\quad\text{for any }f=\sum_{m=1}^M V_m\chi_{\tilde{e}_m}\in \widetilde{PC}_{\ast}.
   \end{equation}
 This depends on the fact that
 $E(c\chi_e)=0$ for any $c\in \R$.
 It is easy to see that $\| E^{-1} \|_{\mathcal{L}_{\ast}(\widetilde{PC}_{\ast},PC_{\ast})}\leq 1$.
 
 Therefore, for any operator $\mathcal{R}\in \mathcal{L}(PC_{\ast})$, we have that
 $$\| \mathcal{R}\|_{\mathcal{L}(PC_{\ast})}  \leq\|E\circ\mathcal{R}\circ Q\|_{\mathcal{L}(L^2_{\ast}(\partial\Omega))} \leq \theta \|\mathcal{R}\|_{\mathcal{L}(PC_{\ast})} $$
 The second inequality follows from \eqref{Qnorm} and \eqref{Enorm}, whereas the first follows from the fact that for any $\tilde{I}\in PC_{\ast}$ we have
$$\mathcal{R}[\tilde{I}]=\mathcal{R}\big[Q[\tilde{I}]\big]=(E^{-1}\circ E)\left[ \mathcal{R}\big[Q[\tilde{I}]\big]  \right].$$
 
Finally, the following relationship between $P$, $E$ and $Q$ hold.
   \begin{lem}
   \label{sec4-sub3-lemma1}
  We have that $E\circ P\in\mathcal{L}(L_{\ast}^2(\partial\Omega))$ is the adjoint of $Q\in\mathcal{L}(L_{\ast}^2(\partial\Omega))$. In particular, it follows that
      \begin{equation}\label{adjoint}
         \|\mathbbm{1}-E\circ P\|_{\mathcal{L}(H_{\ast}^{1/2}(\partial\Omega),L_{\ast}^2(\partial\Omega))} =
         \|\mathbbm{1}-Q\|_{\mathcal{L}(L_{\ast}^2(\partial\Omega),H_{\ast}^{-1/2}(\partial\Omega))}.
      \end{equation}
   \end{lem}
 \proof{.}
      Let $f,\,g\in L_{\ast}^2(\partial\Omega)$. Then we have the following equality
      \begin{multline*}
         \langle (E \circ P)[f],g\rangle\\=
         \int_{\partial\Omega}\left(
         \sum_{m=1}^M \left( \frac{1}{\mathcal{H}^{N-1}(e_m)} \int_{e_m}f(x)d\mathcal{H}^{N-1}(x) - \fint_ef \right)\chi_{\tilde{e}_m}(y)\right)g(y)d\mathcal{H}^{N-1}(y) \\
         = \sum_{m=1}^M\left( \int_{\partial\Omega} \left(\int_{\partial\Omega} \frac{1}{\mathcal{H}^{N-1}(e_m)} f(x)\chi_{e_m}(x)\chi_{\tilde{e}_m}(y)g(y)d\mathcal{H}^{N-1}(x)\right)dy\right) \\-
         \int_{\partial\Omega} \left( \fint_e f \right)g(y)d\mathcal{H}^{N-1}(y) \\
         = \int_{\partial\Omega} f(x)\left(\sum_{m=1}^M \left( \frac{1}{\mathcal{H}^{N-1}(e_m)}\int_{\tilde{e}_m}g(y)d\mathcal{H}^{N-1}(y)\right)\chi_{e_m}(x)\right)dx=
         \langle f,Q[g]\rangle.
      \end{multline*}
      The thesis follows.\cvd
      
   \bigskip

We now want to estimate $\|\mathbbm{1}-E\circ P\|_{\mathcal{L}(H_{\ast}^{1/2}(\partial\Omega),L_{\ast}^2(\partial\Omega))}$.

 \begin{lem}\label{extensionestlemma}
 Under the previous assumptions on the electrodes, the following inequality holds
      \begin{equation}
      \label{sec4-sub3-eq6}
         \|\mathbbm{1}-E\circ P\|_{\mathcal{L}(H_{\ast}^{1/2}(\partial\Omega),L_{\ast}^2(\partial\Omega))} \leq \tilde{C}_2\delta^{1/2},
         \end{equation}
where 
\begin{equation}\label{tildeC2}
\tilde{C}_2= C_T\left(\eta^N\theta^{N/(N-1)}\tilde{H}^{1/(N-1)}\right)^{1/2},
\end{equation}
with $\tilde{H}$ is as in \eqref{measurevsdiam}.
Moreover, for the same constant $\tilde{C}_2$,
      \begin{equation}
      \label{sec4-sub3-eq6bis}
         \|\mathbbm{1}-Q\|_{\mathcal{L}(L_{\ast}^2(\partial\Omega),H_{\ast}^{-1/2}(\partial\Omega))} \leq \tilde{C}_2\delta^{1/2}.
         \end{equation} 
   \end{lem}
   \proof{.} Let $f\in H_{\ast}^{1/2}(\partial\Omega)$. Then we have that
   \begin{multline*}
   \|(\mathbbm{1}-E\circ P)[f]\|_{L^2_{\ast}(\partial\Omega)}\\\leq
   \sum_{m=1}^M\int_{\tilde{e}_m}\left(f(x)-\frac{1}{\mathcal{H}^{N-1}(e_m)}\int_{e_m}f(y)d\mathcal{H}^{N-1}(y)   \right)^2d\mathcal{H}^{N-1}(x)\\=
   \sum_{m=1}^M\int_{\tilde{e}_m}\left(\frac{1}{\mathcal{H}^{N-1}(e_m)}\int_{e_m}(f(x)-f(y))d\mathcal{H}^{N-1}(y)   \right)^2d\mathcal{H}^{N-1}(x)\\
   \leq  \sum_{m=1}^M\frac{1}{\mathcal{H}^{N-1}(e_m)}\int_{\tilde{e}_m}\left(\int_{e_m}|f(x)-f(y)|^2d\mathcal{H}^{N-1}(y)   \right)d\mathcal{H}^{N-1}(x)
   \\\leq \sum_{m=1}^M\frac{(\mathrm{diam}(\tilde{e}_m))^{N}}{\mathcal{H}^{N-1}(e_m)}
   \int_{\tilde{e}_m}\left(\int_{\tilde{e}_m}\frac{|f(x)-f(y)|^2}{|x-y|^{N}}d\mathcal{H}^{N-1}(x)\right)d\mathcal{H}^{N-1}(y)
   \\
   \leq \sum_{m=1}^M\eta^N(\theta^N\mathcal{H}^{N-1}(e_m))^{1/(N-1)}\!\int_{\tilde{e}_m}\!\!\left(\int_{\tilde{e}_m}\!\frac{|f(x)-f(y)|^2}{|x-y|^{N}}d\mathcal{H}^{N-1}(x)\right)d\mathcal{H}^{N-1}(y)
   \\
   \leq \eta^N\theta^{N/(N-1)}\tilde{H}^{1/(N-1)}\delta  |f|^2_{B^{1/2,1/2}(\partial\Omega)}.
     \end{multline*}
Inequality \eqref{sec4-sub3-eq6} is proved, whereas \eqref{sec4-sub3-eq6bis} immediately follows from \eqref{adjoint}.\cvd
   
\bigskip

The approximation result is the following.

\begin{teo}
   \label{sec4-sub3-theorem1}
   Let $\Omega\in\mathcal{A}(r,L,R)$ and let $A\in \mathcal{M}(\lambda_0,\lambda_1)$. Let $X$ be a closed subspace of $H^1(\Omega)$ containing constants and let $\mathcal{R}_X(A)$ be defined as in \eqref{mathcalR}.
      Under the previous assumptions on the electrodes, the following inequality holds
      \begin{equation}
         \|\mathcal{N}_X(A)-E\circ \mathcal{R}_X(A)\circ Q\|_{\mathcal{L}(L_{\ast}^2(\partial\Omega))} \leq \tilde{C}_3\delta^{1/2}
      \end{equation}
      where
      \begin{equation}\label{tildeC3}
      \tilde{C}_3=\frac{\hat{C}_P}{\lambda_0}\tilde{C}_2 \left(1+  2\tilde{C}_1 \theta^{1/2}\right).
      \end{equation}
      Here $\tilde{C}_1$ is given in \eqref{tildeC1} and $\tilde{C}_2$ is given in \eqref{tildeC2}.
        \end{teo}

   \proof{.}
   Let us fix 
      $f\in L_{\ast}^2(\partial\Omega)$ and let $v=v(F_f)\in X_{\ast}$ solve \eqref{Neuweak} with $F$ replaced by $F_f$. Let $(u,U)\in H_{\ast}$ be the solution to \eqref{eqform} with $g=0$ and $I=\Phi ( Q[f])$. Then
\begin{multline*}
\|\mathcal{N}_X(A)[f]-(E\circ \mathcal{R}_X(A)\circ Q)[f]\|_{L_{\ast}^2(\partial\Omega)}\leq
\|v-(E\circ P)[u]\|_{L_{\ast}^2(\partial\Omega)}
\\
\leq \|v-u\|_{L_{\ast}^2(\partial\Omega)}
+\|(\mathbbm{1}-E\circ P)[u]\|_{L_{\ast}^2(\partial\Omega)}.
\end{multline*}
    
 Recalling that $u=\mathcal{N}_X(A)\big[\tilde{\mathcal{K}}_X(A)[Q[f]]\big]$,   
the second term can be estimated as follows
$$\|(\mathbbm{1}-E\circ P)[u]\|_{L_{\ast}^2(\partial\Omega)}=\|(\mathbbm{1}-E\circ P)[u]\|_{L_{\ast}^2(\partial\Omega)}
\leq \frac{\hat{C}_P}{\lambda_0}\tilde{C}_2\delta^{1/2}   \tilde{C}_1 \theta^{1/2}\|f\|_{L^2(\partial\Omega)}.$$
Here we have used \eqref{L2L2N0}, \eqref{L2L2}, thus $\tilde{C}_1$ is as in \eqref{tildeC1}, \eqref{Qnorm}, and, finally, \eqref{sec4-sub3-eq6},
thus $\tilde{C}_2$ is as in \eqref{tildeC2}.
      
For the first term, we use the other representation of $u$, namely,
$$u=\mathcal{N}_X(A)\left[\sum_{m=1}^M\frac{1}{z_m}(U_m-u)\chi_{e_m}\right],$$
that is,
$$\phi=\left[\sum_{m=1}^M\frac{1}{z_m}(U_m-u)\chi_{e_m}\right]=\tilde{\mathcal{K}}_X(A)[Q[f]].
 $$
Therefore, again by \eqref{L2L2N0},
\begin{multline*}
\|v-u\|_{L_{\ast}^2(\partial\Omega)}\leq \frac{\hat{C}_P}{\lambda_0}\left\|f-  \sum_{m=1}^M\frac{1}{z_m}(U_m-u)\chi_{e_m}  \right\|_{H_{\ast}^{-1/2}(\partial\Omega)}\\
\leq \frac{\hat{C}_P}{\lambda_0}\left(\|f- Q[f] \|_{H_{\ast}^{-1/2}(\partial\Omega)}
+\left\|Q[f]-  \sum_{m=1}^M\frac{1}{z_m}(U_m-u)\chi_{e_m}  \right\|
_{H_{\ast}^{-1/2}(\partial\Omega)}\right)\\
\leq \frac{\hat{C}_P}{\lambda_0}\tilde{C}_2\delta^{1/2}\|f\|_{L^2_{\ast}(\partial\Omega)}+\frac{\hat{C}_P}{\lambda_0}\left\|Q[f]-  \sum_{m=1}^M\frac{1}{z_m}(U_m-u)\chi_{e_m}  \right\|_{H_{\ast}^{-1/2}(\partial\Omega)}.
\end{multline*}
For the last inequality we used \eqref{sec4-sub3-eq6}. We need to treat the last term.
By \eqref{Uexp2}, we have that
$$Q[f]=Q\left[\sum_{m=1}^M\frac{1}{z_m}(U_m-u)\chi_{e_m}\right],$$
so, again by \eqref{sec4-sub3-eq6}, \eqref{L2L2} and \eqref{Qnorm},
\begin{multline*}
\left\|Q[f]-  \sum_{m=1}^M\frac{1}{z_m}(U_m-u)\chi_{e_m}  \right\|_{H_{\ast}^{-1/2}(\partial\Omega)}\\\leq
\tilde{C}_2\delta^{1/2}\left\|\sum_{m=1}^M\frac{1}{z_m}(U_m-u)\chi_{e_m}  \right\|_{L_{\ast}^2(\partial\Omega)}
=  \tilde{C}_2\delta^{1/2}\left\| \tilde{\mathcal{K}}_X(A)[Q[f]]\right\|_{L_{\ast}^2(\partial\Omega)}
\\\leq \tilde{C}_2\delta^{1/2}\tilde{C}_1\theta^{1/2}  \|f\|_{L^2_{\ast}(\partial\Omega)}.
\end{multline*}
The proof is concluded.\cvd

\bigskip

\subsection{The simplified resistance matrix}\label{simplsec}

We define, for any conductivity tensor $A\in \mathcal{M}(\lambda_0,\lambda_1)$ and any 
subspace $X$ containing constants, a new $M\times M$ matrix $\hat{R}_X(A)$, which we call the \emph{simplified resistance matrix}.
For any $I\in \R^M_{\ast}$, we define
\begin{equation}\label{simplresmat}
\hat{R}_X(A)I= (\Phi^{-1}\circ P\circ \mathcal{N}_X(A)\circ \Phi)[I].
\end{equation}
We assume that $\hat{R}_X(A)[1]=0$.
This corresponds to the following experiment. To any current density $I\in \R^M_{\ast}$, we consider the solution
$v=v(F_{\tilde{I}})\in X_{\ast}$ solving \eqref{Neuweak} with $F$ replaced by $F_{\tilde{I}}$, $\tilde{I}$ being as usual $\Phi^{-1}(I)$. We consider
$$V_m=\int_{e_m}v+c\mathcal{H}^{N-1}(e_m)$$
where $c$ is chosen in such a way that $V\in \R^M_{\ast}$.

Correspondingly, we define the operator
$\hat{\mathcal{R}}_X(A):PC_{\ast}\to PC_{\ast}$ such that
\begin{equation}\label{mathcalhatR}
\hat{\mathcal{R}}_X(A)(\tilde{I})=    \Phi   \circ \hat{R}_X(A)  \circ\Phi^{-1}= P\circ \mathcal{N}_X(A) \quad \text{for any }\tilde{I}\in PC_{\ast}.
\end{equation}

We have the following result.

\begin{prop}\label{simplvsnot} Under the previous assumptions on the electrodes, we have
\begin{equation}
\|\hat{\mathcal{R}}_X(A)-\mathcal{R}_X(A)\|_{\mathcal{L}(PC_{\ast})}\leq \tilde{C}_4\delta^{1/2},
\end{equation}
where
\begin{equation}\label{tildeC4}
\tilde{C}_4=\frac{\hat{C}_P}{\lambda_0}\tilde{C}_2 \left(1+  (2\tilde{C}_1+1) \theta^{1/2}\right).
\end{equation}

\end{prop}

\proof{.} We begin by noting that
\begin{multline*}
\|\hat{\mathcal{R}}_X(A)-\mathcal{R}_X(A)\|_{\mathcal{L}(PC_{\ast})}\leq
\|E\circ (\hat{\mathcal{R}}_X(A)-\mathcal{R}_X(A))\circ Q\|_{\mathcal{L}(PC_{\ast},L^2_{\ast}(\partial\Omega))}\\\leq
\|E\circ (P\circ \mathcal{N}_X(A))\circ Q-\mathcal{N}_X(A)\circ Q\|_{\mathcal{L}(L^2_{\ast}(\partial\Omega))}\\+
\|\mathcal{N}_X(A)-E\circ\mathcal{R}_X(A)\circ Q|\|_{\mathcal{L}(L^2_{\ast}(\partial\Omega))}.
\end{multline*}
The second term is estimated in Theorem~\ref{sec4-sub3-theorem1}, the first one, by \eqref{Qnorm}, \eqref{L2L2N0} and \eqref{sec4-sub3-eq6},
 is bounded by
$\tilde{C}_2\delta^{1/2}\dfrac{\hat{C}_P}{\lambda_0}\theta^{1/2}$. The thesis follows.\cvd

\bigskip

We note that the analogous estimate holds
\begin{equation}
\|(\hat{\mathcal{R}}_X(A)-\mathcal{R}_X(A))\circ Q\|_{\mathcal{L}(L^2_{\ast}(\partial\Omega)}\leq \tilde{C}_4\delta^{1/2},
\end{equation}
$\tilde{C}_4$ as in \eqref{tildeC4}.
 
\section{The discretization estimate}\label{discrsec}

We consider the discretization of the conductivities and of the solutions to the direct problem. For conductivities
we follow mainly the arguments of \cite{Ron16}, see Proposition~\ref{Adiscrprop}. The analysis for solutions to the direct problem is new and it is actually much more involved, see Proposition~\ref{vdiscrprop}.

Let $\Omega\in \mathcal{A}(r,L,R)$ be such that $\Omega$ is 
discretizable with constant $s>0$, as in Definition~\ref{sdiscretizable}. For any $h$, $0<h\leq 1$, let $\mathcal{T}_h$, $X^h$
and $\Pi_h$ as in Definition~\ref{sdiscretizable}. We call $X^h[\mathbb{M}]$ any $N\times N$ matrix valued function whose elements all belong to $X^h$. We also note that $\Pi_h$ is assumed to operate element by element in the case of matrix or vector valued functions. We point out that symmetric matrices are mapped onto symmetric matrices by $\Pi_h$ and that $\Pi_h$ preserves the ellipticity conditions, with the same constants, as well.

Let  $\eta$ be a positive symmetric mollifier,
that is,
$$\eta\in C^{\infty}_0(B_1),\quad \eta\geq 0, \quad \displaystyle \int_{B_1}\eta=1,$$ and such that $\eta(x)$ depends only on $\|x\|$ for any $x\in B_1$. Clearly $\eta\in C^{\infty}_0(\mathbb{R}^N)$ by extending it to $0$ oustide $B_1$.
For any $\gamma>0$, we call
$$\eta_{\gamma}(x)=\gamma^{-N}\eta(x/\gamma)\quad x\in\mathbb{R}^N,$$
and, for any $u\in L^1_{loc}(\mathbb{R}^N)$, we call
$$u_{\gamma}=\eta_{\gamma}\ast u,$$
where as usual $\ast$ denotes the convolution.
If $A\in L^1_{loc}(\mathbb{R}^N,\mathbb{M}^{N\times N}(\R))$, we call
$$A_{\gamma}=\eta_{\gamma}\ast A,$$
where the convolution is done element by element. 
We note that if $A$ is symmetric, then $A_{\gamma}$ is still symmetric.
 Moreover,
if $\lambda_0\leq u\leq\lambda_1$, we have that still $\lambda_0\leq u_{\gamma}\leq\lambda_1$, and, if
$\lambda_0 I_N\leq A\leq\lambda_1 I_N$, we have that still $\lambda_0 I_N\leq A_{\gamma}\leq\lambda_1I_N$.

We know that $u_\gamma \in C^{\infty}(\R^N)$ and, for any multiindex $\alpha$, we have
$D^{\alpha}u_{\gamma}=(D^{\alpha}\eta_{\gamma})\ast u$. Clearly we have that $A_\gamma\in C^{\infty}(\R^N,\mathbb{M}^{N\times N}(\R))$
and a similar property holds. We have that $u_\gamma$
converges to $u$ in $L^1_{loc}(\R^N)$ and that 
$A_\gamma$ converges to $A$ in $L^1_{loc}(\R^N,\mathbb{M}^{N\times N}(\R))$. In the sequel we make such a convergence much more explicit for %the following two cases: $u\in W^{1,p}(\Omega)$, for some $1\leq p<+\infty$, and 
$A\in\mathcal{M}(\lambda_0,\lambda_1)$ such that $A\in BV(\Omega,\mathbb{M}^{N\times N}(\mathbb{R}))$.

%In the first case, for any $u\in W^{1,p}(\Omega)$, we assume that $u$ is actually defined all over $\R^N$ by identifying $u$ with $Su$, with the extension operator $S$ given in Theorem~\ref{extension}. We have that $u\in W^{1,p}(\R^N)$,
%$u\equiv 0$ outside $B_{R+1}$ and
%$$\|u\|_{W^{1,p}(\R^N)}\leq C_E\|u\|_{W^{1,p}(\Omega)}.$$
%We recall that $C_E$ depends on $r$, $L$, $R$ and $p$ only.
% 
%In the second case, 
By Theorem~\ref{extension}, more precisely by Remark~\ref{extensionwithbound}, we assume that $A$ is actually defined all over $\R^N$ by identifying $A$ with $\tilde{S}A$. So, $A\in L^{\infty}(\R^N,\mathbb{M}_{sym}^{N\times N}(\R))$, satisfies  $\lambda_0 I_N\leq A\leq \lambda_1 I_N$ almost everywhere in $\R^N$,
the total variation of $A$ on $\partial\Omega$ is $0$,
 and
 $$|A|_{BV(\mathbb{R}^N,\mathbb{M})}\leq C_E\|A\|_{BV(\Omega,\mathbb{M})}.$$
We recall that $C_E$ depends on $r$, $L$ and $R$ only. Moreover,
 $$\|A\|_{BV(\Omega,\mathbb{M})}=\|A\|_{L^1(\Omega,\mathbb{M})}+|A|_{BV(\Omega,\mathbb{M})}\leq \lambda_1|\Omega|+|A|_{BV(\Omega,\mathbb{M})}.$$

By \cite[Lemma~3.24]{Amb-Fus-Pal}, the following $L^1$ convergence result holds. For a constant $\tilde{C}_E$ depending on $r$, $L$ and $R$ only, we have,
for any $\gamma$, $0<\gamma\leq 1$,
\begin{equation}\label{firstestima}
\|A-A_{\gamma}\|_{L^1(\Omega,\mathbb{M}^{N\times N}(\R))}\leq \tilde{C}_E\|A\|_{BV(\Omega,\mathbb{M})}\gamma.
\end{equation}
Moreover, by \cite[Proposition~1.15]{Giu}, we conclude that
\begin{equation}\label{strictconvergence0}
|A_{\gamma}|_{BV(\Omega,\mathbb{M})}\leq |A|_{BV(B_{R+1},\mathbb{M})}\quad\text{for any }0<\gamma\leq 1.
\end{equation}
and
\begin{equation}\label{strictconvergence}
\lim_{\gamma\to 0^+}|A_{\gamma}|_{BV(\Omega,\mathbb{M})}=|A|_{BV(\Omega,\mathbb{M})}.
\end{equation}

For any $1\leq q<+\infty$, and any $0<\gamma\leq 1$, we also have
$$\|D^{\alpha}A_{\gamma}\|_{L^q(\Omega)}\leq 
C(|\alpha|)
\|A \|_{L^q(B_{R+1})}
\gamma^{-|\alpha|}\leq C(|\alpha|)|B_1|^{1/q}(R+1)^{N/q}\lambda_1\gamma^{-|\alpha|},$$
where $C(|\alpha|)$ is an absolute constant depending, through our choice of $\eta$, on $|\alpha|$ only.

We conclude that, for some absolute constant $C(N)$, depending, through our choice of $\eta$, on $N$ only,
for any $0<\gamma\leq 1$,
\begin{equation}\label{W2pbound}
\|A_{\gamma}\|_{W^{2,q}(\Omega)}\leq C(N)|B_1|^{1/q}(R+1)^{N/q}\lambda_1\gamma^{-2}.
\end{equation}
This prepares the proof of the following result.

\begin{prop}\label{Adiscrprop}
 Let $A\in\mathcal{M}(\lambda_0,\lambda_1)$ be such that $A\in BV(\Omega,\mathbb{M}^{N\times N}(\mathbb{R}))$. We fix $\alpha$ such that $0<\alpha<1/2$. Then
for any $h$, $0<h\leq 1$, we can find $\tilde{A}_h\in X_h[\mathbb{M}]\cap \mathcal{M}(\lambda_0,\lambda_1)$ such that
\begin{equation}\label{finecotruzione}
\|A-\tilde{A}_h\|_{L^1(\Omega,\mathbb{M}^{N\times N}(\R))}\leq C_Dh^{\alpha}\quad\text{and}\quad
|\tilde{A}_h|_{BV(\Omega,\mathbb{M})}\leq \tilde{C}_D.
\end{equation}
The constants $C_D$ and $\tilde{C}_D$ are defined, respectively, in \eqref{CDdefin} and in \eqref{tildeCDdefin} below.
Moreover,
\begin{equation}\label{finecostruzionebis}
\lim_{h\to 0^+}|\tilde{A}_h|_{BV(\Omega,\mathbb{M})}=|A|_{BV(\Omega,\mathbb{M})}.
\end{equation}
\end{prop}

\proof{.} Fix $q>N/2$. Then by using \eqref{W2pbound} and Theorem~\ref{Ciarletestimteo}, we have that, for any $0< h,\,\gamma\leq 1$,
\begin{equation}\label{Ciarletest2}
\|A_{\gamma}-\Pi_h(A_{\gamma})\|_{L^{q}(\Omega,\mathbb{M}^{N\times N}(\R))}\leq C_1(N)\tilde{C}(q)|B_1|^{1/q}(R+1)^{N/q}\lambda_1h^2\gamma^{-2},
\end{equation}
and, for any $i,\,j=1,\ldots,N$, if $\sigma=A_{ij}$,
\begin{equation}\label{Ciarletest2bis}
\|\nabla \left(\sigma_{\gamma}-\Pi_h(\sigma_{\gamma})\right)\|_{L^q(\Omega,\mathbb{R}^N)}\leq C_1(N)\tilde{C}(q)|B_1|^{1/q}(R+1)^{N/q}\lambda_1sh\gamma^{-2},
\end{equation}
where $C_1(N)$ depends on $N$ only and $\tilde{C}(q)$ is as in \eqref{Ciarletest}. We have that
$\Pi_h(A_{\gamma})\in X_h[\mathbb{M}]$. Furthermore,
\begin{multline*}
\|A-\Pi_h(A_{\gamma})\|_{L^1(\Omega,\mathbb{M}^{N\times N}(\R))}\\\leq
\|A-A_{\gamma}\|_{L^1(\Omega,\mathbb{M}^{N\times N}(\R))}+
\|A_{\gamma}-\Pi_h(A_{\gamma})\|_{L^1(\Omega,\mathbb{M}^{N\times N}(\R))}\\\leq
\tilde{C}_E\|A\|_{BV(\Omega,\mathbb{M})}\gamma
+C_1(N)\tilde{C}(q)|B_1|R(R+1)^{(N-1)/q}\lambda_1h^2\gamma^{-2}.
\end{multline*}
By picking $\gamma=h^{\alpha}$, we call $\tilde{A}_h=\Pi_h(A_{h^{\alpha}})$ and
the first part of \eqref{finecotruzione} is proved, with
\begin{equation}\label{CDdefin}
C_D:=\tilde{C}_E\|A\|_{BV(\Omega)}
+C_1(N)\tilde{C}(q)|B_1|R(R+1)^{(N-1)/q}\lambda_1.
\end{equation}
Finally,
$$|\tilde{A}_h|_{BV(\Omega,\mathbb{M})}\leq |\tilde{A}_h-A_{\gamma}|_{BV(\Omega,\mathbb{M})}+|A_{\gamma}|_{BV(\Omega,\mathbb{M})}$$
hence, by \eqref{Ciarletest2bis} with $\gamma=h^{\alpha}$ and \eqref{strictconvergence0}, for a constant $C(R)$ depending on $R$ only,
$$|\tilde{A}_h|_{BV(\Omega,\mathbb{M})}\leq  C(R)\lambda_1sh^{1-2\alpha}+C_E(\lambda_1|B_1|R^N+|A|_{BV(\Omega,\mathbb{M})})
$$
so the second part of \eqref{finecotruzione} is proved with
\begin{equation}\label{tildeCDdef}
\tilde{C}_D:=C(R)\lambda_1s+C_E(\lambda_1|B_1|R^N+|A|_{BV(\Omega,\mathbb{M})}).
\end{equation}
About \eqref{finecostruzionebis}, we have
\begin{multline*}
\left||\tilde{A}_h|_{BV(\Omega,\mathbb{M})}-|A|_{BV(\Omega,\mathbb{M})}  \right|\\\leq 
\left||\tilde{A}_h|_{BV(\Omega,\mathbb{M})}-|A_\gamma|_{BV(\Omega,\mathbb{M})}\right|+\left||A_\gamma|_{BV(\Omega,\mathbb{M})} -  |A|_{BV(\Omega,\mathbb{M})}  \right|.
\end{multline*}
The first term of the right hand side is controlled by 
$ |\tilde{A}_h-A_{\gamma}|_{BV(\Omega,\mathbb{M})}$, thus goes to $0$ by \eqref{Ciarletest2bis}
with $\gamma=h^{\alpha}$. The second term of the right hand side goes to $0$ by \eqref{strictconvergence}. The proof is concluded.\cvd

\bigskip

We need a similar approximation for solutions to the Neumann boundary value problems.

\begin{prop}\label{discreteSobolev}\label{vdiscrprop}
Let $A\in\mathcal{M}(\lambda_0,\lambda_1)$ be such that $A\in BV(\Omega,\mathbb{M}^{N\times N}(\mathbb{R}))$.
Let $g\in L^2_{\ast}(\partial\Omega)$ and let $v=v[F_g]$ be the solution to \eqref{Neuweak} with $X=H^1(\Omega)$.

Then
for any $h$, $0<h\leq 1$, we can find $\tilde{v}_h\in X_h$ such that
\begin{equation}\label{finecostrsobolev}
\|v-\tilde{v}_h\|_{H^1(\Omega)}\leq \hat{C}_Dh^{\alpha_1}\|g\|_{L^2_{\ast}(\partial\Omega)}.
\end{equation}
The constant $\alpha_1$, $0<\alpha_1<1$ is defined in \eqref{alfa1def}.
The constant $\hat{C}_D$ is defined below in \eqref{tildeCDdefin}.
\end{prop}

\proof{.} By Theorem~\ref{Meyers}, we have that $v\in W^{1,Q_1}(\Omega)$, therefore we can assume, identifying $v$ with the extension $Sv$ of
Theorem~\ref{extension}, that
$v\in W^{1,Q_1}(\R^N)$, $v\equiv 0$ outside $B_{R+1}$ and
$$\|v\|_{W^{1,Q_1}(\R^N)}\leq \hat{C} \|g\|_{L^2_{\ast}(\partial\Omega)}.$$
Let $Q_2=3$ for $N=2$, and, for $N\geq 3$, assuming without loss of generality that $Q_1<N$,  $Q_2=\dfrac{Q_1N}{N-Q_1}>\dfrac{2N}{N-2}\geq 2$.
By Sobolev immersion we obtain that
$$\|v\|_{L^{Q_2}(\R^N)}\leq \hat{C}_1 \|g\|_{L^2_{\ast}(\partial\Omega)}$$
where $\hat{C}$, $\hat{C}_1$ are constants depending on $Q_1$, $r$, $L$, $R$, $\lambda_0$ and $\lambda_1$ only.

The idea is similar to the one in the proof of Proposition~\ref{Adiscrprop}. We choose a suitable $\gamma$, depending on $h$, and we pick
 $\tilde{v}_h=\Pi_h(v_{\gamma})$. Therefore, we need to estimate the $H^1(\Omega)$ norm of $v-v_{\gamma}$ and $v_{\gamma}-\Pi_h(v_{\gamma})$.

The second term can be estimated as follows. Pick $q>\max\{N/2,2\}$, to be chosen later. Then, by Theorem~\ref{Ciarletestimteo},
\begin{multline*}
\|v_{\gamma}-\Pi_h(v_{\gamma})\|_{H^1(\Omega)}\leq C(q)\|v_{\gamma}-\Pi_h(v_{\gamma})\|_{W^{1,q}(\Omega)}\\
\leq C(q)\tilde{C}(q)\max\{s,1\}h\|D^2v_{\gamma}\|_{L^q(\Omega,\mathbb{M})},
\end{multline*}
where $C(q)$ depends on $q$ and $R$ only.

In order to estimate $\|D^2v_{\gamma}\|_{L^q(\Omega,\mathbb{M})}$, we need to distinguish two cases. If $Q_2>N/2$, which holds at least for $N\leq 6$, then we just pick $q=Q_2$ and obtain
$$\|D^2v_{\gamma}\|_{L^q(\Omega,\mathbb{M})}\leq 
C(N)
\|v \|_{L^q(B_{R+1})}
\gamma^{-2}\leq C(N) \hat{C}_1 \|g\|_{L^2_{\ast}(\partial\Omega)}\gamma^{-2},$$
where $C(N)$ is an absolute constant depending, through our choice of $\eta$, on $N$ only. If $N\geq 7$, instead, 
let $q_1=\dfrac{2N}{N+6}>1$. Then let $q$ be such that
$$1+\frac1q=\frac1{Q_2}+\frac1{q_1}.$$
We have that $q>\max\{N/2,2\}$ and depends on $Q_1$ and $N$ only. Then
$$\|D^2v_{\gamma}\|_{L^q(\Omega,\mathbb{M})}\leq 
C_1(N)
\|D^2(\eta_{\gamma}) \|_{L^{q_1}(\R^N,\mathbb{M})}
\hat{C}_1 \|g\|_{L^2_{\ast}(\partial\Omega)},$$
where $C_1(N)$ is an absolute constant depending on $N$ only. We conclude that
$$\|D^2v_{\gamma}\|_{L^q(\Omega,\mathbb{M})}\leq 
C(N)
\hat{C}_1 \|g\|_{L^2_{\ast}(\partial\Omega)}\gamma^{-(N-2)/2},$$
where $C(N)$ is an absolute constant depending, through our choice of $\eta$, on $N$ only. We conclude that, for any $N\geq 2$,
\begin{multline}\label{eq1a}
\|v_{\gamma}-\Pi_h(v_{\gamma})\|_{H^1(\Omega)}\\ 
\leq C(q)\tilde{C}(q)\max\{s,1\}   C(N)
\hat{C}_1 \|g\|_{L^2_{\ast}(\partial\Omega)}h\gamma^{-\max\{(N-2)/2,2\}} 
\end{multline}
where $C(N)$ is an absolute constant depending, through our choice of $\eta$, on $N$ only.

For what concerns $v-v_{\gamma}$, we just know that $v\in W^{1,Q_1}(\R^N)$ with bounded support, for some $Q_1>2$. Hence we can immediately obtain that, for any $0<\gamma\leq 1$,
$$\|\nabla (v-v_{\gamma})\|_{L^{Q_1}(\Omega)}\leq 2\|\nabla v\|_{L^{Q_1}(B_{R+1})}\leq
\hat{C}\|g\|_{L^2_{\ast}(\partial\Omega)}$$
and
$$\|v-v_{\gamma}\|_{L^2(\Omega)}\leq \|v-v_{\gamma}\|_{L^2(B_{R+1})}\leq \|\nabla v\|_{L^2(\R^N,\R^N)}\gamma\leq
\hat{C}\|g\|_{L^2_{\ast}(\partial\Omega)}\gamma
.$$
Here $\hat{C}$ depends on $r$, $L$, $R$, $\lambda_0$ and $\lambda_1$ only.

Estimating the $H^1(\Omega)$ norm of $v-v_{\gamma}$, however, is more delicate. We need to use in an essential way the fact that $v$ is a solution to an elliptic equation. For any $x\in \Omega$, we call $d(x):=\mathrm{dist}(x,\partial\Omega)$.
 We split $\Omega$ into two regions, one close to the boundary and the other far away from the boundary. Namely, for any $0<t\leq 1$, let $\Omega_t:=\{x\in\Omega:\ d(x)>t\}$. Note that, for any $y\in \R^N$ with $\|y\|\leq t$, we have
 $\Omega_t+y\subset\Omega$.
 
 Let us fix $0<t\leq 1$. For any $\varphi\in C^{\infty}_0(\Omega_t)$ and any $y\in \R^N$ with $\|y\|\leq t$, we have that $\tilde{\varphi}(\cdot)=\varphi(\cdot+y)\in C^{\infty}_0(\Omega)$, hence
 \begin{multline*}
 \int_{\Omega_t}A(x-y)\nabla v(x-y)\nabla \varphi(x)dx=
 \int_{\Omega_t-y}A(z)\nabla v(z)\nabla \tilde{\varphi}(z)dz\\=
 \int_{\Omega}A(z)\nabla v(z)\nabla \tilde{\varphi}(z)dz=0.
 \end{multline*}
By density the same property holds for any $\varphi\in H^1_0(\Omega_t)$. Let $0\leq \gamma\leq t$. By multiplying with $\eta_{\gamma}(y)$ and integrating in $dy$, we obtain that
$$\int_{\Omega_t}(A\nabla v)_{\gamma}\nabla\varphi=0\quad\text{for any }\varphi\in H^1_0(\Omega_t).$$
We now consider
$$\int_{\Omega_t}A\nabla(v-v_{\gamma})\cdot\nabla \varphi=\int_{\Omega_t}\big((A\nabla v)_{\gamma}-A\nabla v_{\gamma}   \big)\cdot\nabla\varphi.$$
We note that, for any $x\in\Omega_t$,
$$
(A\nabla v)_{\gamma}(x)-A(x)\nabla v_{\gamma}(x)\\=\int_{B_{\gamma}} \eta_{\gamma}(y)\big(A(x-y)-A(x)\big) \nabla v(x-y)dy.
$$
Hence, picking $q=Q_1$ and $p=\dfrac{2Q_1}{Q_1-2}$, we obtain for any $i,\,j=1,\ldots,N$,
\begin{multline*}
\int_{\Omega_t}\big|(A_{ij}\partial_j v)_{\gamma}(x)-A_{ij}(x)\partial_j v_{\gamma}(x)\big|^2dx
\\\leq\int_{\Omega_t}\left(\int_{B_{\gamma}}\eta_{\gamma}(y)\big|A_{ij}(x-y)-A_{ij}(x)\big|^2 |\partial_j v(x-y)|^2dy\right)dx\leq CD
\end{multline*}
where, using twice the H\"older inequality,
\begin{multline*}C=
\left(\int_{\Omega_t}\left(\int_{B_{\gamma}}\eta_{\gamma}(y)\big|A_{ij}(x-y)-A_{ij}(x)\big|^pdy\right)dx\right)^{2/p}\quad\text{and}\\
D=\left(\int_{\Omega_t}\left( \int_{B_{\gamma}}\eta_{\gamma}(y)|\partial_j v(x-y)|^qdy\right)dx\right)^{2/q}.
\end{multline*}
It is easy to infer that
$$D\leq \|\partial_jv\|_{L^q(\Omega)}^2.$$
 About $C$, we have that it can be bounded from above by
 $$ (2\|A_{ij}\|_{L^{\infty}(\Omega)})^{(2p-2)/p}\left(\int_{\Omega_t}\left(\int_{B_{\gamma}}\eta_{\gamma}(y)\big|A_{ij}(x-y)-A_{ij}(x)\big|dy\right)dx\right)^{2/p},$$
 which, again by \cite[Lemma~3.24]{Amb-Fus-Pal}, 
can be estimated by 
 $$ (2\lambda_1)^{(2p-2)/p}\left(\gamma|DA_{ij}|(\Omega)\right)^{2/p}\leq 
 (2\lambda_1)^{(2p-2)/p}\left(\gamma|A|_{BV(\Omega,\mathbb{M})}\right)^{2/p}.$$
 In conclusion, we obtain that, for a constant $C(Q_1)$ depending on $Q_1$ only,
 $$\|(A\nabla v)_{\gamma}-Av_{\gamma}\|_{L^2(\Omega_t)}\leq C(Q_1)\hat{C}\lambda_1^{\frac{Q_1+2}{2Q_1}}|A|_{BV(\Omega,\mathbb{M})}^{\frac{Q_1-2}{2Q_1}} \|g\|_{L^2_{\ast}(\partial\Omega)}\gamma^{\frac{Q_1-2}{2Q_1}}.
 $$

We fix $\tilde{\chi}\in C^{\infty}(\R)$ such that $\tilde{\chi}$ is increasing, $\tilde{\chi}\equiv 0$ in $(-\infty,0)$ and $\tilde{\chi}\equiv1$ in ($1,+\infty)$. Then we define, for any $0<t_1<t_2\leq 2$,
$$\chi_{t_1,t_2}(x)=\tilde{\chi}\left(\frac{d(x)-t_1}{t_2-t_1}\right) \quad\text{for any }x\in\Omega.$$
We have that $\chi_{t_1,t_2}\equiv 0$ in $B_{t_1}(\partial\Omega)\cap\Omega$ and $\chi_{t_1,t_2}\equiv 1$ in $\Omega_{t_2}$. Moreover,
$\chi_{t_1,t_2}$ is Lipschitz in $\Omega$ and
$$\|\nabla \chi_{t_1,t_2}\|_{L^{\infty}(\Omega,\R^N)}\leq \frac{C_c}{t_2-t_1}$$
with $C_c$ an absolute constant depending on our choice of $\tilde{\chi}$. We consider $0<\gamma\leq t_1$ and choose $\varphi=\chi^2_{t_1,t_2} (v-v_\gamma)$.  We use an argument similar to the one used to prove the Caccioppoli inequality.
Then, calling $\chi=\chi_{t_1,t_2}$ and $w=v-v_{\gamma}$,
$$
\int_{\Omega_{t_1}}\!\!A\nabla w\cdot\nabla \varphi\\=\int_{\Omega_{t_1}}\!\!\chi^2A\nabla w\cdot\nabla w+2\int_{\Omega_{t_1}}\!\!\chi wA\nabla w\cdot\nabla\chi,
$$
therefore
\begin{multline*}
\int_{\Omega_{t_1}}\!\!\chi^2A\nabla w\cdot\nabla w\leq
C(Q_1)\hat{C}\lambda_1^{\frac{Q_1+2}{2Q_1}}|A|_{BV(\Omega,\mathbb{M})}^{\frac{Q_1-2}{2Q_1}} \|g\|_{L^2_{\ast}(\partial\Omega)}\gamma^{\frac{Q_1-2}{2Q_1}}\|\nabla\varphi\|_{L^2(\Omega,\R^N)}
\\-2\int_{\Omega_{t_1}}\!\!\chi wA\nabla w\cdot\nabla\chi=C_1\|\nabla\varphi\|_{L^2(\Omega,\R^N)}+D_1.
\end{multline*}
We have that
$$D_1=-2\int_{\Omega_{t_1}}\!\!\chi wA\nabla w\cdot\nabla\chi\leq
\frac12\int_{\Omega_{t_1}}\!\!\chi^2 A\nabla w\cdot \nabla w+8\int_{\Omega_{t_1}}\!\!w^2 A\nabla \chi\cdot \nabla \chi,
$$
hence
$$
\int_{\Omega_{t_1}}\!\!\chi^2A\nabla w\cdot\nabla w\leq 2C_1\|\nabla\varphi\|_{L^2(\Omega,\R^N)}+16\int_{\Omega_{t_1}}\!\!w^2 A\nabla \chi\cdot \nabla \chi.$$
We have that
$$2C_1\|\nabla\varphi\|_{L^2(\Omega,\R^N)}\leq 3C_1^2+
\frac{4C^2_c}{(t_2-t_1)^2}\int_{\Omega_{t_1}} w^2+ \frac12\int_{\Omega_{t_1}} \chi^2\nabla w\cdot\nabla w$$
and we conclude that
\begin{multline*}
\int_{\Omega_{t_2}}\!\!\|\nabla w\|^2\leq 
\frac1{\lambda_0}\int_{\Omega_{t_2}}\!\!A\nabla w\cdot\nabla w\leq \frac1{\lambda_0}\int_{\Omega_{t_1}}\!\!\chi^2A\nabla w\cdot\nabla w\\\leq
\frac1{\lambda_0}\left(6C_1^2+
\frac{40C^2_c}{(t_2-t_1)^2}\int_{\Omega_{t_1}} w^2\right),
\end{multline*}
that is,
\begin{multline*}
\|\nabla(v-v_{\gamma})\|_{L^2(\Omega_{t_2},\R^N)}\\
\leq \frac8{\lambda_0^{1/2}}\left(C(Q_1)\lambda_1^{\frac{Q_1+2}{2Q_1}}|A|_{BV(\Omega,\mathbb{M})}^{\frac{Q_1-2}{2Q_1}}+\frac{C_c}{t_2-t_1}\right)    \hat{C}\|g\|_{L^2_{\ast}(\partial\Omega)}\gamma^{\frac{Q_1-2}{2Q_1}}.
\end{multline*}
In order to estimate $\nabla ((v-v_{\gamma}))$ near the boundary, we just use the fact that $v-v_{\gamma}\in W^{1,Q_1}(\R^N)$. Therefore,
\begin{multline*}
\|\nabla(v-v_{\gamma})\|_{L^2(\Omega\backslash \Omega_{t_2},\R^N)}\leq  |\Omega\backslash \Omega_{t_2}|^{\frac{Q_1-2}{2Q_1}}\|\nabla(v-v_{\gamma})\|_{L^{Q_1}(\Omega,\R^N)}\\\leq 2
|\Omega\backslash \Omega_{t_2}|^{\frac{Q_1-2}{2Q_1}}\hat{C}\|g\|_{L^2_{\ast}(\partial\Omega)}.
\end{multline*}
We note that there exists a constant $C_b$, depending on $r$, $L$ and $R$ only, such that, for any $0< t\leq 2$ we have
$|\Omega\backslash \Omega_{t}|\leq C_bt$. Finally, we pick, $0<t_1=\gamma\leq 1$ and $t_2=\gamma+\gamma^{\frac{Q_1-2}{3Q_1-2}}\leq
2\gamma^{\frac{Q_1-2}{3Q_1-2}}\leq 2$ and we conclude that, for a constant $C_1(Q_1)$ depending on $Q_1$ only,
\begin{multline}\label{eq1b}
\|v-v_{\gamma}\|_{H^1(\Omega)}\\\leq C_1(Q_1)\left(\frac{\lambda_1^{\frac{Q_1+2}{2Q_1}}}{\lambda_0^{1/2}}|A|_{BV(\Omega,\mathbb{M})}^{\frac{Q_1-2}{2Q_1}}+\frac{1}{\lambda_0^{1/2}}+1\right)  \hat{C}\|g\|_{L^2_{\ast}(\partial\Omega)}\gamma^{\frac{(Q_1-2)^2}{2Q_1(3Q_1-2)}}.
\end{multline}
Now we choose $0<\alpha_0<1$ and, for any $0<h\leq 1$, $\gamma=h^{\alpha_0}$ so that
\begin{equation}\label{alfa1def}
\alpha_1:=\alpha_0\frac{(Q_1-2)^2}{2Q_1(3Q_1-2)}=1-\alpha_0\max\{(N-2)/2,2\}
\end{equation}
and we use \eqref{eq1a} and \eqref{eq1b}. With
the constant 
$\hat{C}_D$ given by
\begin{multline}\label{tildeCDdefin}
\hat{C}_D:=C(q)\tilde{C}(q)\max\{s,1\}   C(N)
\hat{C}_1 
\\+C_1(Q_1)\left(\frac{\lambda_1^{\frac{Q_1+2}{2Q_1}}}{\lambda_0^{1/2}}|A|_{BV(\Omega,\mathbb{M})}^{\frac{Q_1-2}{2Q_1}}+\frac{1}{\lambda_0^{1/2}}+1\right)\hat{C},
\end{multline}
the proof is concluded.\cvd

\bigskip

By C\'ea's Lemma, Theorem~\ref{cea}, we obtain the following estimate for the corresponding Neumann-to-Dirichlet maps.

\begin{cor}\label{Neumanncor}
Let $A\in\mathcal{A}(\lambda_0,\lambda_1)$ be such that $A\in BV(\Omega,\mathbb{M}^{N\times N}(\mathbb{R}))$.

For any $h$, $0<h\leq 1$, let $\mathcal{N}_h(A)=\mathcal{N}_{X^h}(A)$. Then
\begin{equation}\label{corND1}
\|\mathcal{N}(A)-\mathcal{N}_h(A)\|_{\mathcal{L}(L^2_{\ast}(\partial\Omega))}\leq \sqrt{\frac{\lambda_1}{\lambda_0}}C_P\hat{C}_Dh^{\alpha_1}
\end{equation}
and
\begin{equation}\label{corND2}
\|\mathcal{N}(A)-\mathcal{N}_h(A)\|_{\mathcal{L}(L^2_{\ast}(\partial\Omega),H^{1/2}_{\ast}(\partial\Omega) )}\leq \sqrt{\frac{\lambda_1}{\lambda_0}}C_T\hat{C}_Dh^{\alpha_1}.
\end{equation}
As before,
the constant $\alpha_1$, $0<\alpha_1<1$ is defined in \eqref{alfa1def} and
the constant $\hat{C}_D$ is defined in \eqref{tildeCDdefin}.
\end{cor}

\section{The main theorem}\label{mainsec}

In this section we finally summarize all our results, by stating and proving our main approximation theorem, Theorem~\ref{mainteo}.

Let $\Omega\subset\R^N$ such that $\Omega\in\mathcal{A}(r,L,R)$. Let the unknown conductiovity $A_0\in \mathcal{M}(\lambda_0,\lambda_1)$ satisfy
\begin{equation}\label{BVunknown}
|A_0|_{BV(\Omega,\mathbb{M})}<+\infty.
\end{equation}

Let $F_0:\mathcal{M}(\lambda_0,\lambda_1)\to [0,+\infty]$ be defined as follows for any $A\in \mathcal{M}(\lambda_0,\lambda_1)$
$$F_0(A)=\left\{
\begin{array}{ll}|A|_{BV(\Omega,\mathbb{M})}
 &\text{if }\mathcal{N}(A)=\mathcal{N}(A_0)\\
+\infty&\text{otherwise}.\end{array}
\right.$$
By our assumption \eqref{BVunknown}, $F_0$ is not identically equal to $+\infty$.

We call $S=\{A\in \mathcal{M}(\lambda_0,\lambda_1):\ \mathcal{N}(A)=\mathcal{N}(A_0)\}$ the set of solutions to the inverse problems and
$$\hat{S}=\left\{A\in S:\ |A|_{BV(\Omega,\mathbb{M})}=\min_{\tilde{A}\in S}|\tilde{A}|_{BV(\Omega,\mathbb{M})}\right\}$$
the set of optimal solutions. We note that $\hat{S}$ is not empty and compact in $L^1(\Omega,\mathbb{M}^{N\times N}(\R))$, since it coincides with the set of minimizers of $F_0$.

We fix positive constants $s_0$,  $Z_1\leq Z_2$, $\mu_0$, $\theta_0$, $\eta_0$ and we make the following assumptions.

Let $\Omega$ be a discretizable set with constant $s\leq s_0$, as in Definition~\ref{sdiscretizable}.
For any $h$, $0<h\leq 1$, let $\mathcal{T}_h$, $X^h$
and $\Pi_h$ as in Definition~\ref{sdiscretizable}. 
We call $\mathcal{N}_h=\mathcal{N}_{X^h}$.
We call $X^h[\mathbb{M}]$ any $N\times N$ matrix valued function whose elements all belong to $X^h$. The positive parameter $h$ is called the \emph{mesh size parameter} or \emph{mesh parameter}.

Let us consider electrodes $e_m$, constants $z_m$ and extended electrodes $\tilde{e}_m$, $m=1,\ldots,M$, such that \eqref{zl} is satisfied and $\mu\leq \mu_0$, $\mu$ as in \eqref{mudef}. The electrodes are characterized by the positive \emph{electrode size parameter} or \emph{electrode parameter} $\delta$ defined in \eqref{deltadef}. We assume that
$\theta\leq \theta_0$ and $\eta\leq\eta_0$, $\theta$ as in \eqref{thetadef} and $\eta$ as in \eqref{etadef}.

The exact experimental measurements are given by the $M\times M$ \emph{resistance matrix} $R^{\delta}_0=R^{\delta}(A_0)$ associated to the electrodes. We assume that, for some $\varepsilon$, $0<\varepsilon\leq 1$, the \emph{noise level}, the available data are given by the matrix $R_{\varepsilon}^{\delta}$ such that
\begin{multline}\label{noise}
\|R^{\delta}_{\varepsilon}-R_0^{\delta}\|^2\leq\|R^{\delta}_{\varepsilon}-R_0^{\delta}\|^2_{M\times M}=\sum_{i,j=1}^M(R^{\delta}_{\varepsilon}-R^{\delta}_0)_{ij}^2\leq M^2\varepsilon^2
\\\leq 
\left(\tilde{c}_2
\eta_0^{N-1}\theta_0\mu_0\right)^2\delta^{-2(N-1)}\varepsilon^2,
\end{multline}
by \eqref{Lest}, with $\tilde{c}_2$ depending on $r$, $L$ and $R$ only.
We recall that, for all these and the next matrices, we assume that their product with the vector $[1]$ is $0$.

We can define our fully discretized functional to be minimized. Given $0<\varepsilon\leq 1$, $0<h\leq 1$, $\delta>0$ and $a>0$, let 
$F_{\varepsilon,h,a}^{\delta}:\mathcal{M}(\lambda_0,\lambda_1)\to [0,+\infty]$ be defined as follows for any $A\in \mathcal{M}(\lambda_0,\lambda_1)$
$$F_{\varepsilon,h,a}^{\delta}(A)=\left\{
\begin{array}{ll}\dfrac{\|\hat{R}^{\delta}_h(A)-  R_{\varepsilon}^{\delta}\|^2}{a}+|A|_{BV(\Omega,\mathbb{M})}
 &\text{if }A\in X_h[\mathbb{M}]\cap \mathcal{M}(\lambda_0,\lambda_1)\\
+\infty&\text{otherwise}.\end{array}
\right.$$
where the \emph{simplified resistance matrix} $\hat{R}^{\delta}_h(A):=\hat{R}_{X^h}(A)$ is defined in \eqref{simplresmat}.

The regularized and discretized minimization problem one has to solve is $\displaystyle{\min_{A\in\mathcal{M}(\lambda_0,\lambda_1)}F_{\varepsilon,h,a}^{\delta}(A)}$, that is,
\begin{equation}\label{finalminpbm}
\min\left\{ \dfrac{\|\hat{R}^{\delta}_h(A)-  R_{\varepsilon}^{\delta}\|^2}{a}+|A|_{BV(\Omega,\mathbb{M})}:\
A\in X_h[\mathbb{M}]\cap \mathcal{M}(\lambda_0,\lambda_1)  \right\}.
\end{equation}
It is easy to note that \eqref{finalminpbm} admits a solution, that is,
$F_{\varepsilon,h,a}^{\delta}$ admits a (possibly not unique) minimizer that we call $A_{\varepsilon,h,a}^{\delta}$,
with $F_{\varepsilon,h,a}^{\delta}(A_{\varepsilon,h,a}^{\delta})<+\infty$. In fact, $F_{\varepsilon,h,a}^{\delta}\left(\frac{\lambda_0+\lambda_1}2I_N\right)<+\infty$, so a minimizing sequence is bounded in $BV(\Omega,\mathbb{M})$, thus without loss of generality we can assume it converges in $L^1(\Omega,\mathbb{M}^{N\times N}(\R))$. Since  $X_h[\mathbb{M}]\cap \mathcal{M}(\lambda_0,\lambda_1)$ is closed with respect to $L^1$ convergence and $A\mapsto \hat{R}^{\delta}_h(A)$ is continuous with respect to $L^1$ convergence, we obtain existence of the minimum.

In our main theorem we prove that, if we suitably choose the parameters $a$, $h$ and $\delta$ with respect to the noise level $\varepsilon$, we have that, up to subsequences, 
$ A_{\varepsilon,h,a}^{\delta}$ converges, as $\varepsilon\to 0^+$, to $\hat{A}\in \hat{S}$. 

\begin{teo}\label{mainteo}
Under the previous assumptions, for any $\varepsilon$, $0<\varepsilon\leq 1$, let
\begin{equation}\label{hdef-deltadef}
a=a(\varepsilon)=c\varepsilon^{\gamma},\quad h=h(\varepsilon)\leq c_0\varepsilon^{a_1}\quad\text{and}\quad c_1 \varepsilon^{a_2}\leq \delta=\delta(\varepsilon)\leq c_2\varepsilon^{a_2}
\end{equation}
where $c$, $c_0$, $c_1<c_2$, $\gamma$, $a_1$ and $a_2$ are positive constants. We assume that
\begin{equation}\label{coeffdef}
\gamma<a_2;\quad \gamma+2(N-1)a_2<2;\quad \gamma<2a_1\alpha_1;\quad\gamma<\beta_1,
\end{equation}
with $\alpha_1$ defined in \eqref{alfa1def} and $\beta_1$ defined in \eqref{beta1def}.

We call $F_{\varepsilon}=F_{\varepsilon,h(\varepsilon),a(\varepsilon)}^{\delta(\varepsilon)}$.
Then  there exists $\displaystyle{\min_{\mathcal{M}(\lambda_0,\lambda_1)} F_{\varepsilon}}$, for any $\varepsilon$, $0\leq\varepsilon\leq 1$, and
$$\min_{\mathcal{M}(\lambda_0,\lambda_1)} F_0=\lim_{\varepsilon\to 0^+}\min_{\mathcal{M}(\lambda_0,\lambda_1)} F_{\varepsilon}<+\infty.$$

Let $\{A_{\varepsilon}\}_{0<\varepsilon\leq 1}\subset \mathcal{M}(\lambda_0,\lambda_1)$ satisfy $\displaystyle{\limsup_{\varepsilon\to 0^+} F_{\varepsilon}(A_{\varepsilon})<+\infty}$. Then
\begin{equation}
\lim_{\varepsilon\to 0^+}\mathrm{dist}(A_{\varepsilon},S)=0.
\end{equation}
Here and in the sequel, $\mathrm{dist}$ is the distance  with respect to the $L^1(\Omega,\mathbb{M}^{N\times N}(\R))$ norm.
Moreover, for any sequence
$\{\varepsilon_n\}_{n\in\mathbb{N}}$ of numbers such that $0<\varepsilon_n\leq 1$ for any $n\in\N$ and
$\displaystyle{\lim_n\varepsilon_n= 0}$, we have that $\{A_{\varepsilon_n}\}_{n\in\N}$ converges, up to a subsequence, in the $L^1(\Omega,\mathbb{M}^{N\times N}(\R))$ norm to $A\in S$, that is, $A\in\mathcal{M}(\lambda_0,\lambda_1)$
such that $\mathcal{N}(A)=\mathcal{N}(A_0)$.

Let $\{A_{\varepsilon}\}_{0<\varepsilon\leq 1}\subset \mathcal{M}(\lambda_0,\lambda_1)$ satisfy $\displaystyle{\lim_{\varepsilon\to 0^+} \left(F_{\varepsilon}(A_{\varepsilon})-\min_{\mathcal{M}(\lambda_0,\lambda_1)} F_{\varepsilon}\right)=0}$. Then
\begin{equation}
\lim_{\varepsilon\to 0^+}\mathrm{dist}(A_{\varepsilon},\hat{S})=0
\end{equation}
Moreover, for any sequence
$\{\varepsilon_n\}_{n\in\mathbb{N}}$ of numbers such that $0<\varepsilon_n\leq 1$ for any $n\in\N$ and
$\displaystyle{\lim_n\varepsilon_n= 0}$, we have that $\{A_{\varepsilon_n}\}_{n\in\N}$ converges, up to a subsequence, in the $L^1(\Omega,\mathbb{M}^{N\times N}(\R))$ norm to $\hat{A}\in \mathcal{M}(\lambda_0,\lambda_1)$
such that $\hat{A}$ is a minimizer of $F_0$, that is, $\hat{A}\in \hat{S}$.

If $S=\{A_0\}$ and
$\{A_{\varepsilon}\}_{0<\varepsilon\leq 1}\subset \mathcal{M}(\lambda_0,\lambda_1)$ satisfy $\displaystyle{\limsup_{\varepsilon\to 0^+} F_{\varepsilon}(A_{\varepsilon})<+\infty}$, 
then we have that
$$\lim_{\varepsilon\to 0^+}\int_{\Omega}\|A_{\varepsilon}-A_0\|=0.$$
\end{teo}

\proof{.} It follows immediately, through Theorem~\ref{fundthm}, by the following two propositions, Propositions~\ref{Gammaconv} and
\ref{equicoerc}.\cvd

\begin{oss}\label{osserv1}
If we consider the noise level $\varepsilon$ as a relative error, that is, we assume
$$\|R^{\delta}_{\varepsilon}-R_0\|  \leq \varepsilon$$
independently on the number $M$ of electrodes, then the assumptions of Theorer~\ref{mainteo} simplify as follow. We can replace \eqref{hdef-deltadef} and \eqref{coeffdef} with
\begin{equation}\label{hdef-deltadefbis}
a=a(\varepsilon)=c\varepsilon^{\gamma},\quad h=h(\varepsilon)\leq c_0\varepsilon^{a_1}\quad\text{and}\quad 
\delta=\delta(\varepsilon)\leq c_2\varepsilon^{a_2}
\end{equation}
and
\begin{equation}\label{coeffdefbis}
\gamma<a_2;\quad \gamma<2;\quad \gamma<2a_1\alpha_1;\quad\gamma<\beta_1.
\end{equation}
\end{oss}

\begin{oss}
Concerning the choice of $h$ and $\delta$, it is convenient to choose $h$ and $\delta$ as large as possible within the constraints of \eqref{hdef-deltadef} or \eqref{hdef-deltadefbis}. About $h$, this should help to reduce the numerical instability as well as the numerical complexity of the problem to be solved. About $\delta$, this reduces the numbers of electrodes to be used.
\end{oss}

\begin{oss} The main theorem, Theorem~\ref{mainteo}, as well as the previous remark~\ref{osserv1}, holds exactly the same if we replace $\mathcal{M}(\lambda_0,\lambda_1)$ with
$\mathcal{M}_{scal}(\lambda_0,\lambda_1)$. In this case, if $N=2$, we indeed have that $S=\{A_0\}$. 
\end{oss}

We begin by proving the $\Gamma$-convergence result.

\begin{prop}\label{Gammaconv} Under the assumptions and notation of Theorem~\textnormal{\ref{mainteo}}, we have that, as $\varepsilon\to 0^+$,
$F_{\varepsilon}$ $\Gamma$-converges to $F_0$, with respect to the $L^1$ norm on $\mathcal{M}(\lambda_0,\lambda_1)$.
\end{prop}

\proof{.}  
We begin with the $\Gamma$-liminf inequality.
Let us assume that $\{A_{\varepsilon}\}_{\varepsilon\in(0,1]}\subset\mathcal{M}(\lambda_0,\lambda_1)$ is such that
$\displaystyle{\lim_{\varepsilon\to 0^+}A_{\varepsilon}=A}$ in the $L^1$ norm and $\displaystyle{\liminf_{\varepsilon\to 0^+}F_{\varepsilon}(A_{\varepsilon})<+\infty}$. Then
$$|A|_{BV(\Omega,\mathbb{M})}\leq \liminf_{\varepsilon\to 0^+}|A_{\varepsilon}|_{BV(\Omega,\mathbb{M})}   \leq \liminf_{\varepsilon\to 0^+}F_{\varepsilon}(A_{\varepsilon})<+\infty.$$
Hence, it remains to prove that $\mathcal{N}(A)=\mathcal{N}(A_0)$.
Without loss of generality, we can assume that, for a constant $C_0$ we have
\begin{equation}\label{BVuniform}
|A_{\varepsilon}|_{BV(\Omega,\mathbb{M})}\leq C_0\quad\text{for any }0<\varepsilon\leq 1.
\end{equation}
 We have that
\begin{multline*}
\|E\circ ( \hat{\mathcal{R}}^{\delta}_h(A_{\varepsilon})-  \mathcal{R}_{\varepsilon}^{\delta}  )\circ Q\|_{\mathcal{L}(L^2_{\ast}(\partial\Omega))}\leq \theta_0 
\|  \hat{\mathcal{R}}^{\delta}_h(A_{\varepsilon})-  \mathcal{R}_{\varepsilon}^{\delta}  \|_{\mathcal{L}(PC_{\ast})}
\\\leq \theta_0\mu_0^{1/2}\|  \hat{R}^{\delta}_h(A_{\varepsilon})-  R_{\varepsilon}^{\delta}  \|
%\leq \theta_0\mu_0^{1/2} \|\hat{R}^{\delta}_h(A_{\varepsilon})-  R_{\varepsilon}^{\delta}\|^2_{M\times M}
\to 0\quad\text{as }\varepsilon\to 0^+.
\end{multline*}

We next show that
\begin{equation}\label{firstapproxliminf}
\|\mathcal{N}(A_0)-E\circ  \mathcal{R}_{\varepsilon}^{\delta}  \circ Q\|_{\mathcal{L}(L^2_{\ast}(\partial\Omega))}\to 0\quad\text{as }\varepsilon\to 0^+.
\end{equation}
We have that
\begin{multline*}
\|\mathcal{N}(A_0)-E\circ  \mathcal{R}_{\varepsilon}^{\delta}  \circ Q\|_{\mathcal{L}(L^2_{\ast}(\partial\Omega))}\\\leq
\|\mathcal{N}(A_0)-E\circ  \mathcal{R}^{\delta}(A_0)  \circ Q\|_{\mathcal{L}(L^2_{\ast}(\partial\Omega))}\\+\|E\circ(\mathcal{R}^{\delta}(A_0)    -  \mathcal{R}_{\varepsilon}^{\delta}  )\circ Q\|_{\mathcal{L}(L^2_{\ast}(\partial\Omega))}.
\end{multline*}
The first term on the right hand side can be estimated, by using Theorem~\ref{sec4-sub3-theorem1}, as follows
$$\|\mathcal{N}(A_0)-E\circ  \mathcal{R}^{\delta}(A_0)  \circ Q\|_{\mathcal{L}(L^2_{\ast}(\partial\Omega))}\leq \tilde{C}_3\delta(\varepsilon)^{1/2}\to 0\quad\text{as }\varepsilon\to 0^+.$$
The second one is bounded by
\begin{multline*}
\|E\circ(\mathcal{R}^{\delta}(A_0)    -  \mathcal{R}_{\varepsilon}^{\delta}  )\circ Q\|_{\mathcal{L}(L^2_{\ast}(\partial\Omega))}\leq
\theta_0\mu_0^{1/2} \|R^{\delta}(A_0)-  R_{\varepsilon}^{\delta}\|^2_{M\times M}\leq
\theta_0\mu_0^{1/2}M\varepsilon\\\leq \tilde{c}_2
\theta_0^2\mu_0^{3/2}
\eta_0^{N-1}\delta(\varepsilon)^{-(N-1)}\varepsilon\to 0\quad\text{as }\varepsilon\to 0^+.
\end{multline*}
Hence we have
\begin{multline}\label{0estimate}
\|\mathcal{N}(A_0)-E\circ  \mathcal{R}_{\varepsilon}^{\delta}  \circ Q\|_{\mathcal{L}(L^2_{\ast}(\partial\Omega))}\\\leq 
\tilde{C}_3\delta(\varepsilon)^{1/2}+\tilde{c}_2
\theta_0^2\mu_0^{3/2}
\eta_0^{N-1}\delta(\varepsilon)^{-(N-1)}\varepsilon
\end{multline}
and \eqref{firstapproxliminf} is proved.

If we prove that
$$\|E\circ  \hat{\mathcal{R}}^{\delta}_h(A_{\varepsilon})\circ Q-\mathcal{N}(A)\|_{\mathcal{L}(L^2_{\ast}(\partial\Omega))}\to 0\quad\text{as }\varepsilon\to 0^+,$$
the $\Gamma$-liminf inequality is proved. We split the problem into four different terms, namely,
\begin{multline*}
\|E\circ  \hat{\mathcal{R}}^{\delta}_h(A_{\varepsilon})\circ Q-\mathcal{N}(A)\|_{\mathcal{L}(L^2_{\ast}(\partial\Omega))}\leq
\|E\circ (\hat{\mathcal{R}}^{\delta}_h(A_{\varepsilon})-\mathcal{R}^{\delta}_h(A_{\varepsilon}))\circ Q\|_{\mathcal{L}(L^2_{\ast}(\partial\Omega))}
\\
+\|E\circ \mathcal{R}^{\delta}_h(A_{\varepsilon})\circ Q-\mathcal{N}_{h}(A_{\varepsilon})\|_{\mathcal{L}(L^2_{\ast}(\partial\Omega))}+
\|\mathcal{N}_{h}(A_{\varepsilon})-\mathcal{N}(A_{\varepsilon})\|_{\mathcal{L}(L^2_{\ast}(\partial\Omega))}
\\+
\|\mathcal{N}(A_{\varepsilon})-\mathcal{N}(A)\|_{\mathcal{L}(L^2_{\ast}(\partial\Omega))}=(I)+(II)+(III)+(IV).
\end{multline*}
The term $(I)$ is bounded by Proposition~\ref{simplvsnot}. In fact,
\begin{multline}\label{Iestimate}
\|E\circ (\hat{\mathcal{R}}^{\delta}_h(A_{\varepsilon})- \mathcal{R}^{\delta}_h(A_{\varepsilon}))\circ Q\|_{\mathcal{L}(L^2_{\ast}(\partial\Omega))}\\\leq\theta_0\|\hat{\mathcal{R}}^{\delta}_h(A_{\varepsilon})- \mathcal{R}^{\delta}_h(A_{\varepsilon})\|_{\mathcal{L}(PC_{\ast})}\leq \tilde{C}_4\theta_0\delta(\varepsilon)^{1/2}\to 0\quad\text{as }\varepsilon\to 0^+.
\end{multline}
The term $(II)$ is bounded by Theorem~\ref{sec4-sub3-theorem1}. In fact
\begin{equation}\label{IIestimate}
\|E\circ \mathcal{R}^{\delta}_h(A_{\varepsilon})\circ Q-\mathcal{N}_{h}(A_{\varepsilon})\|_{\mathcal{L}(L^2_{\ast}(\partial\Omega))}\leq
\tilde{C}_3\delta(\varepsilon)^{1/2}\to 0\quad\text{as }\varepsilon\to 0^+.
\end{equation}
The term $(III)$ is bounded by Corollary~\ref{Neumanncor}. In fact
\begin{equation}\label{IIIestimate}
\|\mathcal{N}_{h}(A_{\varepsilon})-\mathcal{N}(A_{\varepsilon})\|_{\mathcal{L}(L^2_{\ast}(\partial\Omega))}
\leq
\sqrt{\frac{\lambda_1}{\lambda_0}}C_P\hat{C}_Dh(\varepsilon)^{\alpha_1}\to 0\quad\text{as }\varepsilon\to 0^+.
\end{equation}
Note that $\hat{C}_D$ does not depend on $A_{\varepsilon}$ by \eqref{BVuniform}.
The term $(IV)$ is bounded by Theorem~\ref{pcontinuity}. In fact
\begin{multline}\label{IVestimate}
\|\mathcal{N}(A_{\varepsilon})-\mathcal{N}(A)\|_{\mathcal{L}(L^2_{\ast}(\partial\Omega))}\\\leq
C_T(1+C_P^2)^{1/2}
\tilde{D}_1C(Q_1)(2\lambda_1)^{1-\beta_1}
\|A_{\varepsilon}-A\|^{\beta_1}_{L^1(\Omega,\mathbb{M}^{N\times N}(\R))}\to 0\quad\text{as }\varepsilon\to 0^+.
\end{multline}

The construction of a recovery sequence uses a similar argument. Let $\tilde{A}\in\mathcal{M}(\lambda_0,\lambda_1)$ be such that
$F_0(\tilde{A})<+\infty$, that is, $\mathcal{N}(\tilde{A})=\mathcal{N}(A_0)$ and $|\tilde{A}|_{BV(\Omega,\mathbb{M})}<+\infty$.
Let us apply Proposition~\ref{Adiscrprop} with $A$ replaced by $\tilde{A}$. For any $0<h\leq 1$, we can find $\tilde{A}_h\in X_h[\mathbb{M}\cap\mathcal{M}(\lambda_0,\lambda_1)$ such that
\begin{equation}\label{finecotruzionenuova}
\|\tilde{A}-\tilde{A}_h\|_{L^1(\Omega,\mathbb{M}^{N\times N}(\R))}\leq C_Dh^{\alpha}\quad\text{and}\quad
|\tilde{A}_h|_{BV(\Omega,\mathbb{M})}\leq \tilde{C}_D
\end{equation}
and
\begin{equation}\label{finecotruzionenuovabis}
\lim_{h\to 0^+}|\tilde{A}_h|_{BV(\Omega,\mathbb{M})}=|\tilde{A}|_{BV(\Omega,\mathbb{M})}.
\end{equation}
Here $\alpha$, $0<\alpha<1/2$, is a fixed constant. We wish to show that $A_{\varepsilon}=\tilde{A}_{h(\varepsilon)}$, $0<\varepsilon\leq 1$ is the looked for recovery sequence for $\tilde{A}$. By \eqref{finecotruzionenuovabis}, it is enough to show that
$$\dfrac{\|\hat{R}^{\delta}_h(A_{\varepsilon})-  R_{\varepsilon}^{\delta}\|^2}{\varepsilon^{\gamma}}\to 0\quad\text{as }\varepsilon\to 0^+.$$
But
$$\|\hat{R}^{\delta}_h(A)-  R_{\varepsilon}^{\delta}\|\leq\mu_0^{1/2}
\|E\circ ( \hat{\mathcal{R}}^{\delta}_h(A_{\varepsilon})-  \mathcal{R}_{\varepsilon}^{\delta}  )\circ Q\|_{\mathcal{L}(L^2_{\ast}(\partial\Omega))}.$$
With the same argument as before, using \eqref{0estimate}, \eqref{Iestimate}, \eqref{IIestimate}, \eqref{IIIestimate} and \eqref{IVestimate}, we conclude that
\begin{multline}
\|E\circ ( \hat{\mathcal{R}}^{\delta}_h(A_{\varepsilon})-  \mathcal{R}_{\varepsilon}^{\delta}  )\circ Q\|_{\mathcal{L}(L^2_{\ast}(\partial\Omega))}\\\leq
(2\tilde{C}_3+\tilde{C}_4\theta_0)\delta(\varepsilon)^{1/2}+\tilde{c}_2
\theta_0^2\mu_0^{3/2}
\eta_0^{N-1}\delta(\varepsilon)^{-(N-1)}\varepsilon
\\+
\sqrt{\frac{\lambda_1}{\lambda_0}}C_P\hat{C}_Dh(\varepsilon)^{\alpha_1}+
C_T(1+C_P^2)^{1/2}
\tilde{D}_1C(Q_1)(2\lambda_1)^{1-\beta_1}C_D^{\beta_1}h(\varepsilon)^{\alpha\beta_1}\\\leq
C\left(\delta(\varepsilon)^{1/2}+\delta(\varepsilon)^{-(N-1)}\varepsilon+h(\varepsilon)^{\alpha_1}+h(\varepsilon)^{\alpha\beta_1}
\right).
\end{multline}
with $C$ independent on $\varepsilon$, $0<\varepsilon\leq1$. Then the conclusion follows by choosing $\alpha$ sufficiently close to $1/2$.
.\cvd

\bigskip

Next we prove the equicoerciveness.

\begin{prop}\label{equicoerc}
Under the assumptions and notation of Theorem~\textnormal{\ref{mainteo}}, there exists a set $K\subset \mathcal{M}(\lambda_0,\lambda_1)$, which is compact with respect to the $L^1$ norm, such that for any $0<\varepsilon\leq 1$, we have
$$\inf_K F_{\varepsilon}=\inf_{\mathcal{M}(\lambda_0,\lambda_1)}F_{\varepsilon}.$$
\end{prop}

\proof{.} 
By Proposition~\ref{Gammaconv}, we can find a constant $C$ such that
\begin{equation}\label{equiminim}
\min_{\mathcal{M}(\lambda_0,\lambda_1)}F_{\varepsilon}\leq C\quad\text{for any }0<\varepsilon\leq 1.
\end{equation}

We define $K=\{A_{\varepsilon}\}_{0<\varepsilon\leq 1}$,
where $A_{\varepsilon}$ is a minimizer for $F_{\varepsilon}$, for any $0<\varepsilon\leq 1$.
By \eqref{equiminim}, we obtain that, for some constant $C_1$,
$|A_{\varepsilon}|_{BV(\Omega,\mathbb{M})}\leq C_1$ for any $0<\varepsilon\leq 1$. Then $K$ is relatively compact in 
$L^1(\Omega,\mathbb{M}^{N\times N}(\R))$ by the compact immersion of $BV$ into $L^1$.\cvd

\bigskip

\begin{oss}\label{discreteelectrorem} In this final remark we show that the minimization problem can be made completely discrete. Let
$\{\mathcal{T}_n\}_{n\geq 0}$ be the sequence of triangulations defined in Remark~\ref{oss2}. Let $\tilde{T}_n=(\mathcal{T}_n)|_{\partial\Omega}$ be the corresponding triangulation of the boundary. We recall that, by \eqref{htri}, for some positive constants $C_0$ and $C_1$, we have
$$C_02^{-n}\leq   h_n\leq C_12^{-n}.$$

Let us fix $m\geq 1$ and  $k$, $1\leq k\leq 2^{m(N-1)}$.
Fixed $n\geq m$, for any $K\in\tilde{\mathcal{T}}_{n-m}$, we define $e_K$ as $k$ elements of $\tilde{\mathcal{T}}_{n}$ which are contained in $K$ and
$\tilde{e}_K=K$.
The geometric quantities characterizing the electrodes that we need to check are $\delta_n$, $\mu_n$, $\theta_n$ and $\eta_n$, as in \eqref{deltadef}, \eqref{mudef}, \eqref{thetadef} and \eqref{etadef}, respectively.
We have, for suitable positive constants $\tilde{C}_0$, $\tilde{C}_1$ and $C_2$,
$$\tilde{C}_0h_n\leq \delta_n\leq \tilde{C}_12^{m}h_n$$
and $\mu_n\leq C_2$. We have that
$$\theta_n=\frac{2^{m(N-1)}}{k}$$
and, finally, for another constant $C_3$, $\eta_n\leq C_3$.

We pick $\gamma$, $a_1$ and $a_2$, with $a_1=a_2$, such that \eqref{coeffdef}, or \eqref{coeffdefbis} respectively, is satisfied.
Let $\{\varepsilon_j\}_{j\in \N}$ be a sequence of real numbers such
that $0<\varepsilon_{j+1}<\varepsilon_j\leq 1$ for any $j\in \N$ and $\displaystyle{\lim_j\varepsilon_j=0}$. Then choose subsequences $\{n_k\}_{k\in\N}$ and
$\{\varepsilon_{j_k}\}_{k\in \N}$
such that for any $k\in\N$ we have
$$2^{-(n_k+1)}<\varepsilon_{j_k}^{a_1}\leq 2^{-n_k}$$
so that the sequences
$\{\varepsilon_{j_k}\}_{k\in \N}$, $\{a_k=\varepsilon^{\gamma}_{j_k}\}_{k\in \N}$,
$\{h_{n_k}\}_{k\in\N}$ and $\{\delta_{n_k}\}_{k\in \N}$ satisfy \eqref{hdef-deltadef}, or \eqref{hdef-deltadefbis} respectively.

\end{oss}

\section{Conclusions and perspectives}\label{conclsec}

We have rigorously shown that a good approximate solution to the inverse conductivity problem can be obtained by solving a completely discrete
regularized minimum problem. Moreover, given the noise level, we provide optimal choices for the size of the electrodes to be used to collect the measurements and for the size of the mesh used to discretize the inverse, as well as the direct, problem. Due to the exponential ill-posedness of the problem at hand, one could have conjectured that, to keep the instability under control, 
the size of the mesh should decay logarithmically with respect to the noise level. Instead, we confirm the result already found in \cite{Ron16} that 
the size of the mesh should decay
polynomially with respect to the noise level, a fact that is
absolutely not obvious a priori. Concerning electrodes, the most important finding is that the number of electrodes $M$ should grow also polynomially as the noise level diminishes.

We wish to point out that we are able to deal with a rather severe noise, since we just assume that the measurement at each electrode can be perturbed by a level of noise whose amplitude is independent on the value of the measurement there, thus it is not a relative error. In particular the errors on each electrode can sum up and this could be extremely bad given the fact that the number of electrodes is increasing.

Another interesting consequence of our analysis is that, when dealing with reconstruction, we can use as direct problem the usual Neumann boundary value problem. In particular we do not need to use the experimental measurements direct problem, which is less standard, and, in particular, we 
can completely neglect the value of the contact  impedance at each electrode.

The only restrictions in the current work are due to the geometry of the domain $\Omega$ and the geometry of the electrodes. In fact, we need to assume that $\Omega$ is polyhedral. Moreover, in order to have a completely discretized problem, we need the electrodes to be given by unions of elements of the triangulation of the boundary, as described in Remark~\ref{discreteelectrorem}. In the applications, this is often not the case, both for the domain and the electrodes. Currently we are working to extending the result in these two directions.
For domains with curved boundaries, the issue is to find a suitable approximation with polyhedral domains. We recall that an analysis of this type was carried over in \cite{Ge-Ji-Lu}, but it does not seem enough to achieve the approximation we are interested in. In order to deal with electrodes with more arbitrary shape, another suitable approximation of electrodes with unions of elements of a triangulation of the boundary need to be devised and analyzed. We wish to mention that, in \cite{KLO,KLO2,KLOS,AKLOS}, a very nice line of research, which is strictly related to our plan, has been developed. In fact, in these papers the aim is to tackle, mainly on the numerical side,
modelling errors on the shape of the domain and on the shape of electrodes, as well as on their contact impedance. 

The last step in our analysis, which is however extremely challenging in this general framework, would be to find quantitative convergence estimates of our discrete regularized minimizer to an optimal solution of the Calder\'on problem.

 \appendix
 
 \section{Appendix: geometric properties of Lipschitz domains and related results}\label{geomlemma}
 
 Let $\Omega\in \mathcal{A}(r,L,R)$ for positive constants $r$, $L$ and $R$. We begin by illustrating the following geometric properties of the domain $\Omega$, which depend on $r$, $L$ and $R$ only. This geometric construction is crucial in proving all the results of Subsection~\ref{Sobolevsec}.
 
 For any $x\in \partial \Omega$ we can find positive constants $r_1$ and $r_2$, depending on $r$ and $L$ only,
such that the following holds. Up to a rigid change of coordinates, we assume that $x=0$ and that
\begin{multline*}
\Omega\cap \left([-r_1,r_1]^{N-1}\times [-r_2,r_2]\right)\\=\big\{y=(y_1,\ldots,y_{N-1},y_N)\in [-r_1,r_1]^{N-1}\times [-r_2,r_2]:\\ y_N<\varphi(y_1,\ldots,y_{N-1})\big\}
\end{multline*}
where $\varphi: [-r_1,r_1]^{N-1}\to [-r_2/2,r_2/2]$ is a Lipschitz function with Lipschitz constant bounded by $L$.
Let $S=[-r_1,r_1]^{N-1}\times \R$ and
let $T:S\to S$ be the bi-Lipschitz transformation such that
$T(y',y_N)=(y',y_N-\varphi(y'))$ for any $y'=(y_1,\ldots,y_{N-1})\in [-r_1,r_1]^{N-1}$ and any $y_N\in\R$. We have that
$$T\left(\partial\Omega\cap \left([-r_1,r_1]^{N-1}\times [-r_2,r_2]\right)\right)=\{y\in S:\ y_N=0\}$$ and that there exists $r_3>0$ such that
$$B_{r_3}(x)\subset  U_x=T^{-1}\left((-r_1,r_1)^{N-1}\times (-r_2/2,r_2/2)\right).$$
In particular,
$$B_{r_3}(x)\cap\Omega\subset  U^-_x=T^{-1}\left((-r_1,r_1)^{N-1}\times (-r_2/2,0)\right)\subset \Omega$$
and
$$B_{r_3}(x)\backslash\overline{\Omega}\subset  U^+_x=T^{-1}\left((-r_1,r_1)^{N-1}\times (0,r_2/2)\right)\subset \R^N\backslash\overline{\Omega}.$$
We note that $r_3$ and the Lipschitz constants of $T$ and $T^{-1}$ depend on $r$ and $L$ only.

We prove that there exist $x_1,\ldots,x_m\in \partial\Omega$ such that $\partial\Omega\subset \bigcup_{i=1}^mB_{r_3/8}(x_i)$, with $m\leq m_0$, where $m_0$ is a constant depending on $r_3$ and $R$ only. The argument is the following.
Fix $x_1\in\partial\Omega$ arbitrarily. Then proceed by induction as follows. Given $x_1,\ldots,x_n$, with $n\in \N$, we have two cases. If 
$\partial\Omega  \backslash\bigcup_{i=1}^nB_{r_3/8}(x_i)$ is empty, we have
$\partial\Omega\subset \bigcup_{i=1}^nB_{r_3/8}(x_i)$. Otherwise, let $x_{n+1}\in\partial\Omega \backslash\bigcup_{i=1}^nB_{r_3/8}(x_i)$.
We have that, for any $i\neq j$, $B_{r_3/16}(x_i)\cap B_{r_3/16}(x_j)=\emptyset$. Therefore, after a number of steps that depends on $R$ and $r_3$ only, such a construction has to stop.
It is clear that $\overline{B_{r_3/8}(\partial\Omega)}\subset \bigcup_{i=1}^mB_{r_3/4}(x_i)$.

With a completely analogous argument, we find $x_{m+1},\ldots,x_{m+l}\in \Omega$ such that $B_{r_3/8}(x_i)\subset \Omega$ and
\begin{equation}\label{cover}
K=\overline{B_{r_3/8}(\overline{\Omega})}\subset \left(\bigcup_{i=1}^mB_{r_3/4}(x_i)\right)\cup\left(\bigcup_{i=m+1}^{m+l}B_{r_3/16}(x_i)\right).
\end{equation}
Again, we have that $l\leq l_0$, where $l_0$ is a constant depending on $r_3$ and $R$ only. For convenience, let us call
$r(x_i)=r_3/4$ for any $i=1,\ldots,m$ and $r(x_i)=r_3/16$ for any $i=m+1,\ldots,m+l$.

We construct a partition of unity on $K$ as follows. For any $i=1,\ldots,m+l$, we find
$\psi_i\in C_0^{\infty}(B_{2r(x_i)}(x_i))$ such that
$0\leq \psi_i\leq 1$ in $\mathbb{R}^N$ and
$\sum_{i=1}^{m+l}\psi_i(x)=1$ for any $x\in K$. The construction is classical. Take $\varphi\in C_0^{\infty}(B_{2})$ with $0\leq \varphi\leq 1$ and such that $\varphi\equiv 1$ in $B_1$.
Then we call $\varphi_i(x)=\varphi(  (x-x_i)/r(x_i)  )$, $i=1,\ldots,m+l$, so that $\varphi_i\in C_0^{\infty}(B_{2r(x_i)}(x_i))$, $0\leq \varphi_i\leq 1$ and such that $\varphi_i\equiv 1$ in $B_{r(x_i)}(x_i)$. Then call $\psi_1=\varphi_1$, $\psi_2=\varphi_2(1-\varphi_1)$, $\psi_3=\varphi_3(1-\varphi_1)(1-\varphi_2)$, and so on up to
$\psi_{m+l}=\varphi_{m+l}(1-\varphi_1)\cdots(1-\varphi_{m+l-1})$.
It is not difficult to conclude that, for any $i=1,\ldots,m+l$ and any $k\in \N$, the $C^k$ norm of $\varphi_i$ is bounded by a constant depending
$k$, $r$, $L$ and $R$ only.

 \bigskip
 
 \proof{ of Theorems~\textnormal{\ref{traceineqteo}} and \textnormal{\ref{extension}}.}
 We illustrate the $BV$-case and we leave the analogous Sobolev case to the reader. More details can be found, for instance, in \cite{Leo}, see in particular Chapter~15, where these results are proved except for the precise dependence of the constants involved on the geometric properties of $\Omega$.
 
 Let $u\in BV(\Omega)$.
We use the partition of unity constructed before and define $u_i=u\psi_i$. The $BV$ norm of $u_i$ is controlled by the $BV$ norm of $u$, through constants depending on $r$, $L$ and $R$ only. For any $i=1,\ldots,m$, we apply the transformation $T$, related to $x_i$,
and define $v_i=u_i\circ T^{-1}$ in the parallelepiped $R=(-r_1,r_1)^{N-1}\times (-r_2/2,0)$. We call $S$ its upper side, that is,
$S=(-r_1,r_1)^{N-1}\times \{0\}$. Since $T$ is a bi-Lipschitz transformation, by
\cite[Theorem~3.16]{Amb-Fus-Pal} we infer that
$v_i$ is still a $BV$ function, whose $BV$ norm is controlled by that of $u_i$.

For what concerns the trace estimate, we have that
$v_i|_{S}$ belongs to $L^1(S)$ and, for a constant $C$ depending on $r_1$ and $r_2$ only, we have
$$\|v_i\|_{L^1(S)}\leq C\|v_i\|_{BV(R)}.$$
By going back through $T$, we infer that
$$\|u_i\|_{L^1(\partial\Omega)}\leq C_1\|u\|_{BV(\Omega)}.$$
Since on the boundary $u=\sum_{i=1}^mu_i$ and $m$ is controlled by the geometric constants of $\Omega$, the trace inequality is proved with the right dependence on $r$, $L$ and $R$.

The reverse inequality is quite similar. We start with $\varphi\in L^1(\partial\Omega)$.
Let $\chi_i\in C^{\infty}_0(B_{r_3}(x_i))$ be such that $0\leq \chi_i\leq1$ and $\chi_i\equiv 1$ on $B_{r_3/2}(x_i)$. These may just be equal up to a translation.
We let, for any $i=,1\ldots,m$, by using $T$,
$\tilde{\varphi}_i=(\chi_i\varphi)\circ T^{-1}$ which is an $L^1$ function on $S$ and it is also compactly supported in $S$.
We can find $\hat{v}_i$, a $W^{1,1}$ function on the lower half space $\{y_n<0\}$, whose trace on $\{y_n=0\}$ is $\tilde{\varphi}_i$ and whose $W^{1,1}$ norm is controlled by an absolute constant times the $L^1$ norm of $\tilde{\varphi}_i$. We consider
$u_i=\psi_i(\hat{v}_i\circ T)$. Then $u=\sum_{i=1}^mu_i$ has the desired properties.

For the extension, we use the previous construction and
define, for any $i=1,\ldots,m$, through the usual transformation $T$, $v_i=u\circ T^{-1}$
in the parallelepiped $R$. By an even reflection, we extend $v_i$ to a function $\tilde{v}_i$ on $(-r_1,r_1)^{N-1}\times (-r_2/2,r_2/2)$. Then, we let $\tilde{u}_i=\tilde{v}_i\circ T$, for $i=1,\ldots,m$, and $\tilde{u}_i=u_i$, for $i=m+1,\ldots,m+l$, and obtain the desired extension by summing the functions $\psi_i\tilde{v}_i$. Such a construction provides the desired result in the $BV$ case and in the Sobolev case and it is completely independent from the particular space it is used for.\cvd
 
 \bigskip
 
\begin{oss}\label{extensionwithbound}
Let us assume that $\lambda_0\leq u\leq \lambda_1$ in $\Omega$.
We note that, in the previous construction for the extension, we have that $\tilde{u}_i$ still satisfies 
$\lambda_0\leq \tilde{u}_i\leq \lambda_1$, for any $i=1,\ldots,m+l$.
We note that the union of the domains of $\tilde{u}_i$ contains $K$.
Since the condition is convex, and $\psi_i$ is a partition of unity on $K$, we conclude that the extension $Su$ still satisfies
$\lambda_0\leq Su\leq\lambda_1$ in $K$. Let $\chi:\R\to\R$ be such that $0\leq\chi\leq1$, $\chi\equiv 1$ on $(-\infty,r_3/16]$ and $\chi\equiv 0$ on $[r_3/8,+\infty)$. Then let $\tilde{\chi}=\chi(d)$ where $d$ is the distance from $\overline{\Omega}$. By taking
$$\tilde{S}u:=\tilde{\chi}Su+(1-\tilde{\chi})\frac{\lambda_0+\lambda_1}2$$
we obtain an alternative extension of $u$ to the whole $\R^N$ such that $\lambda_0\leq \tilde{S}u\leq \lambda_1$ in $\R^N$, still
$|D(\tilde{S}u)|(\partial\Omega)=0$, and
$$|D(\tilde{S}u)|(\R^N)\leq C_E\|u\|_{BV(\Omega)}$$
or, for $1\leq p<+\infty$,
$$\|\nabla (\tilde{S}u)\|_{W^{1,p}(\R^N,\R^N)}\leq C_E\|u\|_{W^{1,p}(\Omega)},$$
with $C_E$ with the same dependence as in Theorem~\ref{extension}.

Let $A\in BV(\Omega,\mathbb{M}^{N\times N}(\R))$ be such that $\lambda_0 I_N\leq A\leq \lambda_1 I_N$. Then let $SA$ be
the extension of $A$ obtained by extending all the elements of $A$ with $S$. Then by taking
$$\tilde{S}A:=\tilde{\chi}SA+(1-\tilde{\chi})\frac{\lambda_0+\lambda_1}2I_N$$
we obtain an alternative extension of $A$ to the whole $\R^N$ such that $\lambda_0I_N\leq \tilde{S}A\leq \lambda_1I_N$ in $\R^N$,
the total variation of $\tilde{S}A$ on $\partial\Omega$ is $0$,
 and
$$|\tilde{S}A|_{BV(\R^N,\mathbb{M})}\leq C_E\|A\|_{BV(\Omega,\mathbb{M})}.$$
Finally, if $A$ is symmetric, then both $SA$ and $\tilde{S}A$ are symmetric as well.
\end{oss}

 \proof{ of Theorem~\textnormal{\ref{Pointeo}}.}
 In \cite[Proposition~3.2]{A-Mor-Ros}, 
\eqref{Poin2} and \eqref{Poin3} are proved in the case $p=2$. 
The only difference, for $p$ different from $2$, is in the very last part of the proof of \cite[Lemma~4.1]{A-Mor-Ros}. Following their notation, for some bounded function $w\geq 0$ on $\partial\Omega$ such that
$\int_{\partial \Omega}w=1$, we consider $w(u)=\int_{\partial \Omega}uw$, for any $u\in W^{1,p}(\Omega)$. Then, taking $v=u-w(u)$, it is proved that
$$\|v\|_{L^p(\partial \Omega)}+\|v\|_{L^p(\Omega)}\leq C\| \nabla v\|_{L^p(\Omega,\mathbb{R}^N)}=C\| \nabla u\|_{L^p(\Omega,\mathbb{R}^N)}.
$$
From here, since
\begin{equation}\label{meanest}
\left|\fint_{E}v\right|\leq m(E)^{-1/p}\|v\|_{L^p(E)},
\end{equation}
we easily conclude that
$$\left\|v-\fint_{\partial \Omega}v\right\|_{L^p(\partial \Omega)}+\left\|v-\fint_{\Omega}v\right\|_{L^p(\Omega)}\leq 2C\| \nabla u\|_{L^p(\Omega,\mathbb{R}^N)}.
$$
Since $\displaystyle{u-\fint_{\partial \Omega}u=v-\fint_{\partial \Omega}v}$ and $\displaystyle{u-\fint_{\Omega}u=v-\fint_{\Omega}v}$, \eqref{Poin2} and \eqref{Poin3} are true for any $1\leq p<+\infty$.

In order to obtain \eqref{Poin1},
for $v=u-\displaystyle{\fint_{\Omega}u}$, we have that, by \eqref{Poin2} and \eqref{traceineq} (or \eqref{traceineq2} when $p=1$),
$$\|v\|_{L^p(\partial \Omega)}\leq C_T\|v\|_{W^{1,p}(\Omega)}\leq C_T(1+\tilde{C}_P^p)^{1/p}\| \nabla v\|_{L^p(\Omega,\mathbb{R}^N)}.$$
Therefore,
$$\left\|u-\fint_{\partial \Omega}u\right\|_{L^p(\Omega)}=\left\|v-\fint_{\partial \Omega}v\right\|_{L^p(\Omega)}\leq \|v\|_{L^p(\Omega)}+\left|\fint_{\partial \Omega}v\right|
|\Omega|^{1/p},
$$
hence, using \eqref{meanest},
$$\left\|u-\fint_{\partial \Omega}u\right\|_{L^p(\Omega)}\leq \left(\tilde{C}_P+\left(\frac{|\Omega|}{\mathcal{H}^{N-1}(\partial \Omega)}\right)^{1/p}
C_T(1+\tilde{C}_P^p)^{1/p}\right)\| \nabla v\|_{L^p(\Omega,\mathbb{R}^N)}$$
and \eqref{Poin1} is fully proved with the help of \eqref{measure}.\cvd

 \bigskip
 
 \proof{ of Theorem~\ref{Meyers}.} The case $p=2$ has already been proved, therefore we can assume $p>2$ in what follows. We provide just a hint for the proof, leaving the details to the reader.

From \cite{NMey}, one can deduce the following interior estimate. Let $0<\tilde{r}_1<\tilde{r}_2$ and let $U\in H^1_{loc}(B_{\tilde{r}_2})$ be a weak solution to
$$-\mathrm{div}(A\nabla U)=F\quad\text{in }B_{\tilde{r}_2}$$
for some $F\in (H^1(B_{\tilde{r}_2}))'$. Then there exist a constant $Q>2$, depending on $r_2$, $\lambda_0$ and $\lambda_1$ only, and a constant $C$,
depending on $\tilde{r}_1$, $\tilde{r}_2$, $\lambda_0$, $\lambda_1$ and $p$ only, such that, if $2<p\leq  Q$, then
$$
\|\nabla U\|_{L^p(B_{\tilde{r}_1},\R^N)}\leq C\left(\|F\|_{(W^{1,p}(B_{\tilde{r}_2}))'}+\|U\|_{L^p(B_{\tilde{r}_2})}+\|\nabla U\|_{L^2(B_{\tilde{r}_2},\R^N)}
\right).
$$

This interior estimate is the key point of the proof. By the geometric construction outlined at the beginning of this Appendix, for any $x\in \partial\Omega$, we use the bi-Lipschitz transformation $T$ and transform 
the function $u-\varphi$, in the Dirichlet case, or $v$, in the Neumann case, in a function $U$ satisfying an elliptic equation in
$(-r_1,r_1)^{N-1}\times (-r_2/2,0)$. By a reflection on the upper side of this parallelepiped, reflection which is odd in the Dirichlet case and even in the Neumann one, we obtain a function, which we still call $U$, satisfying an elliptic equation in $(-r_1,r_1)^{N-1}\times (-r_2/2,r_2/2)$. To this function we apply the previous interior estimate. By going back with the tranform $T$, we obtain a suitable estimate on $u-\varphi$ or $v$, respectively, in a neighbourhood of $x$.
If $x$ is in the interior of $\Omega$ we apply directly the interior estimate. By the covering in \eqref{cover}, we need to perform just a finite number of these estimates to control the corresponding norms of $u-\varphi$ or $v$ all over $\Omega$.

Since all parts of this procedure depend only on the geometric constants of $\Omega$, that is, $r$, $L$ and $R$, the proof may be concluded.\cvd

\end{document}